\documentclass[a4paper,reqno]{amsart}
\pdfoutput=1
\usepackage[T1]{fontenc}
\usepackage{lmodern}

\usepackage{amsmath,amsfonts,amssymb,amsthm,stmaryrd,bbm,graphicx,mathtools,enumerate}
\usepackage{mathrsfs}
\usepackage[T1]{fontenc} 
\usepackage[latin1]{inputenc}
\usepackage[english]{babel}
\usepackage[tracking,spacing,kerning,babel]{microtype}
\usepackage{wasysym} 
\usepackage{color}
\usepackage{xcolor}
\usepackage{framed} 

\usepackage{wrapfig} 

\usepackage[normalem]{ulem}

\usepackage[final]{hyperref}   
\hypersetup{
    linktoc=page,
    linkcolor=red,          
    citecolor=blue,        
    filecolor=blue,      
    urlcolor=cyan,
   colorlinks=true           
}

\usepackage{lipsum} 

\usepackage[toc,page]{appendix}

\usepackage{minitoc}
\usepackage{multicol}

\newtheoremstyle{mine}
{\baselineskip}
{\baselineskip}
{\itshape}
{
}
{\bfseries}
{.}
{.5em}
{#1 #2\ifx#3\relax\else~(#3)\fi}

\theoremstyle{mine}

\newtheorem{theorem}{Theorem}
\numberwithin{theorem}{section}
\newtheorem{corollary}[theorem]{Corollary}
\newtheorem{proposition}[theorem]{Proposition}
\newtheorem{lemma}[theorem]{Lemma}

\newtheorem{definition}[theorem]{Definition}

\newtheorem{OP}{Open Problem} 
\numberwithin{equation}{section}

\theoremstyle{remark}
\newtheorem{remark}{Remark}

\colorlet{shadecolor}{blue!10}

\newcommand{\margin}[1]{\textcolor{magenta}{*}\marginpar[\textcolor{magenta} {  \raggedleft  \footnotesize  #1 }  ]{ \textcolor{magenta} { \raggedright  \footnotesize  #1 }  }}

\let\qed=\QED

\renewcommand{\epsilon}{\varepsilon}

\newcommand{\R}{\mathbb{R}}

\newcommand{\Z}{\mathbb{Z}}
\newcommand{\N}{\mathbb{N}}

\def\T{\mathbb{T}}

\def\calC{\mathcal{C}}
\def\calD{\mathcal{D}}

\def\calF{\mathcal{F}}
\def\calG{\mathcal{G}}
\def\calL{\mathcal{L}}
\def\calH{\mathcal{H}}

\def\calN{\mathcal{N}}

\def\diam{\mathrm{diam}}

\def\Var#1{\mathrm{Var}\bigl[ #1\bigr]}

\def\Cov#1{\mathrm{Cov}\bigl[ #1\bigr]}

\def\dist{\mathrm{dist}}

\def\floor#1{\lfloor{#1}\rfloor}

\def\P{\mathbb{P}} 
\def\E{\mathbb{E}} 
\def\md{\mid}

\def \eps {\epsilon}

\def\Bb#1#2{{\def\md{\bigm| }#1\bigl[#2\bigr]}}

\def\Pb{\Bb\P}
\def\Eb{\Bb\E}

\def\FK#1#2#3{{\def\md{\bigm| } \P_{#1}^{\,#2}  \bigl[  #3 \bigr]}}
\def\EFK#1#2#3{{\def\md{\bigm| } \E_{#1}^{\,#2}  \bigl[  #3 \bigr]}}

\def \p {{\partial}}

\def\<#1{\langle #1\rangle}
\newcommand{\red}[1]{{\color{red}#1}}
\newcommand{\blue}[1]{{\color{blue}#1}}

\definecolor{darkgreen}{rgb}{0,0.6,0.05}
\newcommand{\dgreen}[1]{{\color{darkgreen}#1}}

\def\nn{\nonumber}
\def\bi{\begin{itemize}}  
\def\ei{\end{itemize}}
\def\bnum{\begin{enumerate}} 
\def\enum{\end{enumerate}}
\def\ni{\noindent}
\def\bf{\bfseries}

\def\GFF{\mathrm{GFF}}
\def\free{\mathrm{free}}

\def\GFF{\mathrm{GFF}}
\def\IV{\mathrm{IV}}

\def\fracL{u}

\title[Scaling limit of the Discrete Gaussian chain]
{Invisibility of the integers for the discrete Gaussian chain via a Caffarelli-Silvestre  extension of the discrete fractional Laplacian}

\author{Christophe Garban}

\address
{Universit\'e Claude Bernard Lyon 1, CNRS UMR 5208, Institut Camille Jordan, 69622 Villeurbanne, France, and Institut Universitaire de France (IUF)}
\email{garban@math.univ-lyon1.fr}

\begin{document}

\maketitle

\begin{abstract}
The Discrete Gaussian Chain is a model of interfaces $\Psi : \Z \to \Z$ governed by the Hamiltonian
\begin{align*}\label{}
H(\Psi)= \sum_{i\neq j} J_\alpha(|i-j|) |\Psi_i -\Psi_j|^2
\end{align*}
with long-range coupling constants $J_\alpha(k)\asymp k^{-\alpha}$.  For any $\alpha\in [2,3)$ and at high enough temperature, we prove an invariance principle for such an $\alpha$-Discrete Gaussian Chain towards a $H(\alpha)$-fractional Gaussian process where the Hurst index $H$ satisfies $H=H(\alpha)=\frac {\alpha-2} 2$.  
  This result goes beyond a conjecture by Fröhlich and Zegarlinski \cite{frohlich1991phase} which conjectured fluctuations of order $n^{\tfrac 1 2 (\alpha-2) \wedge 1}$ for the Discrete Gaussian Chain.
 \smallskip 
  
 More surprisingly, as opposed to the case of the $2D$ Discrete Gaussian $\Psi : \Z^2 \to \Z$, we prove that the integers do not affect the {\em effective temperature} of the discrete Gaussian Chain at large scales. 
Such an {\em invisibility of the integers} had been predicted by Slurink and Hilhorst  in the special case $\alpha_c=2$ in \cite{slurink1983roughening}. We also identify a similar invisibility of integers when a $2D$ Gaussian Free Field at high temperature is conditioned to take integer values on a dilute enough ``fractal subset'' of $\Z^2$. 
 
 Our proof relies on four main ingredients: 
 \bnum
 \item A Caffareli-Silvestre extension for the discrete fractional Laplacian (which may be of independent interest)
 \item A localisation of the chain in a smoother sub-domain 
 \item A Coulomb gas-type  expansion in the spirit of Fröhlich-Spencer \cite{FS} 
 \item Controlling  the amount of Dirichlet Energy 
 supported by a $1D$ band  for the Green functions of $\Z^2$ Bessel-type random walks 
\enum
\smallskip
Our results also have implications for the so-called {\em Boundary Sine-Gordon field}. Finally, we analyse the (easier) regimes where $\alpha\in(1,2) \cup (3,\infty)$ as well as the $2D$ Discrete Gaussian with long-range coupling constants (for any $\alpha>\alpha_c=4$). 

\end{abstract}

\section{Introduction}

\subsection{Context.}

Motivated by the predictions of Anderson \cite{anderson1970exact}  and Thouless \cite{thouless1969long}, Dyson initiated in \cite{dyson1969non} a celebrated line of research on statistical physics models in $1D$ with long-range interactions (see also \cite{cardy1981one} for a renormalization group viewpoint). 
The most studied case is the Ising model with $1/r^\alpha$ interactions on the $1d$ line $\Z$. Its state space is $\sigma \in \{-1,1\}^\Z$ and its (formal) Hamiltonian is given by 
\begin{align*}\label{}
H(\sigma):= - \sum_{i \neq j \in \Z} \frac {\sigma_i \sigma_j} {|i-j|^\alpha}\,. 
\end{align*}
The model has well-defined infinite volume limits when $\alpha>1$ and it is not difficult to prove that the system is disordered at any temperature when $\alpha>\alpha_c=2$. Dyson proved long-range-order at low temperatures for any $\alpha \in (1,2)$ and for the critical exponent $\alpha_c=2$, Anderson and Thouless had made the striking prediction that not only  long-range-order should hold at low temperature but the phase transition should be discontinuous! The first part of their prediction was first proved in \cite{FS82} while the discontinuity statement was proved in \cite{aizenman1988discontinuity}. See also the recent work \cite{duminil2020long} for a short proof of both facts as well as a review of the literature.
  
\medskip
  
The purpose of this paper is to focus on the $\Z$-valued version of this $1D$ long-range model\footnote{In section \ref{s.visible} we shall also analyze the $2D$ version with same long-range coupling constants.}, which is known as the \textbf{discrete Gaussian chain (DGC)}. Informally its state space is now $\Psi \in \Z^\Z$ and its (formal) Hamiltonian is given by 
\begin{align}\label{e.DG}
H(\Psi):=  \sum_{i \neq j \in \Z} J_\alpha(|i-j|) |\Psi_i - \Psi_j|^2 \,. 
\end{align}
with long-range coupling constants $J_\alpha(k)\asymp k^{-\alpha}$. The most classical choices of coupling constants are 
$J_\alpha(k) = k^{-\alpha}$.
Some of the results below will require a specific choice of coupling constants $J_\alpha(k)\sim c k^{-\alpha}$ in~\eqref{e.J2} and~\eqref{e.Jalpha}.  

We will also consider the following $1D$ {\em sine-Gordon}  model with long-range interactions which interpolates between the discrete Gaussian chain ($\lambda=\infty$) and the Gaussian one ($\lambda=0$): the state space is now made of continuous fields $\phi \in \R^\Z$ with (formal) Hamiltonian  
\begin{align}\label{e.SG}
H(\phi):=  -\lambda \sum_{i \in \Z} \cos(\phi_i) + \sum_{i \neq j \in \Z} J_\alpha(|i-j|) \, |\phi_i - \phi_j|^2 \,. 
\end{align}

%

These two models  (the discrete Gaussian chain and its sine-Gordon version)  have a rich history and arise for example in the following settings:
\bi
\item 
The discrete Gaussian chain~\eqref{e.DG} has been introduced as an effective model of interface between two ordered phases in $2d$ long-range Ising models with same $1/r^\alpha$ decaying coupling constants.
We refer the reader to \cite{kjaer1982discrete, frohlich1991phase, velenik2006localization, coquille2018absence}.

\item The localisation/delocalisation of ~\eqref{e.DG} will be closely related to the analogous localisation/delocalisation of the $2D$ integer-valued GFF (also called {\em Discrete Gaussian model}) $\Psi :\Z^2 \to \Z$. The delocalisation at high temperature of the latter was first proved in the seminal paper \cite{FS} and a beautiful different proof was given in \cite{lammers2022height}. See also the recent works \cite{RonFS, GS1,  GS2, lammers2023delocalisation, bauerschmidt2022discrete,bauerschmidt2022discreteN2,lammers2022dichotomy, park2022central}.

\item The properties of the (Gaussian) harmonic crystal on $\Z^d$ with (transient) long-range interactions has been studied for example in \cite{caputo2000large, caputo2001critical, chiarini2016extremes}. 

\item The 1D sine-Gordon model~\eqref{e.SG} arises in models of quantum tunnelling within dissipative systems such as in 
the so-called {\em Caldeira-Leggett model} \cite{caldeira1983quantum} (see in particular the action 4.27 in \cite{caldeira1983quantum}). It is also related to the so-called {\em Polaron problem}. See 
\cite{feynman1955slow, fisher1985quantum, spohn11985quantum,  spohn1986roughening, spohn1986statistical, spohn1987effective, spohn2005models}.

\item The two models above, ~\eqref{e.DG} and ~\eqref{e.SG}, are closely related to the {\em boundary sine-Gordon model.} See \cite{caux2003two,leclair1995boundary, saleur1995boundary}. 


\item Finally, this model fits naturally in the class of {\em self-attracting random walks}, where the attraction mechanism decays as time passes in $1/T^\alpha$. We refer to the nice lecture notes \cite{bolthausen1999large}. 
\ei

The phase diagram of such non-compact spin systems ($\Z$ and $\R$-valued  instead of $\{-1,1\}$) is very similar to the case of the Ising model, with same $\alpha_c=2$, except the difficulty to analyse each phase turns out to be reversed (the high temperature phase will be more challenging while the low temperature case, at least when $\alpha<\alpha_c=2$ will turn out to be rather soft):

\bi
\item When $\alpha\in (1,2)$, we will show that there is long-range order at any inverse temperature $\beta$, in the sense that the field $\Psi$ will be {\em localized} at any $\beta$. Notice that LRO at any $\beta$ differs here w.r.t to the case of the Ising model. 
This question was asked in \cite{frohlich1991phase} and it turns out one can conclude here by a simple comparison with the Gaussian case thanks to Ginibre inequality \cite{Ginibre}. 

\item When $\alpha=\alpha_c=2$, the situation as in the case of the Ising model is rather delicate and the proofs are significantly more involved in the present non-compact case. There is a localisation/delocalisation phase transition in this case. The high temperature regime (=delocalisation or {\em rough} phase) has been analyzed  in \cite{kjaer1982discrete} as outlined in \cite{frohlich1991phase}. We will give a different proof in this paper which relies instead on \cite{FS}. The advantage of our proof is two-fold: first, it enables us to obtain much more detailed information on the limiting fluctuations, namely an invariance principle in the limit $n\to \infty$. Second, it allows us to identify an intriguing  phenomenon which we call {\em invisibility of integers}: namely we show that at high enough temperature the scaling limit of the discrete Gaussian chain coincides with the scaling limit of the (unconditioned) Gaussian chain and they share the same effective temperature! This is solving a conjecture of Slurink and Hilhorst from \cite{slurink1983roughening}.  

The low temperature case(=localisation or {\em smooth} phase) is difficult and was handled in \cite{frohlich1991phase} by building on the multiscale analysis for the $1D$ long-range Ising case from  \cite{FS82}. 

Similarly as for  the Ising model, one may wonder whether there is a discontinuity between the two phases.  (See Open Problem \ref{op.disc}). 


\item When $\alpha>\alpha_c=2$ as far as we know, nothing is known rigorously. It is conjectured in this case  that the discrete Gaussian chain is delocalized at all temperatures. Our second main result below is a detailed analysis of the fluctuations of the discrete Gaussian chain in the high temperature regime when $\alpha\in(2,3)$, where the phenomenon of {\em invisibility of integers} is also shown to hold,  and in the full temperature regime in the easier case $\alpha>3$. (N.B. We also analyse the second threshold $\alpha=3$ which requires additional  log-corrections, see Proposition \ref{pr.comp1}).
When $2\leq \alpha < 3$, we obtain for suitable choices of coupling constants $J_\alpha(r)$ an invariance principle towards a \textbf{$H$-fractional Brownian motion}, with $H=\frac {\alpha -2}  2$.

\begin{figure}[!htp]
\begin{center}
\includegraphics[width=\textwidth]{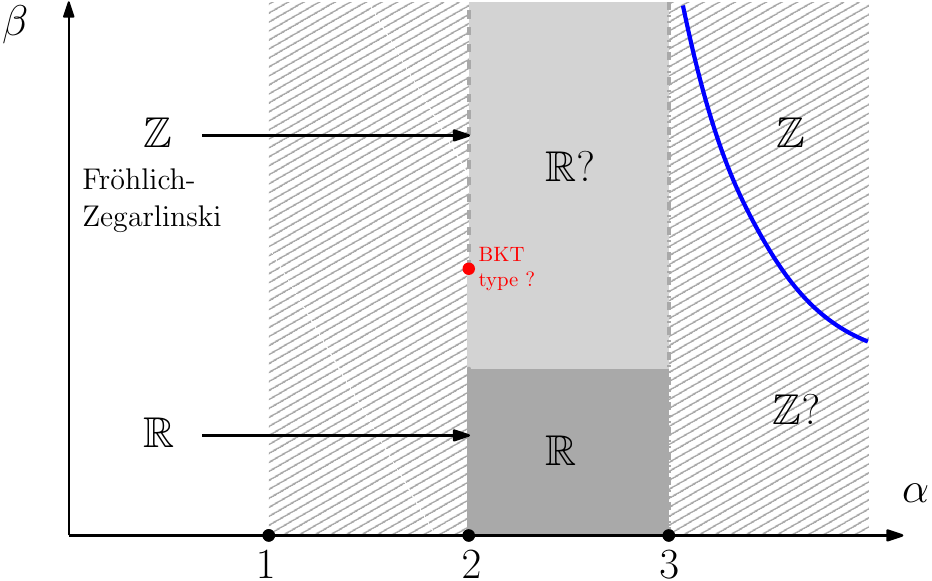}
\end{center}
\caption{The grey shaded region is the expected region with {\em invisibility of the integers} for the Discrete Gaussian Chain. We denote it by $\R$ by contrast with the $\Z$-region where the field conditioned to take integer values is expected to have different large scale fluctuations. The invisibility of integers is proved in this paper for $\beta$ small enough when $\alpha\in [2,3)$ (darkgray region, including a segment along $\{\alpha=2\}$). The red point along the line $\alpha=\alpha_c=2$ should correspond to a BKT-type phase transition. We do not know what to conjecture about the impact of integers when $\alpha=3$ (there are some inevitable $\log$-corrections which appear there as shown in Proposition \ref{pr.comp1}). We expect integers to be visible everywhere outside of the band $\{ 2 \leq \alpha \leq  3\}$. (This is proved very easily when $\alpha>3$ and $\beta$ large enough (blue curve) in Remark \ref{r.NonM}).}\label{f.invisible}
\end{figure}

\item The phenomenon of {\em invisibility of the integers} which we identify here appears to be non-monotone! Indeed, integers are visible at $\alpha=\alpha_c=2$ when $\beta$ is large enough (\cite{frohlich1991phase}) and are invisible when $\beta$ is small enough (Theorem \ref{th.main}). Then they remain invisible when $\alpha\in (2,3)$ and for small enough $\beta$ (also Theorem \ref{th.main}). We expect this should hold for all values of $\beta$ for these $\alpha$ (see Open Problem \ref{OP3}).  And finally they become visible again when $\alpha>3$. See Figure \ref{f.invisible} for the invisibility region (partly conjectural) in the plane $(\alpha, \beta)$ as well as Remark \ref{r.invisible}. 
\ei

\subsection{Main results.}

Consider $\alpha>1$, $\beta>0$ and some fixed coupling constants $J_\alpha=\{J_\alpha(r)\asymp r^{-\alpha}\}_{r\geq 1}$.
 For any $n\geq 1$, let $\Lambda_n$ be the $1D$ interval $\{-n,\ldots,n\}$. The discrete Gaussian chain $\Psi_n$ on $\Lambda_n$, at inverse temperature $\beta$ and with Dirichlet boundary conditions  is defined as follows:
\begin{align*}\label{}
\FK{\Lambda_n, J_\alpha, \beta}{dGC,0}{\Psi_n=h} \propto
\exp\left(-\frac \beta 2  \sum_{i\neq j \in \Z}  J_\alpha(|i-j|)\, |h_i-h_j|^2 \right) 
1_{h_i \in \Z, \forall i }  \, 1_{h_i=0 \text{ on } \Lambda_n^c} \,.
\end{align*}
(N.B. The superscript $0$ stands for Dirichlet boundary conditions). 
\smallskip

Our main result is the invariance principle below. (See Theorem \ref{th.main2} for a more precise version). 

\begin{theorem}\label{th.main}
For any $\alpha\in[2,3)$, there exist $\alpha$-long-range coupling constants $J_\alpha(r) \sim c r^{-\alpha}$ (as defined in~\eqref{e.J2},~\eqref{e.Jalpha}), constants $K>0$ and $\beta^*$ (which may depend on $J_\alpha$) such that the following holds.  If the temperature is high enough, $\beta\leq \beta^*$, then under a suitable rescaling\footnote{See Subsection \ref{ss.TheProof}}, the field $\Psi_n\sim \P_{\Lambda_n, J_\alpha, \beta}^{dGC,0}$,  converges in law towards the fractional Brownian motion
\begin{align*}\label{}
\frac K {\sqrt{\beta}} B_t^H 
\end{align*}
on $(-1,1)$ with Dirichlet boundary conditions  (see Definition \ref{d.fBm3} and \cite{lodhia2016fractional}), and  with Hurst index 
\begin{align*}\label{}
H(\alpha)= \frac{\alpha-2} 2 \,.
\end{align*}

Furthermore, for any choice of coupling constants $J_\alpha(r)\asymp r^{-\alpha}$, we have the following estimates on the fluctuations on the field $\Psi_n$ at low enough $\beta$:
\begin{align}\label{e.UpToC}
\begin{cases}
& \mathrm{Var}_{\Lambda_n, J_\alpha, \beta}^{dGC,0}[\Psi_n(0)] \asymp \log n,  \,\,\,\text{   if } \alpha=\alpha_2=2 \\
& $ $ \\ 
& \mathrm{Var}_{\Lambda_n, J_\alpha, \beta}^{dGC,0}[\Psi_n(0)] \asymp  n^{\alpha-2}, \,\,\,  \text{   if } \alpha\in (2,3) \,.
\end{cases}
\end{align}
\end{theorem}

\ni
In particular, the scaling limit of this integer-valued  field $\Psi_n$ field is exactly the same, including its scaling in $\beta$,  as for the (unconditional) Gaussian chain $\varphi_n$ defined in~\eqref{e.GC} and for which an invariance principle is stated in Proposition \ref{pr.sheffield}. This means the large fluctuations of the limiting field do not feel the $\Z$-conditioning. We call this intriguing phenomenon \uline{invisibility of the integers}. In the physics terminology, this corresponds to $\beta_{eff}(\beta) = \beta$ where $\beta_{eff}(\beta)$ stands for the {\em effective inverse temperature} of the model at $\beta$. (We refer to \cite{bauerschmidt2022discrete,bauerschmidt2022discreteN2,GS2} for a definition of the concept of effective inverse temperature $\beta_{eff}$). 

The above Theorem thus goes beyond a conjecture by Fröhlich and Zegarlinski \cite{frohlich1991phase} and calls for some remarks:
\bnum
\item The fact the inverse temperature $\beta_{eff}\equiv \beta$ in our present setting  is very different from what happens with the $2D$ discrete Gaussian  $\Psi : \Z^2 \to \Z$. Indeed in this latter case, it has been shown recently in the breakthrough paper \cite{bauerschmidt2022discrete} that the discrete Gaussian chain on the $2D$ torus converges to the  Gaussian free field distribution at high enough temperature with a non-explicit effective inverse temperature $\beta \mapsto \beta_{eff}(\beta)$. It is a consequence of \cite{GS2} that $\beta_{eff}(\beta) > \beta$, with quantitative bounds of the type 
\begin{align*}\label{}
\beta_{eff} > \beta + \exp(- c \frac 1 {\beta})
\end{align*}
(see \cite{GS2}). This demonstrates that a GFF on $\Lambda_n^2=\{-n,\ldots,n\}^2$ conditioned to take its values in the integers will fluctuate \underline{strictly} less at macroscopic scales than its unconditioned Gaussian version. 

Our main result thus shows that this is not the case for the discrete Gaussian chain when $\alpha\in [2,3)$ at high enough temperature.  

\item The case $\alpha_c=2$ is rather subtle as by \cite{frohlich1991phase}, it is known that $\beta_{eff}(\beta)=\beta$ cannot hold up to the low temperature phase. Indeed it is proved there that the discrete Gaussian chain at $\alpha_c=2$ is localized at high $\beta$ (which implies in particular that $\beta_{eff}(\beta)=0$ for those $\beta$).

In fact the identity $\beta_{eff}(\beta)=\beta$ at low $\beta$ had been anticipated by Slurink and Hilhorst in the special case $\alpha_c=2$ in \cite{slurink1983roughening} based on numerical simulations.
See also the influential RG analysis by Cardy in \cite{cardy1981one} which identifies a similar phenomenon for long-range Ising model with $\alpha_c=2$. 

As such our results strongly suggest a discontinuous phase transition for the discrete Gaussian chain at $\alpha_c=2$ similarly as in the case of the long-range Ising model at same $\alpha_c=2$.  See the Open Problems. 

\item In the recent work \cite{biskup2023limit}, Biskup and Huang proved an analogous  {\em invisibility of integers} phenomenon which holds for a hierarchical version of the present discrete Gaussian chain. The hierarchical structure allows them to rely on a recursive structure in order to exhibit a robust coupling between the hierarchical integer-valued field and its unconditioned Gaussian version.

Interestingly, it can be easily shown that  
the $\alpha_c=2$ hierarchical discrete Gaussian chain does not see  the peculiar phase transition of localisation/delocalisation of the non-hierarchical one. We prove this in the short Section \ref{s.9}. 


\item We obtain our invariance principle for a particular form of the $\alpha$-long-range coupling constants $J_\alpha$ defined in~\eqref{e.J2} and~\eqref{e.Jalpha}. We expect it should hold for any choice of $J_\alpha(r)$ satisfying $J_\alpha \sim c r^{-\alpha}$.  

\enum

The proof of the above Theorem in fact implies  the following Corollary for the so-called long-range $1D$ {\em sine-Gordon model} defined in~\eqref{e.SG}. (This is obtained essentially by applying Ginibre inequality via the interpolation used for example in \cite{RonFS}).

\begin{corollary}\label{c.BoundarySineGordon}
Let $\alpha\in [2,3)$ and $J_\alpha(r)$ as defined in~\eqref{e.J2} and~\eqref{e.Jalpha}.  If $\beta$ is small enough, then uniformly in the coupling constant $\lambda>0$, the $1D$ $\lambda$-Sine-Gordon model with long-range interactions $J_\alpha(r)$ and with Dirichlet boundary conditions converges after suitable rescaling to the $H(\alpha)-$fractional Brownian motion with Dirichlet boundary conditions on $(-1,1)$. 
\end{corollary} 

The case $\alpha=\alpha_c=2$ in the above theorem, despite being formulated in $1D$, turns out to correspond to the so-called {\em boundary sine-Gordon model} in $2D$ (see for example \cite{caux2003two,leclair1995boundary, saleur1995boundary}). This fact will in fact play a major role in this paper.

 Dirichlet-boundary conditions are rather natural in the $1D$ long-range case (and also they correspond to very explicit coupling constants at $\alpha_c=2$ defined in~\eqref{e.J2}).  See Figure \ref{f.diamond} which illustrates those boundary conditions.
  Yet, in the case of the {\em boundary Sine-gordon model} which is intrinsically of $2D$ nature, another boundary condition which we call the {\em box-boundary condition} is even more natural. (See also Remark \ref{r.BC} for other suitable boundary conditions). 

This allows us to state the result below on the effect of conditioning a $2D$ standard GFF to take its values in $\Z$ on a $1D$ line. 
Let us then consider a standard GFF $\phi$ on $\Lambda_n^2:=\{-n,\ldots,n\}^2 \subset \Z^2$ with Dirichlet boundary conditions along $\p \Lambda_n^2$, i.e. 
\begin{align*}\label{}
\frac{d\mu(\phi)} {d \phi}  \propto \exp\left(-\frac \beta 2 \sum_{x\neq y \in \Z^2} (\phi_x-\phi_y)^2\right) 1_{\phi_z=0, \forall z \in \Z^2 \setminus \Lambda_n^2}\,.
\end{align*}
Now, let $\Psi^{line}=\Psi^{line}_n$  be the Gaussian free field $\phi$ conditioned to take its values in $\Z$ on the $1D$ line $\{-n,\ldots,n\} \times \{0\}$. See Figure \ref{f.diamond}. In the same fashion, one may consider $\varphi$ a GFF on $H_n:=\{-n,\ldots,n\}\times \{0,\ldots,n\}$ equipped with Dirichlet boundary conditions on top/right/left boundaries and free boundary conditions on the bottom boundary $\{-n,\ldots,n\}\times \{0\}$.
In this case, we call $\Psi^{boundary}=\Psi^{boundary}_n$ the GFF $\varphi$ conditioned to take integer values on the bottom boundary $\{-n,\ldots,n\}\times \{0\}$. 

In the Theorem below (see Theorems \ref{th.line2} and \ref{th.FP} for more precise statements), we rescale these two conditioned field $\Psi^{line}_n$ and $\Psi^{boundary}_n$ in order to view them respectively on 
\begin{align*}\label{}
\frac 1 n \Z^2 \cap [-1,1]^2 \text{    and   } \frac 1 n \Z^2 \cap [-1,1]\times [0,1]\,.
\end{align*}

\begin{theorem}\label{th.line}
For  $\beta$ small enough, the rescaled fields $\Psi_n^{line}$ (resp. $\Psi_n^{boundary}$) converge as $n\to \infty$ in the sense of distributions to a GFF on $[-1,1]^2$ with $0$-boundary conditions (resp. GFF on $[-1,1]\times [0,1]$ with free boundary conditions along the bottom side) and with \underline{same}\footnote{This is another instance of the fact that $\beta_{eff}(\beta)=\beta$.} inverse temperature $\beta$.

In other words, the $\Z$-conditioning on a $1D$ line is invisible at large scales when $\beta$ is small enough.
\end{theorem}

\begin{remark}\label{}
In Theorem \ref{th.fractal}, we will prove a {\em invisibility of integers} property for  more general sets on which one requires the GFF field to be $\Z$-valued, such as an horizontal strip 
$\Lambda_n \times [0,H_n]$ or also a {\em Monge-type} fractal set  $F_n\subset \Lambda_n^2$ as pictured in Figure \ref{f.fractal}. 
\end{remark}
\begin{remark}
The localisation result from \cite{frohlich1991phase} shows that when $\beta$ becomes sufficiently large, then suddenly the line becomes visible. It suggests that  the limit of $\Psi_n^{line}$ should then become a GFF with extended Dirichlet boundary conditions, namely 
\begin{align*}\label{}
\p [-1,1]^2 \bigcup [-1,1]\times \{0\}\,.
\end{align*}
\end{remark}


When $\alpha\in(1,2)$, though the question is asked explicitly in \cite{frohlich1991phase}, it turns out that this regime can easily be analysed thanks to Ginibre inequality.   

When $\alpha>3$, we obtain a soft proof of delocalisation at all inverse temperatures $\beta>0$ using a suitable comparison of quadratic forms.

We combine these two regimes in the following  Proposition. 

\begin{proposition}\label{pr.comp1}
The regimes $\alpha\in (1,2)$ and $\alpha\in(3,\infty)$ do not undergo a phase transition in $\beta$:
\bi
\item If $\alpha\in (1,2)$, the discrete Gaussian chain is localised at all inverse temperature $\beta$: there exists $c=c(J_\alpha)>0$ such that for any $\beta>0$,  $\Var{\Psi_n} \leq \frac c \beta$.
\item If $\alpha>3$, then for any $\beta>0$, there exist $0<c_1(\beta)<c_2(\beta)$ s.t. 
 \begin{align*}\label{} 
c_1(\beta)  n  \leq \Var{\Psi_n(x=0)} \leq c_2(\beta) n\,,\,\,\, \text{  for all  } n\geq 1\,.
\end{align*}
\item Furthermore, when $\alpha=3$, we obtain using an estimate from \cite{caputo2009recurrence} the following upper bound. For any $\beta>0$,
\begin{align*}\label{} 
\Var{\Psi_n(x=0)} \leq \frac c \beta \frac {n} {\sqrt{\log(n)}}\,,\,\,\, \text{  for all  } n\geq 1\,.
\end{align*}
While for small enough $\beta$, we obtain a constant $\tilde c>0$ so that for any $\beta\leq \beta*$, 
\begin{align*}\label{}
\Var{\Psi_n(x=0)} \geq \frac {\tilde c} \beta \frac {n} {\log(n)}
\end{align*}
This shows that $\alpha=3$ is somewhat special. 
\ei
\end{proposition}

\begin{remark}\label{}
We expect that the regime $\alpha\in(2,3)$ should not undergo a phase transition in $\beta$ either. 
Our current proof which is based on the Coulomb gas expansion from \cite{FS} would not work at such low temperatures.
See Open Problems \ref{OP3} and \ref{OP5}. 
\end{remark}

\begin{remark}\label{r.invisible}
As illustrated in Figure \ref{f.invisible}, interestingly, the behaviour of $\frac {\beta_{eff}(\beta)} {\beta}$ appears to be "non-monotone" in $\alpha$, more precisely 
\bi
\item At $\alpha_c=2$, it follows from \cite{frohlich1991phase} together with either \cite{kjaer1982discrete} or our present results that $\beta_{eff}(\beta)/\beta$ is $\beta$-dependent. 
\item Then, for $\alpha\in(2,3)$, we conjecture that $\beta_{eff}(\beta)=\beta$ for all $\beta$. (See Theorem \ref{th.main} for a proof at high temperature and Open Problem \ref{OP3}).
\item Finally, when $\alpha>3$, the effective inverse temperature again depends non-trivially on $\beta$. See Remark \ref{r.NonM} in Section \ref{s.visible}. 
\ei
\end{remark}

Finally, using similar considerations as in the above Proposition \ref{pr.comp1}, we obtain the following delocalisation result for the IV-GFF on $\Z^2$ with long-range interactions. 
\begin{proposition}\label{pr.comp2}
Consider the Integer-valued GFF on the 2d grid $\Z^2$ with $1/r^\alpha$ long-range interactions.  Then for all $\alpha>\alpha_c(d=2)=4$, one has delocalisation at high enough temperature. I.e. there exists $c>0$, such that for all $\beta$ small enough, if $\Psi_n$ is the $\alpha$-long-range IV-GFF conditioned to be zero outside the square $\Lambda^2_n=\{-n,\ldots,n\}^2$, then 
\begin{align*}\label{}
\Var{\Psi_n(x=0)} \geq \frac {c} \beta   \log n\,.
\end{align*}
When $\alpha=\alpha_c=4$, we only obtain the following upper bound 
\begin{align*}\label{}
\Var{\Psi_n(x=0)} \leq \frac{O(1)} \beta   \log \log n\,.
\end{align*}
\end{proposition}
We do not know if there is a localisation/delocalisation phase transition at $\alpha_c=4$ and leave it as an interesting open problem. See Open Problem \ref{OP1}.


\subsection{Idea of proof of our main Theorem \ref{th.main}.} 

There exist by now three very different ways of proving the delocalisation of the so-called {\bf discrete Gaussian} $\Psi : \Z^2 \to \Z$. 
\bnum
\item The first proof goes back to the seminal work \cite{FS} by Fröhlich and Spencer. Their proof relies on a Coulomb-gas decomposition. It is then tempting to apply their strategy to the present setting. Instead of working with the standard (nearest-neighbour) Laplacian $\Delta_{\Z^2}$ on $\Z^2$, we would need to work with the discrete fractional Laplacian
$(-\Delta_\Z)^{\fracL}$ on $\Z$ defined for any $\fracL\in(0,1)$ by 
\begin{align}\label{e.FL}
(-\Delta_\Z)^{\fracL} f(x) := -\sum_{y\neq x} \frac{f(y)-f(x)}{|y-x|^{1+2\fracL}}\,.
\end{align}
(see Subsection~\eqref{ss.DFL} for background on the discrete fractional Laplacian and its different possible forms).
One has then the correspondance $\fracL=\frac {\alpha -1} 2$. 
As we shall see below, the approach of \cite{FS} induces major difficulties in the present case due to the fact that the underlying graph is of infinite degree. 

\item A very different approach has been developed  in the recent breakthrough work by Lammers \cite{lammers2022height}. This approach builds on the theory developed by Sheffield in \cite{sheffield2005random}. It appears challenging to extend this framework to non-local Dirichlet forms such as~\eqref{e.FL}. Let us mention the recent relevant work  \cite{lammers2020macroscopic} which studies long-range models of $\R$-valued fields on $\Z^d$. 

Besides the above difficulty, another issue for proving Theorem \ref{th.main} with a proof based on the work \cite{lammers2022height} would be the lack of a quantitative control on the fluctuations. (See also \cite{lammers2022dichotomy} which obtains quantitative fluctuations via RSW type techniques). 


\item Finally, the breakthrough works \cite{bauerschmidt2022discrete,bauerschmidt2022discreteN2} by Bauerschmidt, Park and Rodriguez proved a very non-trivial invariance principle of the $2D$  discrete Gaussian towards a Gaussian free field. Their work relies on a rigorous and tedious Renormalisation group argument.  The $\log n$ delocalisation with same effective temperature is proved in the subsequent work by Park \cite{park2022central}. 

This may also be a  promising approach as rigorous RG flow treatment of spin-systems with long-range interactions have been analyzed for example in \cite{slade2018critical}. 

\enum


In this paper, we rely on the techniques introduced in \cite{FS}. 
But we face here an important difficulty: in the work \cite{FS}, it is of crucial importance that the graph underlying the interactions between Coulomb charges is of \underline{bounded degree}. This assumption is used at two key different places: first in the combinatorial part which decomposes the partition function into a convex combination of Coulomb-type gases. And then at a later stage, in the analytical step which is superposing {\em spin waves} to turn the signed measure into a probability measure (see \cite{FS,RonFS}). 
This is obviously not the case with the fractional Laplacian defined above in~\eqref{e.FL}. 

The need of replacing a non-local operator (-$\Delta_\Z)^\fracL$ by a local elliptic operator (to the cost of  working in dimension $d+1$, i.e. $1+1$ here,  and  loosing some isotropy) is an idea which has been popularized  in the very influential paper by Caffarelli-Silvestre \cite{caffarelli2007extension}. 
We will follow a similar path in the present work by introducing what we shall call a Caffarelli-Silvestre extension of the discrete fractional Laplacian $(-\Delta_\Z)^\fracL$. (See also \cite{slade2018critical} where another useful decomposition is used in the transient case). More accurately, we will introduce a Caffarelli-Silvestre extension of the discrete fractional Laplacian 
\begin{align*}\label{}
(-\Delta_\Z)^{J_\alpha} f(x) := - \sum_{x\neq y} J_\alpha(|x-y|) (f(y)- f(x))\,,
\end{align*}
where the coupling constants $J_\alpha(r)\sim c r^{-\alpha}$ are defined in~\eqref{e.J2} and~\eqref{e.Jalpha}. 


What lead us to such an extension is the special case $\alpha_c=2$ where one may obtain the corresponding fractional Laplacian $\Delta_\Z^{1/2}$ as the restrictions on $1D$ lines of the simple random walk on $\Z^2$. See for example \cite{amir2016one}. In the case of the $\pi/4$ rotated lattice, the coupling constants $J_2(r)$ have a simple explicit form~\eqref{e.J2}. 

Given this context, our proof of Theorem \ref{th.main} is divided into the following steps:
\bnum
\item First, we provide a {\em Caffarelli-Silvestre} extension of fractional discrete Laplacians $(-\Delta_\Z)^\fracL$. In the continuum, such extensions can be understood probabilistically by running a Bessel process in the transverse ``$+1$'' direction. In the case $\alpha=2$, this correponds to the $2d$ Brownian motion and was already observed by Spitzer in \cite{spitzer1958some}. 
Inspired by the continuous setting, we will thus introduce ``Bessel random walks'' on $\Z^2$ whose vertical coordinate follow a $1D$ Bessel walk. To make the analysis slightly simpler, we will also work on a $\pi/4$ rotated lattice $\Z^2$ which we call the {\em diamond graph}: $\calD=  e^{i \pi/4} 2^{-1/2} \Z^2$. For each $\alpha\geq 2$, we will  define 2D $\alpha$-Bessel walks on $\calD$ in Section \ref{s.Bessel}. In turn, these $\alpha$-Bessel walks induce coupling constants $J_\alpha(r)$ (see equation~\eqref{e.Jalpha}). Using \cite{alexander2011excursions}, we show that these coupling constants satisfy the desired asymptotics
\begin{align*}\label{}
J_\alpha(r)\sim \frac c {r^\alpha}\,.
\end{align*}

\item This discrete Caffarelli-Silvestre extension allows us to define a $2D$  Gaussian free field on the graph $\Z^2$ (or $\calD$) in an {\em inhomogeneous} field of (deterministic) conductances. This GFF has the property that its restriction to the line $\Z$ is exactly the Gaussian chain (without the restriction to belong to $\Z$). 

We then apply the analysis of Fröhlich-Spencer in this setting. This leads us to assign Coulomb charges to the vertices of the $1D$ line $\Z$ while letting the rest of the graph $\Z^2 \setminus \Z$ free of charges. 

Following \cite{FS}, we use the power of their multi-scale hierarchy of ``spin-waves''. 

\item The case of Dirichlet boundary conditions corresponds to the geometry of a slit domain (See Figure \ref{f.diamond}). It turns out that the estimates needed for the proof degenerate too much close to the tip of the slit domain. To overcome this issue, we ``localise'' the chain in a sub-domain which has smoother singularities along its boundary (see Figure \ref{f.Hn}). This is handled in Section \ref{s.proof}.

\item From then on, a key point is to notice that since ``vortices'' are bound to a $1D$ line, they have ``less room''  to create fluctuations in the Coulomb gas (indeed as it is shown by the analysis in \cite{GS2}, those are responsible of the non-trivial effective temperature). In the Coulomb-type expansion {\em à la} Fröhlich-Spencer, this corresponds to proving that most of the Dirichlet energy of the Green function associated to these Bessel walks is not confined in the line, but rather spread over the rest of $\Z^2$. This technical part is handled in Section \ref{s.GRAD} via coupling arguments. 
Interestingly, when $\alpha>3$, vortices confined to a line start contributing for a positive fraction of the macroscopic fluctuations and the effective inverse temperature $\beta_{eff}$ ceases to be $\beta$ (See Remark \ref{r.NonM}).  

\enum

\subsection*{Acknowledgments.}
I wish to thank Mickaël Sousa for useful discussions and for pointing to us the reference \cite{frohlich1991phase}. I also wish to thank 
Elie Aïdekon, Loren Coquille, Diederik van Engelenburg, Louis Dupaigne,  Ivan Gentil,  Malo Hillairet, Ron Peled,  Julien Poisat, Pierre-François Rodriguez, Christophe Sabot,  Mikael de la Salle,  Bruno Schapira, Avelio Sepúlveda, Sylvia Serfaty, Gordon Slade, Jean-Marie Stéphan  and Yvan Velenik for very useful  discussions.


The research of C.G. is supported by the Institut Universitaire de France (IUF), the ERC grant VORTEX 101043450 and the French ANR grant ANR-21-CE40-0003.

\section{Preliminaries}

\subsection{Integer gaussians and correlation inequalities.}$ $

\begin{definition}[integer Gaussian]\label{}
Let $k\geq 1$ and $A$ be a $k\times k$ positive definite matrix. The (centred) \textbf{integer Gaussian} on $\Z^k$ induced by the quadratic form $A$ is the (centred) Gaussian vector on $\R^k$ with covariance matrix $A^{-1}$ conditioned to takes its values in $\Z^k$. In other words, if we denote by $\mathbb{P}_{A,k}^{IG}$ the law of this integer Gaussian, then for any vector $m \in \Z^k$, we have 
\begin{align*}\label{}
\mathbb{P}_{A,k}^{IG}(\psi=m) := \frac 1 {Z_A} \exp(-\frac 1 2 \<{m, A m})\,,
\end{align*}
where as usual, $Z_A$ is the normalization constant $Z_A:= \sum_{m \in \Z^k}  \exp(-\frac 1 2 \<{m, A m})$. 
\end{definition}

\begin{definition}[Sine-Gordon vector]\label{}
Let $k\geq 1$, $A$ be a $k\times k$ positive definite matrix 
and $\lambda\in (0,\infty)$ be a positive coupling constant. The $\lambda$-\textbf{Sine-Gordon} random vector $\psi=(\psi_1,\ldots,\psi_k)$ on $\R^k$ induced by the quadratic form $A$ has the following Radon-Nikodym derivative w.r.t  the Lebesgue measure $\mu$ on $\R^k$
\begin{align*}\label{}
\frac{ d\mathbb{P}_{A,k,\lambda}^{SG}}{d\mu}(\psi) := \frac 1 {Z_{A,\lambda}} \exp(-\frac 1 2 \<{\psi,  A \psi}) \exp(\lambda \sum_{i=1}^k  \cos(2\pi \psi_i))\,,
\end{align*}
where $Z_{A,\lambda}:= \int  \exp(-\frac 1 2 \<{\psi, A \psi}) \exp(\lambda \sum_{i=1}^k  \cos(2\pi \psi_i)) \mu(d\psi) $. 
\end{definition}
Clealry, as one varies $\lambda$ from 0 to $\infty$, the probability measure $\P_{A,k,\lambda}^{SG}$ interpolates between the Gaussian vector on $\R^k$ with covariance matrix $A^{-1}$ (denoted by $\calN_{A^{-1}}$) and $\P_{A,k}^{IG}$.  Let us state our first correlation inequality which is based on this interpolation. It is a very useful consequence of Ginibre inequality (\cite{Ginibre}) as detailed for example in \cite{RonFS}: 
\begin{proposition}\label{pr.ginibre}
For any positive definite matrix $A$ on $\R^k$, any coupling constant $\lambda \geq 0$ and any test vector $v\in \R^k$, we have 
\begin{align*}\label{}
\EFK{A,k}{IG}{e^{\<{v,\psi}}} \leq \EFK{A,k,\lambda}{SG}{e^{\<{v,\psi}}} \leq \EFK{\calN_{A^{-1}}}{}{e^{\<{v,\psi}}}\, \left(= e^{\frac 1 2 \<{v, A^{-1} v}} \right)\,.
\end{align*}
\end{proposition}

This inequality says that $k$-dimensional Gaussian vectors which are conditioned to be in $\Z^k$ fluctuate less than the (unconditional) Gaussian vector. This statement is rather intuitive and the goal of this paper is to show the reverse inequality in the case of the discrete Gaussian Chain at high temperature when $\alpha\in[2,3)$.  
\smallskip

We shall need also the following highly useful correlation inequality between two integer-valued Gaussians which is proved in the paper {\em an inequality for Gaussians on lattices} \cite{regev2017inequality} (see also Proposition 2.2 in \cite{aizenman2021depinning}). 

\begin{theorem}[\cite{regev2017inequality}]\label{th.regev}
If $A$ and $B$ are two positive definite matrices on $\R^k$ such that (in the sense of quadratic forms), $A \leq B$, then for any test vector $v\in \R^k$, 
\begin{align*}\label{}
\EFK{B,k}{IG}{e^{\<{v,\psi}}}  \leq  \EFK{A,k}{IG}{e^{\<{v,\psi}}}\,.
\end{align*}
In particular 
\begin{align*}\label{}
\EFK{B,k}{IG}{\<{v,\psi}^2}  \leq  \EFK{A,k}{IG}{\<{v,\psi}^2}\,.
\end{align*}
\end{theorem}


\subsection{Fractional and discrete fractional Laplacian.}\label{ss.DFL}

We start be briefly introducing the fractional Laplacian on $\R$. We refer to \cite{lodhia2016fractional,di2012hitchhiker} for excellent references on $\R^d$, $d\geq 1$. Given a smooth function $f\in C_c(\R)$, and given  $u\in(0,1)$, we may define the \textbf{$u$-fractional Laplacian} of $f$, which we shall denote\footnote{We use the minus sign so that the operator is positive definite on $C_c(\R)$.} by $(-\Delta)^u$ in the following two equivalent ways (see the above references for details):
\begin{align}\label{e.FLC}
 (-\Delta)^u f  &:=  \calF^{-1}\Big( \xi \in \R \mapsto |\xi|^{2u} \calF[f](\xi)  \Big) \\
 (-\Delta)^u f(x)  &:=  -  C(u) \int_\R  \frac{f(y)-f(x)}{|y-x|^{1+2u}} dy \,,\,\, \forall x\in \R
\end{align}
Where $\calF,\calF^{-1}$ denote the Fourier and Fourier inverse operators and where $C(u)>0$ is a positive constant so that both definitions match (\cite{lodhia2016fractional,di2012hitchhiker}).

The fractional Laplacian $(-\Delta)^u$ can be defined beyond $u\in(0,1)$, but we will not consider this regime  in this paper. Note also that most classical references on this topic use the parameter $s$ instead of $u$. We use  $u$ since the parameter $s$  is used throughout in Sections \ref{s.Bessel},\ref{s.proof}.  

In the rest of the paper, the correspondance between the fractional power $u$ of the Laplacian, the parameter $\alpha$ and the Hurst index $H$ will read as follows:
\begin{align}\label{e.uaH}
u = \frac{\alpha -1} 2  = H+\frac 1 2 \,\,\, \text{(for any $\alpha\in[2,3)$).} 
\end{align}

We now turn to the \textbf{discrete $u$-fractional Laplacian} on $\Z$.
For any function $f:\Z \to \R$ with compact support,  we define by analogy with the continuous case 
\begin{align}\label{e.FLu}
(-\Delta_\Z)^u f(x) := - \sum_{y\in \Z, y\neq x} \frac{f(y)-f(x)}{|y-x|^{1+2u}}\,,\,\,\, \forall x\in \Z. 
\end{align}

In the continuous case, there is essentially a unique (up to mult.  constant) natural definition of the $u$-fractional Laplacian. This is no longer the case in the discrete setting due to the lack of ``space rescaling''. For example, Remark \ref{r.FLfourier} below gives another equally legitimate definition of a discrete $u$-fractional Laplacian motivated by Fourier inversion. 

Since there are many natural choices for a $u$-fractional Laplacian, we introduce the following class of such Laplacians. For any $\alpha>1$ (which corresponds to $u=\frac{\alpha-1} 2$ when $\alpha\in [2,3)$),  and any sequence of non-negative coupling constants $\{J_\alpha(r)\}_{r\geq 1}$ satisfying $J_\alpha(r)\sim_{r\to \infty} \frac c {r^\alpha}$ for some $c>0$, we define the \textbf{$J_\alpha$-fractional Laplacian}:
\begin{align}\label{e.FLJ}
(-\Delta_\Z)^{J_\alpha} f(x) := - \sum_{y\in \Z, y\neq x } J_\alpha(|x-y|)\big(f(y)-f(x) \big)\,, \,\,\, \forall x\in \Z. 
\end{align}

In this work, we will work with coupling constants $J_\alpha(r)$ which are induced by certain random walks in $\Z^2$ (or the diamond graph $\calD$) as defined in~\eqref{e.J2} and ~\eqref{e.Jalpha}.

\begin{remark}\label{r.FLfourier}
Recall that for $f:\Z \to \R$ (say with compact support), if we define $\calF_\Z[f](\theta) = \hat f(\theta) = \sum_{n\in \Z} f(n) e^{i n \theta}$, then its Fourier inverse is given by $\calF_\Z^{-1}[\hat f] : n \mapsto \frac 1 {2\pi} \int_0^{2\pi} \hat f (\theta) e^{-i n \theta} d\theta$.

Notice that if $(-\Delta) f(x):= - \frac  {f(x+1)+f(x-1)- 2f(x)} 2$, then 
\begin{align*}\label{}
\calF_\Z[(-\Delta) f](\theta)&  =  (1-\cos(\theta))\,  \calF_\Z[f](\theta) \,\, 
\end{align*}
This then suggests another natural definition of the discrete fractional Laplacian for any $u\in (0,1)$:
\begin{align*}\label{}
(-\Delta_{\Z})_\mathrm{altern.}^u f (x) :=  \calF_\Z^{-1}\Big[ \theta \mapsto (1-\cos \theta)^u \calF_\Z[f](\theta)\Big]\,.
\end{align*}
Let us then define for every $r\geq 1$ the coupling constant 
\begin{align*}\label{}
J(r):= - \frac 1 {2\pi} \int_{-\pi}^\pi (1-\cos(\theta))^u \cos(r \theta) d\theta\,,
\end{align*}
which is easily seen to be positive for all $r\in \N^*$. We readily obtain that $(-\Delta_\Z)_\mathrm{altern.}^u f$ corresponds using the above notations to the operator $(-\Delta_\Z)^J$. Furthermore, if $u\in (0, 1)$, we may indeed check that 
\begin{align*}\label{}
J(r):= - \frac 1 {2\pi} \int_{-\pi}^\pi (1-\cos(\theta))^u \cos(r \theta) d\theta \sim_{r\to \infty}  \frac {c(u)} {r^\alpha}\,,
\end{align*}
with $c(u)>0$ and $\alpha$ such that  $u=\frac {\alpha-1} 2$ (see \cite{ciaurri2015connection}). 

\end{remark}

\subsection{Fractional Brownian motion.}  

Fractional Brownian motion $(B_t^H)_{t\in \R}$ are celebrated stochastic processes which were invented by Kolmogorov in the context of turbulence. See \cite{kolmogorov1940wienersche, mandelbrot1968fractional, whitt2002stochastic}. The parameter $H\in[0,1)$ stands for the {\em Hurst index} and describes the Hölder-regularity of the process (i.e. for any $H\in(0,1)$, the $H$-fractional Brownian motion is a.s. $H-\eps$ Hölder for all $\eps>0$).  
On the real line $\R$, they are defined as follows:

\begin{definition}\label{d.fBm}
For any $H\in(0,1)$, the $H$-fractional Brownian motion  ($H$-fBm) on $\R$ with Hurst index $H$ is the Gaussian process $(B_t^H)_{t\in \R}$ characterized by 
\bnum
\item $B_{t=0}^H =0$
\item $\Cov{B_s^H,B_t^H} = \tfrac 1 2 \left( |s|^{2H} + |t|^{2H} - |s-t|^{2H} \right)$ 
\enum
(N.B. The case $H=0$ will also be of relevance to us and corresponds to a log-correlated field).  
\end{definition}
Notice the case $H=\tfrac 1 2$ corresponds to a standard Brownian motion. For any $H\in(0,1)$, one may construct (by Kolmogorov criterion) versions of this Gaussian process which are a.s. continuous in $t$ and which are $H^-$  Hölder a.s.  

Fractional Brownian motions have many appealing properties:
\bi
\item {\em Self-similarity.}  For any $\lambda>0$, 
\begin{align*}\label{}
\lambda^{-H} (B_{\lambda t}^H)_{t\in \R}  \overset{(d)}= (B_{t}^H)_{t\in \R}  
\end{align*}

\item  {\em Stationary increments.}  (but not independent except when $H=\tfrac 1 2$). This is not quite immediate by looking at the above Covariance formula but it is more transparent when looking at 
\begin{align*}\label{}
\Eb{|B_t^H - B_s^H|^2} = |t-s|^{2H} \,,\,\,\,\text{  for all }s,t\in\R
\end{align*}
\ei

Similarly as for the definition of the $u$-fractional Laplacian, there exist many  equivalent definitions of $H$-fBm. For example another well-known construction of the $H$-fBm is obtained by expressing $B_t^H$ as a stochastic integral of a white noise $W(dx)$ against an explicit kernel.  See for example \cite{whitt2002stochastic}. Among all these equivalent definitions, we will give the one below which is closest to our focus in this paper, namely a Gibbs version of the $H$-fBm when $H\in[0, \tfrac 1 2)$ (notice we include $H=0$ here). We will introduce this viewpoint in an informal way only and we refer to  \cite{lodhia2016fractional} for a rigorous account on this Gibbs formalism in any dimension $d\geq 1$. The discrete version of the next Subsection (which will be rigorously defined) will also shed light on this informal viewpoint. 

To motivate the Gibbs version of $H$-fBm, recall that the Gaussian free field $\Phi$ on the entire $\R^2$ and viewed modulo additive constants corresponds to the following probability measure 
\begin{align*}\label{}
\mu^{\GFF,\free}_{\R^2}\big[ d\Phi \big] \propto \exp\left(- \tfrac 1 2 \int_{\R^2} \| \nabla \Phi(x) \|_2^2 dx \right) D\Phi\,,
\end{align*}
where $D\Phi$ stands for the ``Lebesgue'' measure on all possible fields $\Phi$.

\begin{definition}[Gibbs definition (informal). \cite{lodhia2016fractional}]\label{d.fBm2}
For any $H\in [0,\tfrac 1 2 )$, the $H$-fBm on $\R$ viewed modulo additive constants is the Gaussian measure
\begin{align}\label{e.GibbsfBm}
\mu^{H,\free}_{\R}\big[ dB^H \big] \propto \exp\left(- \tfrac 1 2 \iint_{\R \times \R} \frac{(B_t^H-B_s^H)^2}{|t-s|^{\alpha(H)}} dt ds \right) D B^H\,,
\end{align}
where $DB^H$ stands for the ``Lebesgue'' measure on all possible paths $t\mapsto B_t^H$ and with 
\begin{align*}\label{}
H\in [0,\tfrac 1 2) \text{   and   } \alpha =\alpha(H) = 2H+2\in [2,3)\,.
\end{align*}
(N.B. In definition \ref{d.fBm}, the translation symmetry was broken by fixing $B_{t=0}^H=0$).
\end{definition}
Using~\eqref{e.FLC}, notice that the quadratic form in~\eqref{e.GibbsfBm} is nothing but (up to a mult. constant): 
\begin{align*}\label{}
\<{B^H, (-\Delta)^u B^H} \,\,\,\text{  with  } u = \frac{\alpha -1} 2  = H+\frac 1 2\,,
\end{align*}
as in~\eqref{e.uaH}.  Using the rigorous setting in \cite[Section 3]{lodhia2016fractional}, this implies that Definitions \ref{d.fBm} and \ref{d.fBm2} are equivalent (again up to a mult. constant).

\smallskip
Of special relevance to our present setting, the paper \cite{lodhia2016fractional} extends this definition to the case where the paths are restricted to an open interval $D\subset \R$ endowed with Dirichlet boundary conditions. (In \cite{lodhia2016fractional}, fractional Brownian motion $B^H$ are extended to any dimensions $d\geq 1$, leading to \textbf{fractional Gaussian fields}. See Section 3 in \cite{lodhia2016fractional} as well as Section 4 which deals with fractional Gaussian field on a domain $D\subset \R^d$). Let us stick to the one-dimensional case $d=1$ where the fractional Brownian motion in a bounded open interval $D\subset \R$ is best defined via its Gibbs measure as follows (see \cite{lodhia2016fractional} for the details). 

\begin{definition}[Gibbs definition in a finite interval (informal)]\label{d.fBm3}
Let $D\subset \R$ be a bounded open interval, $\beta>0$ and $H\in [0,\tfrac 1 2 )$.  The $H$-fBm on $D$ with Dirichlet boundary conditions on $\p D$ and with inverse temperature $\beta$ is   the following Gaussian measure
\begin{align}\label{e.fBmD}
\mu^{H,0}_{D,\beta}\big[ dB^H \big] \propto \exp\left(- \tfrac \beta 2 \iint_{\R \times \R} \frac{(B_t^H-B_s^H)^2}{|t-s|^{\alpha(H)}} dt ds \right) 1_{B^H \equiv 0 \text{ on } D^c} D B^H\,.
\end{align}
(Note that the indicator function that $B_t^H=0$ outside of $D$ induces a long-range rooting effect for $B_t^H$ in the bulk of $D$). 
\end{definition}

As observed in \cite{lodhia2016fractional}, if $D=(-1,1)$, the covariance kernel of the above field can be computed explicitly when $H\in(0, \tfrac 1 2)$  thanks to \cite[Corollary 4]{blumenthal1961distribution}: there exists a constant $k(H)>0$ so that for any $x,y \in (-1,1)$, 
\begin{align*}\label{}
\EFK{(-1,1),\beta}{H,0}{B_x^H B_y^H} = \frac{k(H)} \beta |x-y|^{2H}
\int_0^{\frac{(1-|x|^2)(1-|y|^2)}{|x-y|^2}}
(v+1)^{- \tfrac 1 2}v^{H-\tfrac 1 2} dv\,.
\end{align*}
This formula is rather intimidating, yet by letting $x\to y$, notice it gives 
\begin{align}\label{e.VarH}
\Eb{(B_x^H)^2} =  \frac{\tilde k(H)}{\beta} (1-|x|^2)^{2H}\,\, \text{  for any }x\in [-1,1]\,. 
\end{align}
We make two observations: (1) as $H\to \tfrac 1 2$, we recover the variance profile of a Brownian Bridge over $D=(-1,1)$. (2) as $H\to 0^+$, notice that the fluctuations in the bulk are asymptotically insensitive to the distance to the boundary. (This is in agreement with the behaviour of the leading order fluctuations of a $2d$ GFF which corresponds once restricted to a $1d$ line to $H=0$).

\subsection{Discrete Fractional fields and invariance principle.}

Let us now introduce a finite dimensional Gaussian approximation of the Gibbs measure~\eqref{e.fBmD}. We start with a version on $\Z$ before rescaling it in ``time'' and space. 

\begin{definition}[The Gaussian chain]\label{d.GC}
For any $\alpha>1$, any set of coupling constants $\{J_\alpha(r)\}_{r\geq 1}$ s.t. $J_\alpha(r) \sim c r^{-\alpha}$, any $\beta>0$ and any $n\geq 1$, we define the \textbf{$J_\alpha$-Gaussian chain} on $\Lambda_n:=\{-n,\ldots,n\} \subset \Z$ at inverse temperature $\beta$  to be the Gaussian vector $(\varphi_{n}(x))_{x\in \Z}$:
\begin{align}\label{e.GC}
\FK{\Lambda_n,J_\alpha,\beta}{0}{d \varphi_n} \propto 
\exp
\left(-\tfrac \beta 2 \sum_{i,j \in \Z} J_\alpha(|i-j|)\left(\varphi_n(i) - \varphi_n(j)\right)^2
\right) 1_{\varphi_n \equiv 0 \text{  on  } \Lambda_n^c} \prod_{i\in \Lambda_n} d\varphi_n(i) \,.
\end{align}
(N.B.  the superscript $0$ stands for Dirichlet boundary conditions on $\Lambda_n^0$). 

The Gaussian chain with free boundary conditions on $\Z$ is also well-defined (as usual up to a global additive constant). It corresponds to the Gaussian process $(\varphi(x))_{x\in \Z}$
\begin{align}\label{e.GCfree}
\FK{\Z,J_\alpha,\beta}{free}{d \varphi} \propto 
\exp
\left(-\tfrac \beta 2 \sum_{i,j \in \Z} J_\alpha(|i-j|)\left(\varphi(i) - \varphi(j)\right)^2
\right)  \prod_{i\in \Z} d\varphi(i) \,.
\end{align}
\end{definition}

In what follows, \uline{we assume that $\alpha\in[2,3)$}, which corresponds to $H=\frac {\alpha-2} 2\in[0, \tfrac 1 2 )$. We are aiming at an {\em invariance principle} of this Gaussian chain towards an $H$-fractional Brownian motion. Let us start with the case of Dirichlet boundary conditions. 
For any $\beta>0$, $n\geq 0$, $\{J_\alpha(r)\}_{r\geq 1}$ with $J_\alpha(r) \sim \frac c {r^\alpha}$, let us consider the \uline{rescaled process} on $\frac 1 n \Z$:
\begin{align}\label{e.rescaled}
\bar \varphi_n(t):=  \frac 1 {n^{H}}  \varphi_n(n t)\,\,, \,\,\, \forall t\in \frac 1 n \Z\,.
\end{align}
We may now state the invariance principle below which follows from  \cite[Proposition 12.2]{lodhia2016fractional}. 

\begin{proposition}[Proposition 12.2 in \cite{lodhia2016fractional}]\label{pr.sheffield}
For any $\alpha\in[2,3)$ and any sequence of coupling constants $J_\alpha(r)\sim c r^{-\alpha}$, there exists a constant $K=K(J_\alpha)>0$ such that the following holds.
 
If $\bar \varphi_n$ is the rescaled $J_\alpha$-Gaussian chain at inverse temperature $\beta$, then the distribution $h_n:=\sum_{t\in \frac 1 n \Z} \bar \varphi_n(t)  \frac 1 n \delta_t$ converges in law to  $\frac{K}{\sqrt{\beta}} B_t^H$ (where  $B_t^H$ is the Dirichlet $H$-fBm on $D=(-1,1)$ defined in Definition \ref{d.fBm3}). The convergence is in the sense that for any test functions $f_1,\ldots,f_k \in C_c(D)$, 
\begin{align*}\label{}
(\<{h_n,f_1},\ldots,\<{h_n,f_k}) \,\,\, \overset{(d)}\longrightarrow_{n\to \infty} \,\,\,\, \frac{K}{\sqrt{\beta}} (\<{B^H,f_1},\ldots,\<{B^H,f_k})\,.
\end{align*}
\end{proposition}


Proposition 12.2 in \cite{lodhia2016fractional} handles the case where $J_\alpha(r) = \frac 1 {r^\alpha}$. To also handle the present case where $J_\alpha(r) \sim_{r\to \infty} c r^{-\alpha}$, we note that the same proof as in \cite{lodhia2016fractional} applies, except one has to  ensure that the random walk on $\Z$ with Markov kernel 
\begin{align*}\label{}
P_{J_\alpha}(i,j):= \frac{J_\alpha(|i-j|) 1_{i\neq j} }{2\sum_{r\geq 1} J_{\alpha}(r)}
\end{align*}
is also in the bassin of attraction of the symmetric $2u$-stable process with  $2u=\alpha-1$ (recall~\eqref{e.uaH}). This is indeed the case: indeed it is well-known (see for example \cite[Theorem 4.5.2]{whitt2002stochastic}) that if $P(i,j)\sim \frac {\tilde c} {|i-j|^\alpha}$ as $i-j \to \pm \infty$, then it is a sufficient condition for being in the bassin of attraction of the $2u$-stable process (where the scaling of the limiting stable process only depends on $\tilde c$).

\begin{remark}\label{}
We also expect the same invariance principle should  hold also for the Gaussian chain with free boundary conditions on $\Z$ (see~\eqref{e.GCfree}).  In this case the rescaled field $\left(n^{-H} \varphi_n(t)\right)_{t\in \frac 1 n \Z}$ rooted at 0 at $t=0$ should converge to $K \beta^{-1/2} B_t^H$, the $H$-fBm defined either in Definitions \ref{d.fBm} or \ref{d.fBm2} (the constant $K$ will depend on the chosen definition). The proof of Proposition 12.2 from \cite{lodhia2016fractional} would also apply to this case except \cite[Lemma 12.3]{lodhia2016fractional} which would need to be proved for the $2u$-``stable'' walk stopped when it first hits the origin as opposed to when it first hits $D^c$. 
\end{remark}

\begin{remark}\label{}
When $H\in (0,\tfrac 1 2)$, the limiting process $t\in[-1,1] \mapsto B_t^H$ is a.s. a continuous $H^-$ Hölder continuous functions with zero boundary conditions. We thus  expect that the  convergence in the above Proposition also holds for the linear interpolation of $t \mapsto \bar \varphi_n(t)$ under the stronger topology of the uniform convergence of continuous functions on $[-1,1]$. (N.B. Proposition 12.2 in \cite{lodhia2016fractional} handles more general fractional Gaussian fields in $D \subset \R^d$, including some cases where the limiting field is not a.s. a function, for example the present log-correlated case where $H=0$, this explains why \cite{lodhia2016fractional} did not need this stronger notion of convergence).  Combining this with the variance formula~\eqref{e.VarH}, this would then give a precise asymptotics for  $\Var{\varphi_n(x=0)}$.  We point out that it may also follow from the precise asymptotics of the {\em Harmonic potential} of walks in the bassin of attraction of stable processes in \cite{chiarini2023fractional}.
\end{remark}

\subsection{Discrete Gaussian chain, domain considered, infinite volume limits.}\label{ss.IVL}

In the rest of the paper, for any $n\geq 1$, 
\bi
\item $\Lambda_n$ will denote $\{-n,\ldots, n\}$. If we prescribe Dirichlet boundary conditions on $\Lambda_n$, this will not correspond to just fixing the field to be zero on $-n-1$ and $n+1$ but instead on the entire $\Z \setminus \Lambda_n$ (this is due to the long-range interactions). 
\item $\Lambda_n^2$ will denote the $2d$ box $\{-n,\ldots, n\}^2 \subset \Z^2$. In most of the paper (except in the short Section \ref{s.visible} devoted to long-range IV GFF in two-dimension), Dirichlet boundary conditions on $\Lambda_n^2$ will correspond to set the field to be zero on $\p \Lambda_n^2:= \{ y\in \Z^2 \setminus \Lambda_n^2, \text{ s.t. } \exists x\in \Lambda_n^2, x\sim y \}$.  
\ei

Following our definition of the Gaussian chain (Definition \ref{d.GC}) we thus recall  the main object of this paper:

\begin{definition}[The discrete Gaussian chain]\label{d.dGC}
For any $\alpha>1$, any set of coupling constants $\{J_\alpha(r)\}_{r\geq 1}$ s.t. $J_\alpha(r) \sim c r^{-\alpha}$, any $\beta>0$ and any $n\geq 1$, we define the \textbf{$J_\alpha$-discrete Gaussian chain} on $\Lambda_n:=\{-n,\ldots,n\} \subset \Z$ at inverse temperature $\beta$  to be the integer-valued field $(\Psi_{n}(x))_{x\in \Z}$:
\begin{align}\label{e.dGC}
\FK{\Lambda_n,J_\alpha,\beta}{dGC,0}{\Psi_n} \propto 
\exp
\left(-\tfrac \beta 2 \sum_{i,j \in \Z} J_\alpha(|i-j|)\left(\Psi_n(i) - \Psi_n(j)\right)^2
\right) 1_{\Psi_n(i)=0,\,\, \forall i\in \Lambda_n^c} \,.
\end{align}

The discrete Gaussian chain with free boundary conditions on $\Z$ is given by the measure
\begin{align}\label{e.dGCfree}
\FK{\Z,J_\alpha,\beta}{free}{(\Psi(x))_{x\in \Z}} \propto 
\exp
\left(-\tfrac \beta 2 \sum_{i,j \in \Z} J_\alpha(|i-j|)\left(\Psi(i) - \Psi(j)\right)^2
\right)  1_{\Psi(0)=0}\,.
\end{align}
\end{definition}
Let us briefly explain why the second definition is well defined. This is not obvious as it corresponds to an infinite volume limit. The idea is to condition the Gaussian chain defined in~\eqref{e.GCfree} to take integer values only in the finite window $\Lambda_n \subset \Z$. One may then let $n\to \infty$. Using Ginibre's inequality as for Proposition \ref{pr.ginibre} (see the interpolation proof in \cite{RonFS}), one can see that Laplace transforms of any test function $f$ are decreasing functions of $n$. Together with tightness (which is also a consequence of Ginibre), this implies the existence of an infinite volume limit. Since we shall not focus on the free boundary conditions in this paper, we do not provide further details here. 
\medskip

Let us say a few words on the case $\alpha=2$ (i.e. $H=0$). As we shall see in the next section, one interesting way to realise the discrete Gaussian chain at $\alpha=2$ is as follows:
\bi
\item \uline{$\alpha=2$, free boundary conditions.} Consider the free Gaussian free field on $\Z^2$ with inverse temperature $\beta$ (say rooted at the origin) and condition this Gaussian field to take integer-values on the $1D$ line $\Z \setminus \{0\}$.  (N.B. in the case of the discrete Gaussian chain \cite{wirth2019maximum,bauerschmidt2022discrete,bauerschmidt2022discreteN2}, the field is conditioned to take integer-values on the whole lattice $\Z^2$). 

The integer-valued field it produces on the $1d$ line $\Z$ is the same as  the discrete Gaussian chain defined in~\eqref{e.dGCfree} with coupling constants $J_{\alpha=2}(r)$ provided by Appendix.


\item \uline{$\alpha=2$, Dirichlet boundary conditions.} We now consider a Gaussian free field on $\Z^2$ which is rooted at all points of $\Z\times \{0\} \setminus (\Lambda_n \times \{0\})$. See Figure \ref{f.diamond}.
If we condition this Gaussian field to take integer values on the finite set $\Lambda_n \times \{0\}$, we now recover the Dirichlet discrete Gaussian chain~\eqref{e.dGC} with same coupling constant $J_{\alpha=2}(r)$.  
\ei

These observations are consistent with the fact that the $H=0$-fBm from Definition \ref{d.fBm3} is a log-correlated field, more precisely here the restriction of a $2d$ continuous Gaussian free field to a line (see \cite{lodhia2016fractional}).

\begin{remark}\label{}
Definition \ref{d.dGC} provides a natural candidate for a discrete approximation of $H$-fractional Brownian motion when $H\in[0,\tfrac 1 2)$.  
Let us point out that Hammond and Sheffield gave another discrete process converging to $H$-fBm in \cite{hammond2013power}. It is not a Gibbs type construction but instead relies on heavy tailed random variables to generate the next step given all previous steps. Interestingly their construction works in the regime $H\in(\tfrac 1 2, 1)$ while ours works in the regime $H\in[0,\tfrac 1 2]$. 
\end{remark}

%
%



\section{Proof of the main Theorem (Theorem \ref{th.main}) in the case $\alpha_c=2$}\label{s.proof2}

We first focus on the proof of Theorem \ref{th.main} when $\alpha=\alpha_c=2$. This will pave the way for the analysis when $\alpha>2$. Before handling the discrete Gaussian chain model, we will prove quantitative versions of the {\em invisibility of integers} in the context of a $2D$ Gaussian free field on $\Lambda_n^2$ conditioned to take its values in the integers on a $1D$ line (say, $\Lambda_n \times \{0\} \subset \Lambda_n^2$) as well as on more general sub-domains, such as a fractal set $F_n$ in Theorem \ref{th.fractal}. (See also Figure \ref{f.fractal}).  We shall analyze {\em free, periodic and Dirichlet} boundary conditions. The later ones are the most challenging ones because of the presence of {\em non-neutral} charges in the Coulomb gas expansion from \cite{FS, wirth2019maximum}. They are also the ones for which the results are the easiest to state and apprehend, so we will start by those in Subsection \ref{ss.Dir}.

\subsection{The trace of the $2D$ simple random walk.}

\begin{figure}[!htp]
\begin{center}
\includegraphics[width=0.7\textwidth]{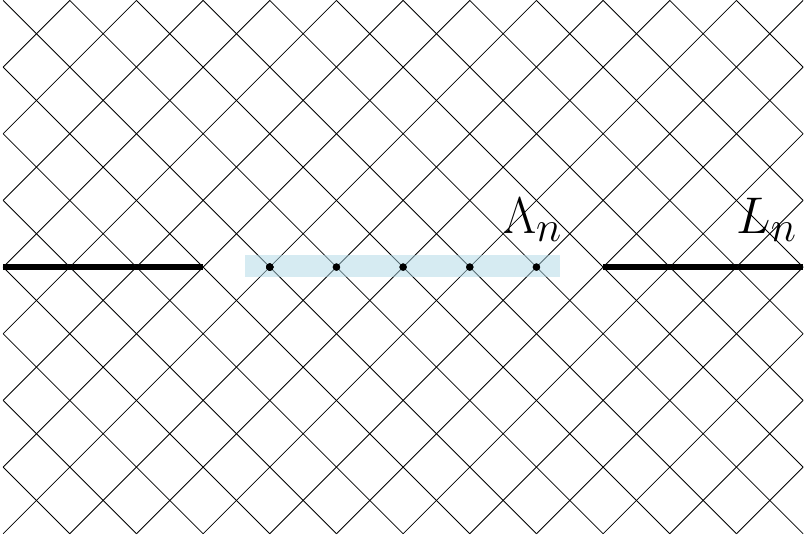}
\end{center}
\caption{}\label{f.diamond}
\end{figure}

Let us consider the {\em diamond graph} $\calD = e^{i \pi/4} 2^{-1/2} \Z^2$ (see Figure \ref{f.diamond}). Notice that the horizontal $1D$ line $\calD\cap \R \times \{0\}$ is given by $\Z \times \{0\}$. With a slight abuse of notation, we will denote this horizontal line by $\Z$. 
  
Let $Z_n$ be a simple random walk on the diamond graph $\calD$ starting at the origin $(0,0)$. We introduce the following sequence of stopping times.
\begin{align*}\label{}
\begin{cases}
& \tau_0:=0 \\
& \forall k\geq 0, \tau_{k+1}:= \inf \{ n>\tau_k, Z_n\in \Z \}  \,.
\end{cases}
\end{align*}
As such, the random times $\{ \tau_k \}_{k\geq 0}$ are the successive return times of the walk $(Z_n)$ to the horizontal line $\Z\times \{0\} =\Z$. 
This procedure induces a $1D$ Markov process $(X_k)_{k\geq 0}$ with long-range jumps defined by 
\begin{align*}\label{}
\begin{cases}
& X_0:=0 \\
& \forall k\geq 1, X_k:= \text{ horizontal position along }\Z \text{ of the walk } Z_{\tau_k}\,. 
\end{cases}
\end{align*}
This process $\{X_k\}_{k\geq 0}$ is thus the trace of the $2D$ walk $Z_n$ on the line $\Z$. It is well-known that the transition kernel of this Markov chain decays like $1/r^2$. 
The reason why we have chosen the diamond graph $\calD$ is  due to the following exact formula which goes back to Spitzer (see \cite[E8.3]{spitzer2013principles} as well as  \cite{amir2016one}):
\begin{align}\label{e.PEx}
\FK{}{(0,0)}{Z_{\tau_{\Z}} = (0,k)} = \FK{}{0}{X_1=k} = 
\begin{cases}
& \frac {2}{\pi(4k^2 -1)}  \text{    if } k\neq 0 \\
& 1 - \frac 2 \pi   \,\, \,\,\,\, \text{        if } k = 0 
 \end{cases}
\end{align}

\subsection{The trace of a $2D$ Gaussian Free Field.}\label{ss.gff}

Let us consider the Gaussian free field $\varphi_n$ on the diamond graph $\calD$ and rooted on the set 
\begin{align}\label{e.Ln}
L_n:=\left(\{ \ldots, -n-1\} \cup \{n+1,\ldots \} \right) \times \{0\} \subset \calD\,.
\end{align}
(See Figure \ref{f.diamond}). We shall use the following normalization
\begin{align}\label{e.gff8}
\FK{}{}{d\varphi_n} \propto \exp\left( - \frac \beta 8 \sum_{i\sim j \in \calD} (\varphi_n(i)-\varphi_n(j))^2\right) 1_{\varphi_n \equiv 0 \text{  on  } L_n}\,,
\end{align}
so that the covariance structure of this field corresponds exactly the Green function of the simple random walk (without further factor $1/4$).

Then, if we denote by a slight abuse of notation the interval $\Lambda_n := \Lambda_n \times \{0\} = \{-n,\ldots,n\} \times \{0\} \subset \calD$, we have the following identity in law. 
\begin{lemma}\label{l.gffs0}
The restriction of the $\beta$-GFF field $\varphi_n$ to the $1D$ interval $\Lambda_n \subset \calD$ is equal in law to the Gaussian vector $\{\bar \varphi_i\}_{i\in \Z}$  with density 
\begin{align*}\label{}
\exp(-\frac {\beta} 2 \sum_{i \neq j} J_{\alpha=2}(i-j) (\bar \varphi_i - \bar \varphi_j)^2) 1_{\bar \varphi_i = 0 \forall i \notin \Lambda_n}\,,
\end{align*}
where the coupling constants $J_2(r)$ are defined by 
\begin{align}\label{e.J2}
J_2(r) :=  \frac {2}{\pi(4k^2 -1)}  \text{    (if  $r\neq 0$)}
\end{align}
\end{lemma}

\ni
{\em Proof.} 
Such a property has already been used implicitly several times in the literature (see for example recently in \cite{aru2022percolation}). To see why this holds, it is sufficient to compute the covariance matrix of the field $\{\bar \varphi_i\}_{i\in \Lambda_n}$ and to check that it is the same as the restriction of the covariance matrix of the larger Gaussian vector $\{(\varphi_n)_x\}_{x\in \calD}$. This is indeed the case as for any $i,j \in \Z$, 
\begin{align*}\label{}
\Eb{\bar \varphi_i \bar \varphi_j} & = \EFK{}{i}{\sharp \text{visits to $j$ for the RW $X_k$ killed when exiting $\Lambda_n$} } \\
& = \EFK{}{(i,0)}{
\sharp\text{visits to $(j,0)$ for the SRW $Z_n$ in $\calD$ killed when hitting $L_n$}} \\
& = \Eb{\varphi_n(i) \varphi_n(j)}\,,
\end{align*}
where the second equality is because $\{X_k\}$ is by definition the trace the $2D$ SRW $Z_n$. This ends the proof by definition of the GFF $\varphi_n$ with $0$-boundary conditions on $L_n:=\left(\{ \ldots, -n-1\} \cup \{n+1,\ldots \}\right) \times \{0\}$.
\qed

\subsection{Fröhlich-Spencer's proof of delocalisation of the discrete Gaussian.}\label{ss.sketchFS}

In this subsection, we briefly explain the main steps of Fröhlich-Spencer's proof of delocalisation of the IV-GFF (or discrete Gaussian) as we will rely extensively on their {\em Coulomb gas expansion} technique. 
We refer the reader to the excellent review \cite{RonFS} as well as to the original paper \cite{FS}. The text below follows closely the concise overview of Fröhlich-Spencer's proof given in \cite{GS1}.  We still include it here, first because it introduces the relevant notations for the rest of the proof and also because the emphasis is a bit different. In \cite{GS1}, the part of Fröhlich-Spencer's proof which relied on Jensen's inequality was troublesome for our setting in \cite{GS1} while here, the emphasis is rather on the quantitative estimates controlling the effective inverse temperature $\beta_{eff}(\beta)$. Also, as we will also consider Dirichlet boundary conditions, we shall also summarise Wirth's recent appendix in \cite{wirth2019maximum} on how to handle these boundary conditions. (In Section \ref{ss.NN}, we will extend the validity of Wirth's analysis to more general test functions).  
\smallskip

To start with, we fix a square domain $\Lambda \subset \Z^2$ and we consider the case of free boundary conditions rooted at some vertex $v\in \Lambda$. (We refer to \cite{RonFS} for a convenient rooting procedure which requires the condition $1_{\phi_v\in[-\pi, \pi)}$ instead of $1_{\phi_v=0}$).  The case of Dirichlet boundary conditions (\cite{wirth2019maximum}) will be discussed at length further in Section \ref{s.EXT}. 
\smallskip

The proof by Fröhlich-Spencer can essentially be decomposed into the following successive steps:

{\em 1) The first step}  is to view the singular conditioning $\{\phi_i \in 2\pi \Z, \forall i \in \Lambda\}$ using Fourier series\footnote{It is slightly more convenient to consider the GFF conditioned to live in $(2\pi \Z)^\Lambda$ rather than $\Z^\Lambda$. Of course by tuning $\beta$,  one may scale back to conditioning in $\Z^\Lambda$ if needed. Following \cite{FS,RonFS}, we will stick to this convention in the remaining of this paper.} thanks to the identity
\[
2\pi \sum_{m\in \Z} \delta_{2\pi m}(\phi)  \equiv 1 + 2 \sum_{q=1}^\infty \cos(q \phi)\,.
\]
To avoid dealing with infinite series, proceeding as in \cite{RonFS}, we consider the following approximate IV-GFF
\[
\P_{\beta,\Lambda,v}[d\phi] := \frac 1 {Z_{\beta,\Lambda,v}}  \prod_{i\in \Lambda} \left(1 + 2 \sum_{q=1}^N \cos( q \phi(i)) \right) \P_{\beta,\Lambda,v}^{\GFF}[d\phi]\,.
\]
In fact, more general measures are considered in \cite{FS,RonFS}:  they fix a family of trigonometric polynomials $\lambda_\Lambda := (\lambda_i)_{i\in \Lambda}$ attached to each vertex $i\in \Lambda$. These trigonometric polynomials are parametrized as follows: for each $i\in \Lambda$, 
\[
\lambda_i(\phi(i)) = 1 + 2 \sum_{q=1}^N \hat \lambda_i(q) \cos(q \phi(i))\,.
\]
It turns out that this more general viewpoint considered in \cite{FS,RonFS} will be of key importance to us. Indeed in \cite{FS,RonFS}, eventually they apply   the same trigonometric polynomial at each vertex of the graph and they let $N\to \infty$ in order to obtain fluctuations bounds on the IV-GFF. 
 In our case, the situation will be dramatically different: points on the middle line of the grid $\Z^2$ and points elsewhere will carry different trigonometric polynomials!  
\medskip 

Let us then consider an arbitrary family of trigonometric polynomials $\lambda_\Lambda = (\lambda_i)_{i\in \Lambda}$ and let us define 
\[
\P_{\beta,\Lambda,\lambda_\Lambda, v}[d\phi] := \frac 1 {Z_{\beta,\Lambda,,\lambda_{\Lambda},v}}  \prod_{i\in \Lambda}  \lambda_i(\phi(i)) \P_{\beta,\Lambda,v}^{\GFF}[d\phi]\,.
\]  

\medskip
{\em 2) The second step} in the proof is to fix a test function $f:\Lambda \to \R$ such that $\sum_{i\in \Lambda} f(i) =0$ and to consider the Laplace transform of $\<{\phi, f}$, $\EFK{\beta,\Lambda,\lambda_\Lambda, v}{}{e^{\<{\phi,f}}}$.


By a simple change of variables, this Laplace transform can be rewritten 
\[
\EFK{\beta,\Lambda,\lambda_\Lambda, v}{}{e^{\<{\phi,f}}} 
= \frac 1 {Z_{\beta,\Lambda,\lambda_\Lambda,v}} \exp(\frac 1 {2 \beta} \<{f, (- \Delta)^{-1} f}) 
\EFK{\beta,\Lambda,v}{\GFF}
{\prod_{i\in \Lambda} \lambda_i(\phi(i) +\sigma(i))}\,,
\]
where the function $\sigma=\sigma_f$ will be used throughout. It is defined by
\begin{align}\label{e.sigma}
	\sigma:= \frac 1 \beta [- \Delta]^{-1} f 
\end{align}

The main difficulty in the proof in \cite{FS} is in some sense to show that the effect induced by the shift $\sigma$ does not have a dramatic effect compared to the exponential term $\exp(\frac 1 {2 \beta} \<{f, (- \Delta)^{-1} f})$ so that ultimately, there exists an $\eps=\eps(\beta)$ which goes to zero as $\beta\to 0$ and which is such that 
\begin{align*}\label{}
	\EFK{\beta,\Lambda}{\IV}{e^{\<{\phi,f}}} \geq  \exp(\frac 1 {2 \beta(1+\eps)} \<{f, - \Delta^{-1} f})\,.
\end{align*}
From such a lower bound on the Laplace transform,  one can easily extract delocalisation properties of the IV-GFF. 

In the next subsection, we will observe that $\eps=\eps(\beta)$ can be taken to be zero in our setting. (Which is known to be wrong for the $2D$ discrete Gaussian as proved in \cite{GS2}). 
\medskip

{\em 3) The third (and  most difficult) step} is to control the effect of the shift $\sigma$ via a highly non-trivial expansion into Coulomb charges which enables us to rewrite the partition function as follows:
\begin{align*}\label{}
	Z_{\beta,\Lambda,\lambda_\Lambda,v} = \sum_{\calN \in \calF} c_\calN \int \prod_{\rho\in \calN} [1+z(\beta,\rho, \calN) \cos(\<{\phi,\bar\rho})] d\mu_{\beta,\Lambda,v}^\GFF(\phi)\,.
\end{align*}

We refer to \cite{FS,RonFS} for the notations used in this expression and in particular for the concept of {\em charges} (i.e. $\rho : \Lambda \to \R$), {\em ensembles} (i.e. sets $\calN$ of mutually disjoint charges $\rho$) etc. 

One important feature of this expansion into charges is the fact that under some (very general) assumptions on the growth of the Fourier coefficients $|\hat \lambda_i(q)|$ (see (5.35) in \cite{FS}), it can be shown that the effective activities $z(\beta, \rho, \calN)$ decay fast. Namely (see (1.14) in \cite{RonFS}),
\begin{align}\label{e.zfree}
|z(\beta, \rho, \calN)|
\leq \exp \left( -\frac c \beta (\|\rho\|_2^2 + \log_2(\mathrm{diam}(\rho) + 1))\right).
\end{align}
As such we see that at high temperature, the partition function corresponds to a sum of positive measures. (Also the weights $c_\calN$ are positive and s.t. $\sum c_\calN =1$). 

\begin{remark}\label{}
	In \cite{RonFS}, the authors have introduced a slightly different definition of the free b.c. GFF which makes the analysis behind this decomposition into charges more pleasant (their definition cures the presence of non-neutral charges $\rho$ very easily). One can switch to their more convenient definition in our setting since in the limit $N\to \infty$, both give the same integer-valued GFF. 
\end{remark}

This crucial third step thus allows us to rewrite the Laplace transform $\EFK{\beta,\Lambda,\lambda_\Lambda, v}{}{e^{\<{\phi,f}}}$ as follows: 
\begin{align*}\label{}
	e^{\frac 1 {2 \beta} \<{f, - \Delta^{-1} f}}
	\frac 
	{
		\sum_{\calN \in \calF} c_\calN \int \prod_{\rho\in \calN} [1+z(\beta,\rho, \calN) \cos(\<{\phi,\bar\rho} + \<{\sigma, \rho})] d\mu_{\beta,\Lambda,v}^\GFF(\phi)
	}
	{
		\sum_{\calN \in \calF} c_\calN \int \prod_{\rho\in \calN} [1+z(\beta,\rho, \calN) \cos(\<{\phi,\bar\rho})] d\mu_{\beta,\Lambda,v}^\GFF(\phi)
	}\,.
\end{align*}
We now rewrite this ratio as  (thus defining $Z_\calN(\sigma)$ and $Z_\calN(0)$)
\begin{align*}\label{}
	\EFK{\beta,\Lambda,\lambda_\Lambda, v}{}{e^{\<{\phi,f}}}
	= e^{\frac 1 {2 \beta} \<{f, - \Delta^{-1} f}}
	\frac
	{\sum_{\calN\in \calF} c_\calN Z_{\calN}(\sigma)}
	{\sum_{\calN\in\calF} c_\calN Z_\calN(0)}
\end{align*}

\medskip
{\em 4) The fourth step is an analysis} for each fixed ensemble $\calN\in \calF$ of the above ratio $\frac {Z_{\calN}(\sigma)} {Z_{\calN}(0)}$. Trigonometric inequalities are used here  in order to obtain for each $\calN$:
\begin{align*}\label{}
	\frac {Z_{\calN}(\sigma)} {Z_{\calN}(0)}  \geq &  
	\exp\big[-D_1 \sum_{\rho\in \calN} |z(\beta,\rho,\calN)|\<{\sigma,\rho}^2\big] \\
	& \times 
	\int \frac
	{e^{S(\calN,\phi)}} 
	{Z_{\calN}(0)}
	\prod_{\rho\in\calN}[1+z(\beta,\rho,\calN) \cos(\<{\phi,\bar \rho})]
	d\mu_{\beta,\Lambda,v}^\GFF(\phi)\,, \nn
\end{align*}
where
\begin{align}\label{e.Sba}
	S(\calN,\phi):=
	-\sum_{\rho\in\calN} \frac
	{z(\beta,\rho,\calN)\sin(\<{\phi,\bar \rho})\sin(\<{\sigma,\rho})}
	{1+z(\beta,\rho,\calN)\cos(\<{\phi,\bar \rho})}
\end{align}
Two crucial observations are made at this stage:
\bnum
\item The functional $\phi\mapsto S(\calN, \phi)$ is odd in $\phi$
\item The  measure $\prod_{\rho\in\calN}[1+z(\beta,\rho,\calN) \cos(\<{\phi,\bar \rho})] d\mu_{\beta,\Lambda,v}^\GFF(\phi)$ is invariant under $\phi \to -\phi$.
\enum
All together this simplifies tremendously the above lower bound, as by using Jensen, one obtains readily 
\begin{align}\label{e.nice1}
	\frac {Z_{\calN}(\sigma)} {Z_{\calN}(0)}  \geq &  
	\exp\big[-D_1  \sum_{\rho\in \calN} |z(\beta,\rho,\calN)|\<{\sigma,\rho}^2\big] \,.
\end{align}
The final step is the following upper bound on the contribution coming from the charges $\rho$ integrated against $\sigma$.  Namely for any large constant $D_2>0$, if $\beta$ is small enough then 
\begin{align}\label{e.1516}
\sum_{\rho\in \calN}  |z(\beta,\rho,\calN)|\<{\sigma,\rho}^2 \leq \frac \beta {D_2} \sum_{i\sim j \in \Lambda} (\sigma_i - \sigma_j)^2\,.
\end{align}
At this point, the fact the effective inverse temperature $\beta_{eff}(\beta)$ does not deviate too much from $\beta$, i.e. $\beta_{eff}(\beta)\leq \beta(1+O(\eps(\beta)))$, is obtained by choosing the constant $D_2$ much larger than $D_1$. 

\medskip

In the proof below, we will need to revisit the way this upper bound~\eqref{e.1516} is obtained by taking into account the fact our charges will belong to a one-dimensional line instead of  the entire $\Lambda\subset \Z^2$. In the special case of Dirichlet boundary conditions, new difficulties arise due to the presence of non-neutral charges in the Coulomb gas expansion. Those have been analyzed in \cite{wirth2019maximum} at least for certain type of test functions $f$.  We will need to extend the set of possible observables handled in \cite{wirth2019maximum} (this will be the purpose of Section \ref{ss.NN}). Both of these extensions of Fröhlich-Spencer's proof are not immediate and this is why  the full Section \ref{s.EXT} will be dedicated to these.



\subsection{The $2D$ GFF conditioned to take integer values on a $1D$ line. The case of Dirichlet boundary conditions.}\label{ss.Dir}

In this subsection, we shall prove the following more general version of Theorem \ref{th.line}. (Recall the notations introduced just before Theorem \ref{th.line}). 
\begin{theorem}\label{th.line2}
Consider $\Psi_n^{line}$ the GFF in $\Lambda_n^2 = \{-n,\ldots,n\}^2$ with Dirichlet boundary conditions and conditioned to take integer values on the line $\Z\times \{0\}$. Then, if $\beta$ is small enough,  the following two properties hold.
\bnum
\item For any $(x,y)\in D=(-1,1)^2$, we have the following convergence in law
\begin{align}\label{e.lessfine}
\frac 1 {\sqrt{\log n}} \Psi_n^{line}(\floor{tx},\floor{ty})  \overset{(d)}\longrightarrow 
\calN(0, \frac{2}{\pi \beta})\,.
\end{align}

\ni
Furthermore, away from the middle line, the variance has the following more precise asymptotics: for any fixed $w=(x,y)\in (-1,1)^2$, with $y\neq 0$, 
\begin{align}\label{e.fine}
\mathrm{Var}_\beta[ \Psi_n^{line}(\floor{x n},\floor{y n}) ]  = \frac 2 {\pi}\left( \log \frac n {r_D(w)}   + \gamma_{\mathrm{EM}} + \frac 1 2 \log 8 \right)   +o(1)
\end{align}
as $n\to \infty$, where $r_D(w)$ is the conformal radius of the domain $D=(-1,1)^2$ seen from the point $w=(x,y)\in D$ and $\gamma_{\mathrm{EM}}$ is the Euler-Mascheroni constant
\footnote{We refer to \cite[Theorem 1.17]{biskup2020extrema} for a proof of this asymptotics in the case of the Green function of the SRW on $\Z^2$ (see also \cite{aidekon2023multiplicative}).}.

\item For any  test function $g\in \calC_c^1([-1,1]^2)$, as $n\to \infty$
\begin{align}\label{e.Line}
\frac 1 {n^2} \sum_{z\in \Lambda_n^2} \Psi_n^{line}(z) g\big(\frac z n\big)  \overset{(d)}\longrightarrow 
\calN(0, \frac{4}{ \beta} \<{g,(-\Delta)^{-1} g})\,,
\end{align}
where $\Delta$ denotes here the continuous Laplacian on $D=(-1,1)^2$ with Dirichlet boundary conditions
\footnote{Note that $\pi$ is absorbed into the Green function of the continuous Laplacian and the factor 4 is due to our normalisation in~\eqref{e.gff8}. See for example \cite{lawler2010random}).}.  
\enum
\end{theorem}

These two properties show that integers are invisible when $\beta$ is small enough: namely $\Psi_n^{line}$ looks essentially as an unconditioned GFF $\varphi_n$ both in the {\em pointwise} sense (statement (1), in particular the refinement~\eqref{e.fine} away from the line) and also when {\em seen as a distribution} (statement (2)).

\medskip
\ni
{\em Proof of Theorem \ref{th.line2}.} 

Let us start with the proof of Item (1). We consider the GFF $\varphi_n$ in the domain $\Lambda_n^2$ with Dirichlet boundary conditions (N.B. This is not quite the same as the GFF $\varphi_n$ considered above in Subsection \ref{ss.gff} for which the domain is instead a ``slit-domain'' $\calD \setminus L_n$).  
In particular, we will follow here both Fröhlich and Spencer's proof but also Wirth's extension to Dirichlet boundary conditions \cite{wirth2019maximum}. 
\smallskip

We shall now condition the field $\varphi_n$ to take its values in $2\pi \Z$ on the line $\Lambda_n$ (the same proof works with a conditioning in $\Z$ rather than $2\pi \Z$, but as explained earlier notations are slightly simpler with this conditioning).

This corresponds to fixing the following family of trigonometric polynomials $\lambda:= (\lambda_i)_{i\in \Lambda_n^2}$ attached to each vertex $i\in \Lambda_n^2$.
\begin{align*}\label{}
\begin{cases}
& \lambda_i(\varphi_i) := 1 + 2 \sum_{q=1}^N \cos(q \varphi(i)) \;\;\;\; \text{  if  }i\in \Z\times \{0\}\cap \Lambda_n^2 \\
& \lambda_i(\varphi_i) := 1  \;\;\;\;\;  \text{  if  }i\notin \Z\times \{0\} \,,
\end{cases}
\end{align*}
and then by letting $N\to \infty$. Indeed such trigonometric polynomials will force the field to be in $2\pi \Z$  along the line (at least when $N\to \infty$), while it does not impose any conditioning outside the line. 
We refer to \cite[Section 5]{RonFS} for details on taking the limit $N\to \infty$. 
 
By choosing as a test function $f(z):=1_{z=(\floor{xn},\floor{yn})}$, if we run the same proof as Fröhlich and Spencer (\cite{FS}) with this choice of trigonometric polynomials,  we end up with the following lower bound on the Laplace transform of $\Psi_n^{line}(\floor{xn}, \floor{yn} )$
\begin{align*}\label{}
\EFK{\beta,\Lambda_n^2,\lambda, v}{}{e^{\<{\phi,f}}}
	= e^{\frac 1 {2 \beta} \<{f, (- \Delta)^{-1} f}}
	\sum_{\calN\in \calF} c_\calN \frac {Z_{\calN}(\sigma)}
	{Z_\calN(0)} 
\end{align*}
where $\sigma:= \frac 1 \beta (-\Delta)^{-1} f$.  Note that since $f:= 1_{z=(\floor{xn},\floor{yn})}$, we have 
\begin{align}\label{e.Green1}
 \<{f, (- \Delta)^{-1} f} = G_{\Lambda_n^2}(\floor{ xn}, \floor{y n}) = \frac {2 } {\pi }  \left( \log \frac n {r_D(x,y)} + c_0 \right) +o(1)\,,
\end{align}
by the result mentioned above from \cite{biskup2020extrema} (and where $c_0$ is explicit and is given in~\eqref{e.fine}).

It remains to bound from below the error term involving the expansion into charges. 
Using the uniform lower bound given in~\eqref{e.nice1}, we get 
\begin{align}\label{e.123}
\EFK{\beta,\Lambda_n^2,\lambda, v}{}{e^{\<{\phi,f}}}
\geq  e^{\frac 1 {2 \beta} \<{f, (- \Delta)^{-1} f}} \sum_{\calN\in \calF} c_\calN  \exp\big[-D_1  \sum_{\rho\in \calN} |z(\beta,\rho,\calN)|\<{\sigma,\rho}^2\big] \,.
\end{align}
The main observation now is that the non-trivial charges are initially  restricted to the line $\Z \times \{0 \}$. Yet the algorithm used in Fröhlich-Spencer's proof (\cite{FS}) does not move charges around, it only groups and splits charges together. 

If boundary conditions were free, then only \underline{neutral} charges would contribute in this sum which would make it easier to lower bound the term
\begin{align*}\label{}
\exp\big[-D_1  \sum_{\rho\in \calN} |z(\beta,\rho,\calN)|\<{\sigma,\rho}^2\big] 
\end{align*}

But in the case of Dirichlet boundary conditions, some of the charges $\rho\in \calN$ may not be neutral.  See \cite{wirth2019maximum}.

Despite the possible presence of non-neutral charges, the following key upper bound is proved in \cite{wirth2019maximum}
for certain test functions $f$:
\begin{lemma}[Claims 15 and 16 in \cite{wirth2019maximum}]\label{l.1516}
For any $D_2>0$, if $\beta$ is chosen small enough, the following holds. For all $n\geq 1$, if $\Lambda_n^2$ is equipped with Dirichlet boundary conditions,  then for any test function $f$ such that the following condition holds
\begin{align}\label{e.condition}
\sigma = \frac 1 \beta [-\Delta_{\Lambda_n^2}]^{-1} f  \text{   vanishes  on the annulus  } \Lambda_n^2 \setminus \Lambda_{\frac n 2}^2\,,
\end{align}
and for any ensemble of charges $\calN$ appearing in the Coulomb gas expansion (with potentially neutral as well as  non-neutral charges), one has 
\begin{align}\label{e.NonN}
D_1 \sum_{\rho\in \calN} |z(\beta, \rho, \calN)| \cdot \<{\sigma,\rho}^2 \leq \frac \beta {D_2} \sum_{j \sim l} (\sigma_j -\sigma_l)^2 
\end{align}

\end{lemma}

The condition~\eqref{e.condition} is rather restrictive. By relying only on such test functions, it turns out one would still manage to obtain a precise control on the fluctuations of the field inside the sub-domain $\Lambda_{\frac n 2}^2$. But one would loose a precise control of the fluctuations in the annulus $\Lambda_n^2 \setminus \Lambda_{\frac n 2}^2$. We will explain in Subsection \ref{ss.NN}   how to extend Wirth's result to all test functions, i.e. we shall prove in Subsection \ref{ss.NN}:

\begin{lemma}[Extension of Claims 15 and 16 in \cite{wirth2019maximum}]\label{l.extension}
For any $D_2>0$, if $\beta$ is chosen small enough, then for any $n\geq 1$, if $\Lambda_n^2$ is equipped with Dirichlet boundary conditions, the bound~\eqref{e.NonN} (i.e. 
Claims 15 and 16 from \cite{wirth2019maximum}) is in fact valid for all test functions $f$ (without assuming the condition~\eqref{e.condition}). 
\end{lemma}

Equally importantly, we will observe (in Section \ref{ss.line} below) that if the charges $\rho\in \calN$ are all supported on the middle line $\Z = \Z\times \{0\}$, then the upper bound is still satisfied with a much smaller Dirichlet energy confined to the line, namely 
\begin{lemma}\label{l.line}
For any $D_2>0$, if $\beta$ is chosen small enough, the following holds. For all $n\geq 1$, if $\Lambda_n^2$ is equipped with Dirichlet boundary conditions, then for any test function $f$ and any 
ensemble $\calN$ made of charges $\rho$ which are all supported on the line $\Z=\Z \times \{0\}$, we have the following sharper upper bound 
\begin{align}\label{e.NonNline}
D_1 \sum_{\rho\in \calN} |z(\beta, \rho, \calN)| \cdot \<{\sigma,\rho}^2 \leq \frac \beta {D_2} \sum_{j \sim l  \in \Z \times \{0\}} (\sigma_j -\sigma_l)^2 
\end{align}
\end{lemma}

Now going back to~\eqref{e.123} and using Lemma \ref{l.line}, we obtain
\begin{align*}\label{}
\EFK{\beta,\Lambda_n^2,\lambda, v}{}{e^{\<{\phi,f}}}
& \geq 
e^{\frac 1 {2 \beta} \<{f, (- \Delta)^{-1} f}} \sum_{\calN\in \calF} c_\calN 
\exp\big[-\frac \beta {D_2} \sum_{j \sim l \in \Z \times \{0\}} (\sigma_j -\sigma_l)^2\big]
\end{align*}
 
We may now analyse the three cases of test functions $f$ listed in Theorem \ref{th.line2}:
\medskip

\ni
\textbf{Case 1a).}  Let us first see what happens if $f_\lambda(z) := \lambda\, 1_{z=(\floor{xn},\floor{yn})}$ when the point $w=(x,y)\in (-1,1)^2$, with $y\neq 0$.  Using~\eqref{e.Green1}, we get 

\begin{align}\label{e.corr1}
\EFK{\beta,\Lambda_n^2,\lambda, v}{}{e^{ \lambda \phi(\floor{ xn}, \floor{y n})}}
& \geq 
 e^{\frac {\lambda^2} {2 \beta} 
\frac {2 } {\pi }  \left( \log \frac n {r_D(x,y)} + c_0 +o(1) \right) 
 }  \\
 & \hskip 3cm \sum_{\calN\in \calF} c_\calN 
\exp\big[-\frac {\beta \lambda^2}  {D_2} \sum_{j \sim l \in \Z \times \{0\}} (\sigma_j -\sigma_l)^2\big]\,,
\end{align}
where
\begin{align*}\label{}
\sigma(\cdot)= \frac 1 \beta [-\Delta_{\Lambda_n^2}]^{-1}(1_{z=(\floor{xn},\floor{yn})})(\cdot) = \frac 1 \beta G_{\Lambda_n^2}(\floor{ xn}, \floor{y n}, \cdot)
\end{align*}
\ni
(Recall that with our convention in~\eqref{e.gff8}, $\Delta$ is the ``probabilistic'' Laplacian here, i.e its inverse is the Green function of the SRW). 

The Green function $G_{\Lambda_n^2}(\floor{ xn}, \floor{y n}, \cdot)$ has a ``log''-singularity at the point $(\floor{ xn}, \floor{y n}, \cdot)$. One then expects that its gradient along edges should decay as one over the distance to the singularity. This is indeed the case and we claim the following upper bound: 
For any $i\sim j \in \Z \times \{0\}$, 
\begin{align}\label{e.ClaimG1}
|G_{\Lambda_n^2}(\floor{ x n}, \floor{y n}, i) - G_{\Lambda_n^2}(\floor{ xn}, \floor{y n}, j)| \leq O(1) \frac 1 {\dist((i,0), (\floor{ xn}, \floor{y n}) +1}
\end{align}
 
In particular, we see that if the point $(x,y)$ is at macroscopic distance from the middle line (i.e. $y\neq 0$), then  at each point $(i,0)\in \Z\times \{0\}$, the gradient of the Green functions is upper bounded by $O(1/n)$. 
For pairs of points $i\sim j$ far from the boundary, this claim follows  from Lemma 6.3.3 in \cite{lawler2010random}. For general pairs of points, it follows from the same techniques as the one used for Proposition \ref{pr.Grad2} which handles the proximity of a more difficult type of  boundary than in the present case. Since this concerns the SRW on $\Z^2$ we do not provide more details here. 

This implies that the correction to the Gaussian behaviour in~\eqref{e.corr1} is at most 
\begin{align*}\label{}
& \sum_{\calN\in \calF} c_\calN 
\exp\big[-\frac {\beta \lambda^2}  {D_2} \sum_{j \sim l \in \Z \times \{0\}} (\sigma_j -\sigma_l)^2\big] \\
& \leq  \sum_{\calN\in \calF} c_\calN 
\exp\big[-\frac { \lambda^2}  {D_2 \beta}  
n*O(\frac 1 {n^2}) \big] \\
& \leq \sum_{\calN\in \calF} c_\calN (1- O(n^{-1}) \leq 1-O(n^{-1})\,,
\end{align*}
since by construction $\sum_{\calN \in \calG} c_{\calN} =1$. 

This implies
\begin{align*}\label{}
\EFK{\beta,\Lambda_n^2,\lambda, v}{}{e^{ \lambda \phi(\floor{ xn}, \floor{y n})}}
& \geq  e^{\frac {\lambda^2} {2 \beta} 
\frac {2 } {\pi }  \left( \log \frac n {r_D(x,y)} + c_0 +o(1) \right) 
 } (1-O(\frac {\lambda^2} {\beta n}))\,.
 \end{align*} 
By letting $\lambda \to 0$ we obtain the desired asymptotic on the variance. (Note that the $o(1)$ is a $o(1)$ as $n\to \infty$ and is not $\lambda$-dependent).  
\medskip

\ni
\textbf{Case 1b).} For general points $(x,y)\in D$, we cannot hope to keep such a precise estimate. Indeed, if $y=0$ and we are measuring the fluctuations at a point on the middle line, say at $(\floor{xn},0)$, it is no longer the case that the Dirichlet energy of $\sigma=\frac 1 \beta [-\Delta]^{-1} f$ has a vanishing contribution coming from the middle line. Instead it is easy to see by the above claim~\eqref{e.ClaimG1} that it  has a non-vanishing contribution of order $O(1)$. This is why we do not expect that the conformal radius correction from~\eqref{e.fine} will still be accurate for points sitting on the line. As long as $(x,y)$ in in the open set $D=(-1,1)^2$ (i.e. middle line or not), the convergence in law~\eqref{e.lessfine} follows easily from the convergence of the Laplace transform of $\frac 1 {\sqrt{\log n}} 1_{z=(\floor{xn},\floor{yn})}$.

\medskip

\ni
\textbf{Case 2).} Fix a  test function $g\in \calC_c^1([-1,1]^2)$. We wish to show that 
\begin{align*}\label{}
\frac 1 {n^2} \sum_{z\in \Lambda_n^2} \Psi_n^{line}(z) g\big(\frac z n\big)  \overset{(d)}\longrightarrow 
\calN(0, \frac{4}{ \beta} \<{g,(-\Delta)^{-1} g})\,, 
\end{align*}
where $(-\Delta)^{-1}$ is the continuous Green function on $D=(-1,1)^2$ with Dirichlet boundary conditions. 
Let us apply once again the Coulomb gas expansion lower-bound~\eqref{e.123} applied to the test function 
\begin{align*}\label{}
f:  \, z\in \Lambda_n^2 \mapsto  \frac 1 {n^2} g(\frac z n)
\end{align*}
This gives us for any $\lambda \in \R$, 
\begin{align*}\label{}
\EFK{\beta,\Lambda_n^2,\lambda, v}{}{e^{\lambda \<{\phi,f}}}
& = 
\EFK{\beta,\Lambda_n^2,\lambda, v}{}{e^{\lambda \frac 1 {n^2} \sum_{z\in \Lambda_n^2} \phi(z) g\big(\frac z n\big)
}
}\\
& \geq  e^{\frac {\lambda^2} {2 \beta} \<{f, (- \Delta_{\Lambda_n^2})^{-1} f}} \sum_{\calN\in \calF} c_\calN  \exp\big[-D_1 \lambda^2  \sum_{\rho\in \calN} |z(\beta,\rho,\calN)|\<{\sigma,\rho}^2\big] \,,
\end{align*}
where as previously $\sigma = \frac 1 \beta [-\Delta_{\Lambda_n^2}]^{-1} f$. With the above choice of test function $f$ (i.e. normalized by the volume), and with our choice of normalisation for the discrete Laplacian one has 
\begin{align*}\label{}
\<{f, (- \Delta_{\Lambda_n^2})^{-1} f} \underset{n\to \infty}\longrightarrow  4 \<{g, (-\Delta_D)^{-1} g}\,,
\end{align*}
which is the reason behind the factor $4$ in~\eqref{e.Line}.

Furthermore, since the charges are still confined to the line $\Z\times \{0\}$, we again rely on  Lemma \ref{l.line}. As in the above two cases, it is thus enough to upper-bound the Dirichlet energy of $\sigma$ localised on the line $\Z \times \{0\}$. Namely,
\begin{align*}\label{}
& \sum_{j\sim l \in \Z \times \{0\}} (\sigma_j - \sigma_l)^2 \\
& = \frac 1 {\beta^2}  \sum_{j\sim l \in \Z \times \{0\}} \left(  [-\Delta_{\Lambda_n^2}]^{-1} f(j) - 
 [-\Delta_{\Lambda_n^2}]^{-1} f(l) \right)^2 \\
& =  \frac 1 {\beta^2}   
\sum_{j\sim l \in \Z \times \{0\}} \left(  \sum_{z\in \Lambda_n^2} G_{\Lambda_n^2}(j,z)f(z) - \sum_{z\in \Lambda_n^2} G_{\Lambda_n^2}(l,z)f(z) \right)^2\,.
\end{align*}
The above sum looks like a discrete integration by parts, but it is not quite so (the discrete gradient is on $j$ instead of $z$). Note that by shifting the domain $\Lambda_n^2$ by one, one may recover a standard discrete integration by parts to the cost of an additional error boundary term of order $O(\frac 1 n)$ when the root is far from the boundaries.

We will instead follow a more direct approach. Let us rewrite the above identity (using the definition of $f$) as 
\begin{align}\label{e.FRACT}
& \sum_{j\sim l \in \Z \times \{0\}} (\sigma_j - \sigma_l)^2 \nn \\
& =  \frac 1 {\beta^2}   
\sum_{j\sim l \in \Z \times \{0\}} \left(  \sum_{z\in \Lambda_n^2} [G_{\Lambda_n^2}(j,z) - G_{\Lambda_n^2}(l,z)] \frac 1 {n^2} g(\frac z n) \right)^2 \nn \\
& \leq \frac 1 {\beta^2}
\sum_{j\sim l \in \Z \times \{0\}}  
\left(
\sum_{z\in \Lambda_n^2} 
|G_{\Lambda_n^2}(j,z) - G_{\Lambda_n^2}(l,z)| 
\cdot 
\frac 1 {n^2} \|g\|_\infty  \right)^2 \nn \\
& 
\leq \frac 1 {\beta^2}
\sum_{j\sim l \in \Z \times \{0\}}  
O(1) \frac {\|g\|_\infty^2}{n^4} \;\;   \left(O(1) \sum_{r=1}^n  r \frac 1 r \right)^2  \\
& \leq \frac {O(1)} {\beta^2}   \frac{\| g \|_\infty^2} {n}   
\end{align}
We used the fact that if $j\sim l$, then $|G_{\Lambda_n^2}(j,z) - G_{\Lambda_n^2}(l,z)| \leq O(1) [\dist(j,z)+1]^{-1}$. See Claim~\eqref{e.ClaimG1}. Note that when $j\sim l$ is close to the boundary, some further care is needed as the distance to the boundary also matters. See again Section \ref{s.GRAD} where the effect of the boundary  being close is taken into account in a more difficult setting.   
Also, in the last inequality we used that  $\sharp \{ i\sim j \in \Z\times \{0\} \}$ is $O(n)$. 

We thus obtain this way a quantitative bound on the speed of convergence to the limiting Gaussian random variable $\calN(0, \frac 4 \beta \<{g, (-\Delta)^{-1} g})$  with sharp optimal speed $O(\frac {1} n)$. 
\qed

\subsection{Invisible subsets of $\Lambda_n^2$.}\label{ss.fractal}
$ $

\begin{figure}[!htp]
\begin{center}
\includegraphics[width=0.7\textwidth]{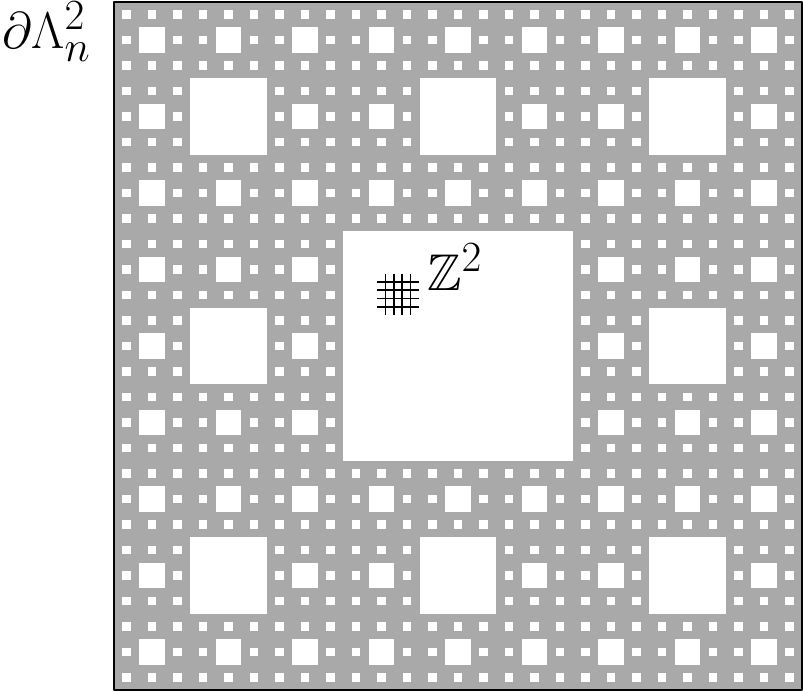}
\end{center}
\caption{The Monge discrete fractal $F_k$ in the $2d$ box $\Lambda_{3^k}^2$.  If the discrete Gaussian free field $\varphi_{3^k}$ on this domain is conditioned to take integer values on each site $x$ of $F_k$, then when the temperature is high enough, the effect of this conditioning is invisible at large scales. }\label{f.fractal}
\end{figure}

In this subsection, we shortly explain how the above quantitative speed of convergence can be extended to prove the invisibility of integers to larger subsets of $\Lambda_n^2$ on which the field will be conditioned to be integer-valued. We will not try to characterize all possible such subsets but will discuss the following two cases:
\bi
\item[a)] An horizontal strip $\Lambda_n\times [0,B_n]$ of width $B_n=o(n )$
\item[b)] A {\em Monge-type} fractal set $F_k \subset \Lambda_{n=3^k}^2$.  This well-known fractal set is built in a recursive way as shown in the self-explanatory Figure \ref{f.fractal}. It spans a polynomial fraction of the volume of $\Lambda_{3^k}^2$. 
\ei

We claim that the same techniques as in the above Sections imply the following statement.
\begin{theorem}\label{th.fractal}
If $\beta$ is small enough, and if for any $k\geq 1$, $\Psi_{3^k}^{Monge}$ denotes the rescaled Dirichlet Gaussian free field on $\frac 1 {3^k}\Lambda_{3^k}^2 \cap [-1,1]^2$ conditioned to take integer values on the (rescaled) set $\frac 1 {3^k} F_k \subset [-1,1]^2$,
then as $k\to \infty$, this field converges to the $\beta$-GFF on $[-1,1]^2$ with Dirichlet boundary conditions and with same inverse temperature $\beta$.  

The same Holds in the case of a band  $\Lambda_n\times [0,B_n]$ as long as the width $B_n=o(n)$. 
\end{theorem}

We only point out the two main points which require some attention in the case of the Monge set $F_k$ (the band is a more direct consequence):
\bi
\item On important hidden point is that we need the paths $\gamma_{x,y}$ in the analysis of $\<{\rho, \sigma}$ (see for example~\eqref{e.1516}) to be roughly of same diameter as the Eucilidean distance $\|x-y\|_2$. 
This requires the fractal sets we may consider  to be well behaved from this point of view and this is the case with the Monge fractal set $F_k$. 
\item As in the previous sections, the key analytical fact is that one needs to show that the Dirichlet energy of Green functions supported here either by the band or the Monge fractal $F_k$ are negligible with respect to the Green function itself. 
In the case of the fractal set $F_k$, let us observe that if one considers the Green function rooted in the center of the box, then the contribution of the points at macroscopic distance $\Omega(3^k)$ from this root will be upper bounded by 
\begin{align*}\label{}
(3^{k})^{d_H} \times \frac 1 {(3^k)^2}\,,
\end{align*}
where $d_H$ stands for the Hausdorff dimension of the discrete set $F_k$ and where the second term (in distance$^{-2}$) arises thanks to the gradient estimate~\eqref{e.ClaimG1}. 

By using the same observation on each dyadic scale (or rather $3^m$ scales), and arguing as we did in the subsections above when testing against smooth functions, we obtain the desired result. 
\ei

\subsection{Free and Periodic boundary conditions.}\label{ss.torus}

We now discuss the case of periodic and free boundary conditions here. One significant advantage of these boundary conditions is that one does not need to handle non-neutral charges. In particular we do not need to extend Wirth results \cite{wirth2019maximum}.
Another advantage is the fact such boundary conditions are more convenient for a Renormalization Group flow approach (as in \cite{bauerschmidt2022discrete,bauerschmidt2022discreteN2}). 

Let us then consider the domain $\Lambda_n^2=\{-n,\ldots,n\}^2$ with  {\em periodic} or {\em free} boundary conditions. In the later case,  we will write the periodic domain $\T_n^2$ instead of $\Lambda_n^2$ plus identifications and we will denote by $\T_n\times \{0\}$ its middle circle. 


We may now state our main result in this setting:
\begin{theorem}\label{th.FP}
If $\beta$ is low enough and if one conditions a $\beta$-GFF on $\Lambda_n^2$ with free boundary conditions (resp. $\T_n^2$ with periodic boundary conditions)  to take its values in the integers along the middle line $\Lambda_n \times \{0\}$  (resp.  middle circle $\T_n \times \{0\}$)), then \uline{\em integers are invisible} in the sense that this conditioned field has the same global fluctuations as the GFF. Indeed, if $\Psi_n^{line}$ denotes such a field. Then for any test function $g\in C^1([-1,1]^2)$ (resp. $C^1(\T^2)$), such that $\int g(z)dz =0$,  if $g_n$ is the discretisation of $g$ which assigns to any $z\in \Lambda_n^2$ (resp. $\T_n^2$) the averaged value $n^2 \int_{\frac 1 n  (z +  [0,1]^2)} g(u) du$ (so that $\sum_{z} g_n(z)=0$), one has  
\begin{align*}\label{}
\EFK{\beta}{line}{e^{\sum_{z\in \Lambda_n^2} \frac 1 {n^2}\Psi_n^{line}(z) g_n(z)}}
\underset{n\to \infty}\longrightarrow  e^{-\frac 4 \beta \<{g,(-\Delta)^{-1} g}}\,,
\end{align*}
where $\Delta$ is the continuous Laplacian on $[-1,1]^2$ (resp. the torus $\T^2$) with Neumann boundary conditions (resp. periodic boundary conditions).
\end{theorem}

\ni
{\em Proof of Theorem \ref{th.FP}.}
Let us give the details of the proof in the case of the periodic boundary conditions (free boundary conditions are treated in a similar way).
Following \cite{RonFS} it is convenient to pick any vertex $v\in \T_n^2$ and to root the field $\Psi_n^{line}$ to be zero at this vertex. (In fact, rather than strongly rooting the field at $v$, following \cite{RonFS}, we impose that $\phi_v\in[-\pi, \pi)$ and we apply a Sine-Gordon trigonometric potential at $v$). 

Note that such a rooting procedure is harmless as one is only integrating against test functions $g_n$ of zero average. 

For our present application, we are in fact required to pick the vertex $v$ \uline{on the middle circle} $\T_n \times \{0\}$. This is not immediate to see why at first sight. This is due to the fact that the root vertex may  carry non-zero Coulomb charges! Because of this, if $v$ is far from the line, it may violate the main estimate of Lemma \ref{l.line} which controls the error using the Dirichlet energy supported by the line. Though the choice of root vertex eventually does not affect fluctuations (by translation invariance),  we cannot use translation invariance to dissociate the choice of root from the choice of line on which we measure the Dirichlet energy.
  One has to be careful here as the Dirichlet energy on a line intersecting the root is  higher than for a line at  far distance from the root.  (See Subsection \ref{ss.line}  which is specific to free boundary conditions).

Let us then choose $v=v_0$ to be the origin.

Once this root is chosen, we claim that for the $0$-mean test function $g_n$ (the discretisation of $g$ defined in Theorem \ref{th.FP}), one has the identity
\begin{align*}\label{}
\<{g_n, (-\Delta_{\T_n^2})^{-1} g_n} = \<{g_n, (-\Delta_{\T_n^2}^{v_0})^{-1} g_n}\,,
\end{align*}
where the first Laplacian is invertible on zero mean functions while the second one is invertible on all functions $f : \T_n^2 \setminus \{ v_0\}$ and where the Laplacian $\Delta^{v_0}_{\T_n^2}$ has Dirichlet boundary conditions at the single point $v_0$. 
This follows for example by viewing the same  Gaussian vector in two different ways. This follows also from \cite[Section A.2]{GS2} where such rerooting procedures have already been extensively used.

Instead of choosing the function $\sigma$ to be $\frac 1 \beta [-\Delta_{\T_n^2}] f$, we instead choose the viewpoint 
\begin{align*}\label{}
\sigma:= \frac 1 \beta [-\Delta_{\T_n^2}^{v_0}] f
\end{align*}
(In particular, $\sigma$ is not a zero-mean function), where $f$ is as before the following test function:
\begin{align*}\label{}
f:  \, z\in \T_n^2 \mapsto  \frac 1 {n^2} g_n(z)
\end{align*}

It turns out that Lemma \ref{l.line} still holds in he context of free/periodic boundary conditions as it is explained in Subsection \ref{ss.line}. In particular, by running the same analysis as in the Dirichlet case, we end up controlling 
\begin{align*}\label{}
& \sum_{j\sim l \in \T_n \times \{0\}} (\sigma_j - \sigma_l)^2 \\
& = \frac 1 {\beta^2}  \sum_{j\sim l \in \T_n \times \{0\}} \left(  [-\Delta^{v_0}_{\T_n^2}]^{-1} f(j) - 
 [-\Delta^{v_0}_{\T_n^2}]^{-1} f(l) \right)^2 \\
& =  \frac 1 {\beta^2}   
\sum_{j\sim l \in \T_n \times \{0\}} \left(  \sum_{z\in \T_n^2} G^{v_0}_{\T_n^2}(j,z)f(z) - \sum_{z\in \T_n^2} G^{v_0}_{\T_n^2}(l,z)f(z) \right)^2\,.
\end{align*}

We now claim that the gradient of the Green function of the random walk on the torus  $\T_n^2$ killed when first hitting $v_0$ is now decaying as follows:
\begin{align}\label{e.GRADe}
|G^{v_0}_{\T_n^2}(j,z) -  G^{v_0}_{\T_n^2}(l,z)| \leq 
\begin{cases}
& O(1)  \frac {\log(\dist(v_0, z))} { \dist(j,z)+1 }  \,\, \text{ \;\; if  } \dist(j,z) <  \tfrac 1 2  \dist(j,v_0) \\
& O(1)  \frac {\log(\dist(v_0, z))} { \dist(j,v_0)+1 }    \,\, 
\text{ \;\; if  } \dist(j,z) \geq  \tfrac 1 2  \dist(j,v_0) 
\end{cases}
\end{align}
This estimate follows from the same type of coupling techniques as our gradient estimates in Section \ref{s.GRAD}: the ratio $\frac 1 {\dist(j,z)+1 }$ comes from the low probability for the walks starting at $j\sim l$ not to couple before reaching $z$ or $v_0$. The numerator follows using the fact the rooted Green function $G_{\T_n^2}^{v_0}(v,v)$ with periodic or free boundary conditions is upper bounded by the log of the distance to the root, i.e. $O(1) \log (\dist(v,v_0))$.  (See for example \cite[Proposition 2.6]{GS2}). We shall not give further details  here.

\textbf{Technical difficulty:}
It turns out that if one applies this bound readily, we will not be able to show that the effect of the line is negligible. 
Indeed the edges $j \sim l$ which are close to the chosen root $v_0$ produce gradients 
\begin{align*}\label{}
G^{v_0}_{\T_n^2}(j,z) - G^{v_0}_{\T_n^2}(l,z)
\end{align*}
which do not converge to 0 as $\dist(z,v_0)$ goes to infinity.  The fundamental reason why we may still conclude is due to the neutrality of the test function $f$. Even though the gradient $G^{v_0}_{\T_n^2}(j,z) - G^{v_0}_{\T_n^2}(l,z)$ does not asymptotically vanish, it does converge to a constant as $|z| \to \infty$ and we may then rely on the neutrality of $f$. 
\smallskip

We will use the neutrality in a more indirect but sharper way: It turns out that the gradient field 
\begin{align*}\label{}
\{ [-\Delta^{v_0}_{\T_n^2}]^{-1} f(j) - 
 [-\Delta^{v_0}_{\T_n^2}]^{-1} f(l) \}_{j\sim l \text{ in } \T_n^2}
\end{align*}
does not depend on the choice of root. (This is also a consequence of  \cite[Section A.2]{GS2}).  I.e. for any $v_0,v_1 \in \T_n^2$, and any $j\sim l$, 
\begin{align*}\label{}
 [-\Delta^{v_0}_{\T_n^2}]^{-1} f(j) - 
 [-\Delta^{v_0}_{\T_n^2}]^{-1} f(l) = 
 [-\Delta^{v_1}_{\T_n^2}]^{-1} f(j) - 
 [-\Delta^{v_1}_{\T_n^2}]^{-1} f(l)\,.
\end{align*}

Let us then change of root at this stage only (up to now in the analysis it was important that $v_0$ belongs to the line $\T_n \times \{0\}$). Let us fix $v_1$ any root which is at macroscopic distance $\Omega(n)$ from the circle $\T_n \times \{0\}$.
We now proceed similarly as in ~\eqref{e.FRACT}, using furthermore (a) the above invariance of gradients over the choice of root and (b) the estimate~\eqref{e.GRADe} applied to the root $v_1$. This gives us

\begin{align*}\label{}
\sum_{j\sim l \in \text{ line }} (\sigma_j - \sigma_l)^2 
 & =  \frac 1 {\beta^2}   
\sum_{j\sim l \in \T_n \times \{0\}} \left(  \sum_{z\in \T_n^2} [G^{v_0}_{\T_n^2}(j,z) - G^{v_0}_{\T_n^2}(l,z)] f(z)  \right)^2 \nn \\ 
 & \overset{(a)}=  \frac 1 {\beta^2}   
\sum_{j\sim l \in \T_n \times \{0\}} \left(  \sum_{z\in \T_n^2} [G^{v_1}_{\T_n^2}(j,z) - G^{v_1}_{\T_n^2}(l,z)] f(z) \right)^2 \nn \\   
& \leq \frac {\|g\|_\infty^2 } {n^4 \beta^2} \;\; \sum_{j\sim l \in \T_n \times \{0\} }
\left(
\sum_{z\in \T_n^2} 
|G^{v_1}_{\T_n^2}(j,z) - G^{v_1}_{\T_n^2}(l,z)|
\right)^2
 \nn \\      
& 
\overset{(b)}\leq \frac {\|g\|_\infty^2 } {n^4 \beta^2}
\;
\sum_{j\sim l \in \T_n \times \{0\}}
\left(  O(1)  \Big(\sum_{k=1}^n  k \frac 1 k \Big) \log(n)   +  O(1) n^2 \frac {1} n  \log(n) 
\right)^2
 \nn \\   
& 
\leq O(1) \frac {\|g\|_\infty^2 } { \beta^2}   \frac {\log(n)^2}{n}
 \nn 
\end{align*}
which thus concludes the proof (with a quantitative speed of $(\log n)^2 /n$).


\qed

\subsection{Fluctuations of the discrete Gaussian chain with $\alpha_c=2$.}

If one now considers the discrete Gaussian Chain as defined in Definition \ref{d.dGC} at $\alpha_c=2$ and with explicit coupling constants $J_2(r)$ given in~\eqref{e.J2},  Lemma \ref{l.gffs0} shows that we may realise this discrete Gaussian chain by conditioning the Gaussian free field $\varphi_n$ (see~\eqref{e.gff8}) on the diamond graph $\calD$ and rooted on the slit $L_n$ to take its values in the integers on the sites of $\Lambda_n \times \{0\}$. As in the previous sections, this then corresponds to affecting Coulomb charges only to the sites on the line $\Lambda_n \times \{0\}$.   
See Figure \ref{f.diamond}.

We will give a more detailed proof in the general case $\alpha \in [2,3)$ in Section \ref{ss.TheProof}. Let us now list the main steps in the present case $\alpha=\alpha_c=2$.

\bi
\item[(1)] We first need to check that Fröhlich-Spencer expansion as well as Wirth's extension to Dirichlet boundary conditions works on the setting of the slit domain $D_n:=\calD \setminus L_n$ with Dirichlet boundary conditions on $\p D_n = L_n$. 
This will be discussed in Section \ref{ss.TheProof}. 
\ei

At this stage notice that, modulo the above (minor) extension,  the $\log n$ delocalisation of the discrete Gaussian Chain at $\alpha_c=2$ (and high temperature) is already a Corollary of \cite{FS, wirth2019maximum} (due to Dirichlet boundary conditions) as we need to control the behaviour of ``fewer'' charges. 
Note that this is not the path followed in \cite{kjaer1982discrete, frohlich1991phase}.

Now, for the identification of an invariance principle and its effective temperature, we need to investigate more quantitatively what the Coulomb-gas decomposition in \cite{FS} gives us when charges are bound to the line $\Lambda_n \times \{0\}$ as we did in the above $2s$ cases.

\bi
\item[(2)] For this, we then need to extend Wirth's control of non-neutral charges all the way to the boundary as we explain (for other boundary conditions) in the next section. See also Section \ref{ss.TheProof}
\item[(3)] We notice here that the weights $c_{i,j}$ used when controlling the terms $\<{\sigma, \rho}$ (for example in expressions such as~\eqref{e.NonN}) do not need to correspond to edges of the lattice. This is not hard to see but is important in the present setting as we work on the rotated lattice $\calD$. 
\item[(4)] When $\alpha\in (2,3)$ (especially $\alpha\in [\tfrac 5 2, 3)$, see Remark \ref{r.tedious2}) as we will see in Section \ref{s.proof}, we will need to work in a smoother domain as otherwise our estimates on the gradients of Green functions would degenerate. See Remark \ref{r.tedious2}. The situation in the present case $\alpha_2=2$ is easier from this point of view. 
\ei

%

We will explain the details of this proof together with the case $s\in (2,3)$ in Section \ref{ss.TheProof}.

\section{Dirichlet boundary conditions, non-neutral charges and Dirichlet energy supported on a line}\label{s.EXT}

The purpose of this Section is to provide the main technical results we used throughout Section \ref{s.proof2} in order to show that {\em integers are invisible} when $\alpha=\alpha_c=2$ and the temperature is high enough. The next two subsections will focus on the following properties:
\bnum
\item We will first explain in the easier case of the free Boundary that the content of Lemma \ref{l.line} (for Dirichlet boundary conditions) still holds for free/periodic boundary conditions. The purpose of this Lemma and its free/periodic version is in some sense to control the distance to  Gaussian behaviour for our fields conditioned to be integer-valued, in terms of a Dirichlet energy confined to a line. 

We will first need to briefly sketch how the upper bound~\eqref{e.1516} is derived in \cite{FS,RonFS}. 

\item Then, we will explain following Wirth \cite{wirth2019maximum} how to extend Fröhlich-Spencer analysis to the case of Dirichlet boundary conditions.
Compared to \cite{wirth2019maximum}, we will significantly extend the possible test functions which can be analyzed.
This corresponds to a proof of Lemma \ref{l.extension} followed by the case where neutral and non-neutral charges are restricted to a line, i.e. a proof of Lemma \ref{l.line}.

\enum

\subsection{Free boundary conditions: effective temperature and Dirichlet energy confined to the line.}\label{ss.line}
$ $ 

The place where one can read off the discrepancy between the fluctuations of the Gaussian free field and its $\Z$ conditioning  version is best seen in Fröhlich-Spencer in the inequality~\eqref{e.nice1}. 

Fröhlich-Spencer manage to show that this difference is small  (yet of same order by \cite{GS2}) compared to the Gaussian Free Field fluctuations thanks to the comparison of the bound in~\eqref{e.nice1} with the Dirichlet energy of $\sigma$ in~\eqref{e.1516}.

Let us first explain how Fröhlich-Spencer prove the key estimate~\eqref{e.1516} from the expression~\eqref{e.nice1} for free boundary conditions.  We will follow \cite{RonFS} here. And we will then explain, still in the case of free boundary conditions, how to adapt this proof to our case where charges are restricted to a line. 

\subsubsection{Dirichlet energy upper bound on the error term in \cite{FS}.}\label{sss.dir}

The idea to prove~\eqref{e.1516} is for each \underline{neutral} charge $\rho \in \calN$ to rewrite $\<{\sigma, \rho}$ as a sum of nearest-neighbor gradients. For example if on a $1D$ line $\rho= \delta_{10}+\delta_5-2\delta_0$, then $\<{\sigma, \rho}= \sum_{k=0}^4 2(\sigma_{k+1}-\sigma_k) + \sum_{k=5}^9 \sigma_{k+1}-\sigma_k$.  

As explained in \cite{FS,RonFS,wirth2019maximum}, for a general neutral charge $\rho \in \calN$, one can rewrite
\begin{align}\label{e.rewrite}
\<{\sigma,\rho} = \sum_{i\sim j \in D(\rho)} c_{i,j}( \sigma(i)-\sigma(j) )\,,
\end{align}
where the coefficients $c_{i,j}\in \Z$ satisfy $|c_{i,j}| \leq \frac 1 2 \|\rho\|_1 \leq \frac 1 2 \|\rho\|_2^2$, where $D(\rho)$ is a square domain surrounding $\rho$ and with a diameter comparable to the diameter $d(\rho)$ of the support of $\rho$ (see \cite{RonFS}). The proof is straightforward: it is enough to  decompose $\rho$ as an arbitrary union of dipoles $\delta_x - \delta_y$ where $x,y\in \mathrm{Supp}(\rho)\subset D(\rho)$. Now for each such dipole, one can choose a directed path $\gamma_{x,y} \subset D(\rho)$ of nearest-neighbour edges which is connecting $x$ to $y$. As such for each $i\sim j$, $c_{i,j}$ is the number of oriented path using $(i,j)$ minus the number of oriented paths using $(j,i)$. Using this viewpoint, the bound $|c_{i,j}| \leq \frac 1 2 \|\rho\|_1$ is clear. 

Then, one applies Cauchy-Schwarz to get 
\begin{align*}\label{}
|\<{\sigma,\rho}|^2 \leq  |D(\rho)| \cdot  \| \rho\|_2^4  \sum_{i\sim j \in D(\rho)} (\sigma_i - \sigma_j)^2\,.
\end{align*}
At this stage, the combinatorial properties of the different charges $\rho$ which constitute an ensemble $\calN$ have not been used (besides the neutrality of each $\rho$). 
The next step is to upper bound 
\begin{align*}\label{}
\sum_{\rho\in \calN}  |z(\beta,\rho,\calN)|\<{\sigma,\rho}^2  
& \leq \sum_{\rho \in \calN}  |z(\beta,\rho,\calN)|\,  |D(\rho)| \cdot  \| \rho\|_2^4  \sum_{i\sim j \in D(\rho)} (\sigma_i - \sigma_j)^2
\end{align*}

Notice that the term $|D(\rho)| \|\rho\|_2^4$ may be very large especially for a charge $\rho$ with large diameter. 
In order to obtain~\eqref{e.1516}, the idea is to use the following quantitative upper bound on the coefficients $z(\beta,\rho,\calN)$ which follows from complex translation along spin-waves: i.e. for $\beta$ sufficiently small, 
\begin{align*}\label{}
|z(\beta,\rho,\calN)| \leq \exp \left[ -\frac c \beta  \left( \|\rho\|_2^2 + \log_2 |D(\rho)| \right)\right]\,.
\end{align*}
We see that for any fixed $\rho \in \calN$, this small coefficient controls the diverging term $|D(\rho)| \|\rho\|_2^4$.

One additional difficulty here is that for any $\rho,\rho'$ one has by construction $\mathrm{Supp}(\rho)\cap \mathrm{Supp}(\rho') = \emptyset$ but this does not imply that the square domains $D(\rho)$ and $D(\rho')$ are distinct. In fact, typically, many such domains do overlap (and up to $\frac 1 \alpha \log n$ such squares may overlap on top of each other in a domain $\Lambda_n^2$, where $\alpha\in(\frac 3 2 ,2)$ is used in the construction of the expansion of charges).  

To deal with this superposition effect, the idea is to use another specificity of the expansion into charges : for each dyadic scale $2^k$,  distinct charges $\rho_1,\rho_2$ with diameters in $[2^{k-1},2^k]$ cannot overlap. To conclude the proof, it is thus sufficient to group the charges depending on their dyadic scale $k$ and to notice that the sum over dyadic scales $k$ is small when $\beta$ is sufficiently small. See \cite[Section 3]{RonFS} for details.  
 
\subsubsection{The case of charges restricted to a line.}\label{ss.linefree}

Let us consider here the setting of Subsection \ref{ss.torus} where we considered the GFF on a the torus $\T_n^2$ conditioned to take integer values on the circle $\T_n \times \{0\} \subset \T_n^2$. We denote this non-Gaussian field as $\Psi_n^{line}$. 

As mentioned in Section \ref{ss.torus}, we decide to fix the root $v_0$ in this line, say at the origin. 
We choose $N$ large ($N$ is sent to infinity at the end of the proof as in \cite[Section 5]{RonFS}) and we assign to each vertex $i$ in the line $\T_n \times \{0\}$ the trigonometric polynomial 
\begin{align*}\label{}
\lambda_i(\phi(i)) = 1 + 2 \sum_{q=1}^N  \cos(q \phi(i))\,,
\end{align*}
while at each vertex $i$ away from the line, we assign the trivial polynomial $\lambda_i(\phi(i)):=1$. 
With this choice of family of trigonometric polynomials, $\lambda^{line}_{N} := (\lambda_i)_{i\in \T_n^2}$, the same proof as in Fröhlich-Spencer applies. It leads us to a Coulomb-gas expansion where charges are restricted to the line and with same quantitative bounds on the coefficients 
$z(\beta,\rho,\calN)$. The reason here is that, charges restricted to the line are a specific case of the charges handled in \cite{FS}. In particular both the splitting/merging algorithm for charges and the spin-wave analytic part are exactly the same. 

We thus obtain for any test function $f$ of vanishing mean the analog of the identity~\eqref{e.123}
\begin{align}\label{e.1234}
\EFK{\beta,\T_n^2,\lambda^{line}_N, v_0}{}{e^{\<{\phi,f}}}
\geq  e^{\frac 1 {2 \beta} \<{f, (- \Delta)^{-1} f}} \sum_{\calN\in \calF} c_\calN  \exp\big[-D_1  \sum_{\rho\in \calN} |z(\beta,\rho,\calN)|\<{\sigma,\rho}^2\big] \,.
\end{align}

For any ensemble $\calN$, we have that each $\rho\in \calN$ is a neutral charge supported on the line $\T_n\times \{0\}$. 
We may thus adapt the above argument and write 
\begin{align*}\label{}
\<{\sigma, \rho}= \sum_{i\sim j \in D(\rho) \cap \T_n \times \{0\}} c_{i,j} (\sigma(i) -\sigma(j))\,. 
\end{align*}
The difference with~\eqref{e.rewrite} is that we are only summing over the nearest neighbour points along the line. (Note that $D(\rho)$ is still a $2D$ square surrounding the charge $\rho$ as the spin-waves used in the proof are still two-dimensional). This is because for any $x,y \in \mathrm{Supp}(\rho)$, we can find a path $\gamma_{x,y}$ which stays inside the line $\T_n \times \{0\}$. {\em (Note that this would be wrong if we had not chosen the root $v_0=0$ inside the line!)}. 
\medskip

Given this, the rest of the proof (which groups charges depending on their dyadic scales) still holds and we obtain that for $\beta$ small enough,

\begin{align*}\label{}
\EFK{\beta,\T_n^2,\lambda_N^{line}, v_0}{}{e^{\<{\phi,f}}}
& \geq 
e^{\frac 1 {2 \beta} \<{f, (- \Delta)^{-1} f}} \sum_{\calN\in \calF} c_\calN 
\exp\big[-\frac \beta {D_2} \sum_{j \sim l \in \T_n \times \{0\}} (\sigma_j -\sigma_l)^2\big]\,.
\end{align*}
This is the analog of the conclusion of Lemma \ref{l.extension} except it holds here for periodic boundary conditions. 
\medskip

\begin{remark}\label{r.calD}
Note that it is not needed in this proof that for each $x,y \in \rho$, the chosen path $\gamma_{x,y}$ is using nearest-neighbor points. In particular the paths $\gamma_{x,y}$ may go from $x$ to $y$ by making jumps of distance $2$. This remark will be relevant when dealing with the discrete Gaussian chain at $\alpha=2$ with explicit coupling constants~\eqref{e.J2} arising from the rotated diamond graph $\calD$. 
\end{remark}

\subsection{Dirichlet boundary and handling non-neutral charges.}\label{ss.NN}
$ $


In this Subsection, we start by briefly reviewing how the work \cite{wirth2019maximum} managed to handle non-neutral Coulomb charges. Then we prove our extension Lemma \ref{l.extension} as well as the main technical Lemma \ref{l.line}.

\subsubsection{The analysis of non neutral charges in \cite{wirth2019maximum}.}

In \cite{wirth2019maximum}, Wirth managed to extend the Coulomb gas expansion from Fröhlich-Spencer \cite{FS} in order to cover the case of Dirichlet boundary conditions. The proof idea of such an extension was highlighted in \cite[Appendix D]{wirth2019maximum}, but a full rigorous treatment of this expansion was missing until \cite{wirth2019maximum}. The main difficulty compared to the argument we sketched above in Section \ref{ss.sketchFS} is the presence of non-neutral Coulomb charges which also contribute to the Laplace transforms.

We may thus proceed exactly as in Section \ref{ss.sketchFS}, except the entire boundary $\p \Lambda$ is now rooted at 0. Wirth adapts the expansion into Coulomb charges in such a way that the partition function of the IV-GFF with Dirichlet boundary conditions may be rewritten 
\begin{align*}\label{}
	Z_{\beta,\Lambda,\lambda_\Lambda}^0 = \sum_{\calN \in \calF} c_\calN \int \prod_{\rho\in \calN} [1+z^0(\beta,\rho, \calN) \cos(\<{\phi,\bar\rho})] d\mu_{\beta,\Lambda,v}^{0,\GFF}(\phi)\,,
\end{align*}
where the superscript $0$ stands for Dirichlet boundary conditions. The main difference with {\em free} or {\em periodic} boundary conditions is that the ensemble of charges $\calN$ may now contain some non-neutral charges $\rho \in \calN$ (i.e. s.t. $Q(\rho):=\sum_{x\in \Lambda} \rho(x) \neq 0$).

\medskip
 \uline{The  first main essential difference between \cite{wirth2019maximum} and \cite{FS} is as follows}: the splitting/merging algorithm from \cite{FS} in order to obtain such an expansion into charges requires some key adjustments as we shall see below.  The important feature of this (new) expansion into charges is also to obtain as in the free/periodic case a powerful enough  upper bound on ``activities'' $z$ under some (very general) assumptions on the growth of the Fourier coefficients $|\hat \lambda_i(q)|$. The upper bound below is the content of  \cite[Proposition 22]{wirth2019maximum} after a suitable expansion:
\begin{align}\label{e.goodz0}
|z^0(\beta, \rho, \calN)|
\leq \exp \left( -\frac c \beta (\|\rho\|_2^2 + \log_2(d_\Lambda(\rho) + 1))\right)\,,
\end{align}
where $d_\Lambda(\rho)$ is a key new definition from \cite{wirth2019maximum}. It is one of the important adjustments: instead of measuring the ``size'' of a charge $\rho$ by its diameter as it is done with free/periodic boundary conditions (notice the difference with~\eqref{e.zfree}), the size of a non-neutral charge now depends also on its distance to the boundary $\p \Lambda$. Wirth introduces the following ``modified diameter'':
\begin{align}\label{e.MD}
d_\Lambda(\rho):= 
\begin{cases}
\max \{ \dist(\rho, \p \Lambda), \diam(\rho) \} & \text{ if } Q(\rho)\neq 0 \\
\diam(\rho) & \text{ if } Q(\rho)=0
\end{cases}
\end{align}

Furthermore, this expansion into Coulomb charges comes with some new constraints which need to be satisfied by each ensemble of charges $\calN$. Let us highlight the one which will be of most significance to us: for any ensemble $\calN$,  any two charges $\rho\neq \rho' \in \calN$ need to satisfy
\begin{align}\label{e.constraint1}
\dist(\rho,\rho') \geq M [\min(d_\Lambda(\rho), d_\Lambda(\rho'))]^{\bar \alpha}
\end{align}
where $M\geq 1$ and $\bar \alpha\in(\tfrac 3 2, 2)$ are  parameters\footnote{in these references, the second one is called $\alpha$ but in the present paper, $\alpha$ has another meaning.} which play an important role in the splitting/merging algorithm in \cite{FS,wirth2019maximum}. 
Due to the definition of the modified diameter in~\eqref{e.MD}, this ``repulsion of charges'' in particular implies that it is impossible to have more than 1 non-neutral charge which intersects the sub-domain $\Lambda_{\frac n 2}^2 \subset \Lambda_n^2$. 
\smallskip

Besides the need of a different splitting/merging algorithm,  the use of overlapping spin-waves also needs to be adapted to the unpleasant presence of non-neutral charges in order to obtain a good enough bound on the activities~\eqref{e.goodz0}. \uline{This is the second main key main difference with \cite{FS}} (on which we shall not elaborate further, see \cite{wirth2019maximum}).

Once this is achieved, exactly as in Section \ref{ss.sketchFS}, the Laplace transform $\EFK{\beta,\Lambda,\lambda_\Lambda, v}{0}{e^{\<{\phi,f}}}$ may be rewritten as follows: 
\begin{align*}\label{}
	e^{\frac 1 {2 \beta} \<{f, (-\Delta^0)^{-1} f}}
	\frac 
	{
		\sum_{\calN \in \calF} c_\calN \int \prod_{\rho\in \calN} [1+z^0(\beta,\rho, \calN) \cos(\<{\phi,\bar\rho} + \<{\sigma, \rho})] d\mu_{\beta,\Lambda,v}^{0,\GFF}(\phi)
	}
	{
		\sum_{\calN \in \calF} c_\calN \int \prod_{\rho\in \calN} [1+z^0(\beta,\rho, \calN) \cos(\<{\phi,\bar\rho})] d\mu_{\beta,\Lambda,v}^{0,\GFF}(\phi)
	}\,,
\end{align*}
where $\sigma$ is now defined as 
\begin{align}\label{e.sigma0}
	\sigma:= \frac 1 \beta [- \Delta^0]^{-1} f 
\end{align}
Notice we are now working with the Laplacian $\Delta^0=\Delta_\Lambda^0$ with Dirichlet boundary conditions on $\p \Lambda$. When the context is clear, to the cost of slight abuse of notation,  we will keep denoting it by $\Delta$.

As in the free/periodic case, the goal is thus to obtain a lower-bound on each ratio
\begin{align*}\label{}
\frac {Z^0_{\calN}(\sigma)} {Z^0_{\calN}(0)}\,,
\end{align*}
which are defined as in Section \ref{ss.sketchFS}.  Here, the same analysis as for the free/periodic case (i.e. the use of Jensen's inequality etc.)  leads us to  
\begin{align}\label{e.nice1DIR}
	\frac {Z^0_{\calN}(\sigma)} {Z^0_{\calN}(0)}  \geq &  
	\exp\big[-D_1  \sum_{\rho\in \calN} |z^0(\beta,\rho,\calN)|\<{\sigma,\rho}^2\big] \,.
\end{align}

We thus come to \uline{the third main difference between \cite{wirth2019maximum} and \cite{FS}}.  Since some of the charges $\rho\in \calN$ may not be neutral, we cannot proceed as in Subsection \ref{sss.dir} and decompose $\<{\sigma,\rho}$ as 
\begin{align*}\label{}
\<{\sigma,\rho} = \sum_{i\sim j \in D(\rho)} c_{i,j}( \sigma(i)-\sigma(j) )\,.
\end{align*}
To overcome this issue, Wirth proceeds as follows: 
\bnum
\item First, if $\Lambda=\Lambda_n^2 =\{-n,\ldots,n\}^2 \subset \Z^2$, he has chosen to consider only test functions $f$ which are such that $\sigma:= \frac 1 \beta [- \Delta^0]^{-1} f$ vanishes everywhere on the annulus 
$A_n = \Lambda_n^2 \setminus \Lambda_{\tfrac n 2}^2$.  (Another aspect ratio is used in \cite{wirth2019maximum} but this is a minor technicality). 
\item By the constraint which followed~\eqref{e.constraint1}, this implies that at most one non-neutral charge $\rho$ may give a non-trivial contribution in~\eqref{e.nice1DIR}. 
\item If there exists a (unique) such non-neutral charge $\rho$, a vertex $v^*$ is picked in $\Lambda^2_{\tfrac 2 3 n} \setminus \Lambda_{\tfrac n 2}^2$. The key point is to notice is that one may correct the lack of neutrality of $\rho$ by defining 
\begin{align*}\label{}
\rho^*:= \rho + Q \delta_{v^*}
\end{align*}
and to notice that since $\sigma_{\md (\Lambda^2_{\frac n 2})^c} \equiv 0$, one still has 
\begin{align*}\label{}
\<{\sigma, \rho} = \<{\sigma, \rho^*}\,.
\end{align*}
To apply the same analysis as in Subsection \ref{sss.dir}, it is necessary to check that the  added vertex $v^*$ does not carry too much charge, namely $Q + \rho(v^*)$. But clearly one has $|Q + \rho(v^*)| = |\rho^*(v^*)| \leq 2 \| \rho \|_1 $ (the facteur $2$ is due to the fact $v^*$ may already carry a non zero charge in $\rho$). 

\item Since $\rho$ is intersecting $\Lambda_{\frac n 2}^2$ and since $\rho$ is assumed to be non-neutral, this implies that its modified diameter $d_\Lambda(\rho)$ (defined in~\eqref{e.MD}) is larger than $\Omega(n)$. As such, one has $v^* \in D_\Lambda(\rho)$  (the square box surrounding $\rho$ of diameter of order $d_\Lambda(\rho)$, see \cite{wirth2019maximum}).  This allows us to write
\begin{align*}\label{}
\<{\sigma, \rho^*} = \sum_{i\sim j \in D_\Lambda(\rho) \cap \T_n \times \{0\}} c_{i,j} (\sigma(i) -\sigma(j))\,.
\end{align*}
Thanks to the previous control on the total charge at $v^*$, one obtains as in Subsection \ref{sss.dir}
\begin{align*}\label{}
|c_{i,j}| \leq \tfrac 1 2  \|\rho^*\|_1 \leq    \tfrac 3 2 \|\rho\|_1 \leq \tfrac 3 2 \|\rho\|_2^2
\end{align*}
\item Still arguing as in Subsection \ref{sss.dir}, one can conclude that 
\begin{align*}\label{}
\sum_{\rho\in \calN}  |z^0(\beta,\rho,\calN)|\<{\sigma,\rho}^2  
& \leq \sum_{\rho \in \calN}  |z^0(\beta,\rho,\calN)|\,  |D(\rho)| \cdot  \| \rho\|_2^4  \sum_{i\sim j \in D(\rho)} (\sigma_i - \sigma_j)^2\,.
\end{align*}
A key point here is that the estimate on the activities~\eqref{e.goodz0} involves the modified diameter $d_\Lambda(\rho)$ instead of $d(\rho)$ (this is crucial, as we might have $d(\rho) \ll d_\Lambda(\rho)$ and the first one would not be sufficient to ``renormalize'' correctly the activities). 

\item The main quantitative estimate to conclude the proof of Lemma \ref{l.1516}, i.e. \cite[Theorem 13]{wirth2019maximum} then follows as in Section \ref{sss.dir} by grouping charges depending on their dyadic scale. Since $\sigma$ was assumed to vanish on $\Lambda_n^2 \setminus \Lambda^2_{\frac n 2}$,  at most one non-neutral charge contributes to this grouping (and it contributes only to the largest macroscopic scale). 

\enum

\subsubsection{Extension to more general test functions.}\label{ss.ext}

Our goal in this Subsection is to show that the above analysis extends to  any test function $f$ without assuming that $\sigma= \frac 1 \beta [- \Delta^0]^{-1} f$ vanishes on $\Lambda_n^2 \setminus \Lambda_{\frac n 2}^2$. 

Let us precise the notations here. If $\Lambda$ is a box in $\Z^2$, we will denote by $\p \Lambda$
its (inner) boundary
\begin{align*}\label{}
\p \Lambda = \{ x\in \Lambda, \text{ s.t. there exists } y\notin \Lambda \text{ with } y\sim x \}\,.
\end{align*}
We will also denote by $\Lambda^o:= \Lambda\setminus \p \Lambda$ its interior. 
For any test function $f : \Lambda \to \R$, there exists a unique fonction $G$ which vanishes on $\p \Lambda$ and which is such that for any interior point $x\in \Lambda^o$,  $(-\Delta^0) G(x) = f(x)$. This function is by definition $(-\Delta^0)^{-1} f$. 
 
In the case of $\Lambda=\Lambda_n^2$, if we test the field $\Psi_n^{line}$ against any such function $f$, note that  $\sigma:= \frac 1 \beta [- \Delta^0]^{-1} f$ is vanishing on $\p \Lambda_n$ (where $\Psi_n^{line}$ is also vanishing due to the Dirichlet boundary conditions). 
We do not need to modify the expansion into charges carried in \cite{wirth2019maximum}. As previously, for each configuration of charges $\calN$, we thus reach the lower-bound
\begin{align*}
	\frac {Z^0_{\calN}(\sigma)} {Z^0_{\calN}(0)}  \geq &  
	\exp\big[-D_1  \sum_{\rho\in \calN} |z^0(\beta,\rho,\calN)|\<{\sigma,\rho}^2\big] \,.
\end{align*}
Now,  if we do not make any assumption on $f$, then the support of $\sigma$ may be the whole interior domain $\Lambda_n^2 \setminus \p \Lambda_n^2$. In particular, many non-neutral charge $\rho \in \calN$ will typically be involved in the above sum.  To deal with all of these at once, we proceed slightly differently from \cite{wirth2019maximum} as follows:
\bnum
\item First, in the above sum, the contribution coming from neutral charges is handled exactly as in Subsection \ref{sss.dir} (this is also the case in \cite{wirth2019maximum}). 
\item We need to prove an upper bound on 
\begin{align*}\label{}
\sum_{\rho\in \calN \setminus \calN^{neutral}} |z^0(\beta,\rho,\calN)|\<{\sigma,\rho}^2\,.
\end{align*}
Since the initial Gaussian free field is already rooted at $\p \Lambda_n^2$, notice that all the charges $\rho$ have their support inside the interior $\Lambda_n^2 \setminus \p \Lambda_n^2$. 
For each fixed non-neutral charge $\rho$ in this sum, we shall modify as previously $\rho$ into a neutral charge $\rho^*$ in such a way that 
\begin{align*}\label{}
\<{\sigma, \rho} = \<{\sigma, \rho^*}\,.
\end{align*}
In the case of \cite{wirth2019maximum}, this operation was done on a single charge which intersected the bulk $\Lambda_{\frac n 2}^2$. 

In the present case, for each $\rho$ non-neutral, we associate a vertex 
\begin{align*}\label{}
\rho \in \calN \setminus \calN^{neutral} \mapsto v^*(\rho) \in \p \Lambda_n^2
\end{align*}
in such a way that $\dist(v^{*}(\rho), \rho) \leq \dist(\rho, \p \Lambda_n^2)$.  
 By definition of the modified diameter $d_{\Lambda}$ from~\eqref{e.MD}, we notice that for each non-neutral $\rho$, the new vertex added $v^*(\rho)$ stays inside the box $D_{\Lambda_n^2}(\rho)$ (see \cite{wirth2019maximum} for the precise definition/centering of such boxes which is not of key importance). The fact $v^*(\rho)\in D_{\Lambda_n^2}(\rho)$ will be crucial  below.  Since furthermore $v^*(\rho)\in \p \Lambda_n^2$ and since $\sigma$ vanishes on the boundary, we may assign any charge on $v^*(\rho)$ without affecting $\<{\sigma, \rho^*}$. We thus define the neutral charge
\begin{align*}\label{}
\rho^* := \rho  - (\sum_{x\in \Lambda_n^2 \setminus \p \Lambda_n^2} \rho(x)) \delta_{v^*(\rho)}
\end{align*}
Proceeding as in Subsection \ref{sss.dir}, we thus obtain 
\begin{align*}\label{}
|\<{\sigma,\rho}|^2 = |\<{\sigma,\rho^*}|^2 \leq  |D_{\Lambda_n^2}(\rho)| \cdot 4  \| \rho\|_2^4  \sum_{i\sim j \in D_{\Lambda_n^2}(\rho)} (\sigma_i - \sigma_j)^2\,.
\end{align*}
\item To end the proof, we need to sum the above bound over all non-neutral charges $\rho \in \calN \setminus \calN^{neutral}$ by taking advantage of the small activities $z^{0}(\beta,\rho, \calN)$ from~\eqref{e.goodz0} (this later bound follows from \cite{wirth2019maximum}).  It is crucial here to make two observations:
\bi
\item[a)] First the upper bound~\eqref{e.goodz0} is controled by the modified diameter~\eqref{e.MD}. This allows us to compensate the fact the added root $v^*(\rho)$ may be rather far and may thus induce a large area after the use of Cauchy-Schwarz (i.e. $ |D_{\Lambda_n^2}(\rho)|$ may be much larger than the term $|D(\rho)|$ in the case of a neutral charge). 
\item[b)] Once each charge is well-controlled using $a)$, we still need to argue that they do not overlap too much. As in the neutral case, we divide the non-neutral charge depending on their modified diameter and we obtain that such an overlapping does not hold thanks to the key contraint in Wirth's Coulomb gas expansion coming from his repulsion condition~\eqref{e.constraint1}.
\ei
Proceeding as in Subsection \ref{sss.dir}, this ends the proof of our extension Lemma \ref{l.extension}. 
\enum

\subsubsection{The case where neutral and non-neutral charges are restricted to a line.}

We stated the extension Lemma \ref{l.extension} first because it is interesting on its own (it provides a control on Laplace transforms of arbitrary functionals of a Dirichlet IV-GFF as opposed to the restricted class in \cite{wirth2019maximum}), but also in order to explain the main idea of adding ``simultaneously'' as many additional vertices $v^*$ as non-neutral charges $\rho$. 

What we really need in this paper is Lemma \ref{l.line} which handles the case where charges are initially restricted to a $1D$ line in $\Lambda_n^2$. 

By combining the adapation to the case of a $1D$ line explained in the case of free/periodic boundary conditions in the Subsection \ref{ss.linefree} with the above way of handling non-neutral charges, we claim that one readily obtains Lemma \ref{l.line}. Neutral charges are handled exactly as in Subsection \ref{ss.linefree}. The only needed adaptation is that for the non-neutral charges $\rho$, we need to pick a vertex $v^*(\rho)$ which can be reached via a path inside the $1D$ line. For this we simply pick the closest point in $\Z\times \{0\} \cap \p \Lambda_n^2$. 
Note that because of this, many distinct non-neutral charges $\rho$ may share the same vertex $v^*(\rho)$. This is not an issue as Cauchy-Schwarz is applied separately for each charge $\rho$ (rather $\rho^*$ here). The only issue could come from a superposition of non-neutral charges of same scale $2^k$. This is ruled out once again by the repulsion \ref{e.constraint1}. 
We obtain this way a proof of Lemma \ref{l.line}.


\section{Bessel type random walks on the diamond graph}
\label{s.Bessel}

For each $s\geq 0$, the purpose of this section is to build natural Markov processes $(Z_n^{(s)})_{n\geq 0}$ on the diamond graph $\calD = e^{i \pi/4} 2^{-1/2} \Z^2$ which have the property that there exists a constant $c>0$ such that for $a,b \in \Z \times \{0\}$, 
\begin{align*}\label{}
\FK{}{a}{Z^{(s)}_{\tau_{\Z}} = b} \sim \frac {c}  {\|a-b\|^{2+s}}\,,
\end{align*}
as $\| a -b \| \to \infty$ and where $\tau_\Z$ is the first return time ($\geq 1$) to the horizontal line $\Z \times \{0 \}$.


Recall that in the case $s=0$, it is sufficient to consider the simple random walk on the diamond graph $\calD$. Furthermore in this case, the return probabilities are explicit and are given by~\eqref{e.PEx}.  When $s>0$ one needs to introduce some confining force driving the process towards the line $\Z$. This will be achieved as follows: at each time $n\geq 0$, the Markov process  $(Z_n^{(s)})_{n\geq 0}$ will:
\bi
\item move vertically up or down (i.e. by $\pm \frac 1 {2}$ since recall $\calD=e^{i \pi/4} 2^{-1/2} \Z^2$) according to a {\em discrete Bessel random walk} as studied for example in \cite{alexander2011excursions}. 
\item And simultaneously make an independent horizontal right or left move (See Figure \ref{f.diamond}). 
\ei

There is a large variety of such Bessel walks considered in \cite{alexander2011excursions}, we will only consider two: this one on the diamond graph $\calD$ as well as the  Bessel walks on the unrotated $\Z^2$ introduced in appendix \ref{s.BesselZ}. We note that any other choice from \cite{alexander2011excursions} would yield a slightly different model of $\alpha$-Discrete Gaussian Chain for which the results of Theorem \ref{th.main} would still hold. (N.B. As explained in details in \cite{alexander2011excursions} some other choices of Bessel walks would induce inevitable slowly varying functions in the Lemmas below which would be rather inconvenient).

The discrete Bessel process  we shall consider will have the following Markov transition kernel $Q_s: \N \times \N \to [0,1]$. 

\bi
\item $Q_s(0,1)=1$\footnote{Note that this first condition ensures that the walk remains inside the upper-half plane $\Z \times \N$ instead of $\Z^2$. This is harmless in our analysis and if one wishes to, one could as well reflect up or down with equal probability at each return to $\Z\times \{0\}$.}.
\item $Q_s(r,r-1)+Q(r,r+1)=1$ whenever $r\geq 1$.
\item $Q_s(r,r+1)= (\tfrac 1 2  - \frac s 4\, \frac 1 r ) \vee \tfrac 1 4$ whenever $r\geq 1$.
\ei

\begin{remark}\label{}
Recall that the classical continuous time Bessel processes on $\R^+$ are given by the SDE
\begin{align*}\label{}
dY_t = \frac a {Y_t} dt + dB_t\,,
\end{align*}
where $a=\frac{d-1}2$ ($d$ being the ``dimension'' of the Bessel process). See for example \cite{revuz2013continuous, lawler2008conformally}. We therefore have the correspondance $s \equiv  - 2 a$ (indeed our discrete Bessel processes are defined in such a way that they have a negative drift when $s\geq 0$, note also that this drift towards 0 is precisely \underline{twice} $\frac s {4r}$). 
 In our definition, we may also use any fixed value of $s\in (-1,0]$ which would then correspond to $a\in[0,1/2)$. This is consistent with the fact that the continuous Bessel process with $a=1/2$ never reaches the origin. 
\end{remark}

We now state a fine analysis of the first return times of such Bessel random walks due to \cite{alexander2011excursions}.
\begin{theorem}[Theorem 2.1 in \cite{alexander2011excursions}]\label{th.ALEX}
For any $s>-1$, there exists a constant $c=c(s)>0$ such that as $n$ even goes to infinity,  one has  
\begin{align*}\label{}
\FK{Q_s}{0}{\tau_0 = n} \sim c(s)\, n^{-(\frac {3+s} 2)}\,, 
\end{align*}
where $\FK{Q_s}{0}{\tau_0 = n}$ is the probability of first return to the origin for the above Bessel random walk on $\N$.
\end{theorem}

We shall now prove that the above constructed process on the diamond graph $\calD$ satisfies the desired property:
\begin{proposition}\label{pr.Bessel}
For any $s>-1$,  (only the case $s\geq 0$ will be useful in this work), there exists a constant $b=b(s)>0$ such that as $k\to \infty$ 
\begin{align*}\label{}
\FK{}{(0,0)}{Z^{(s)}_{\tau_{\Z}} = (0,k)} \sim \frac {b(s)} {k^{2+s}} \,.
\end{align*}
\end{proposition}

\ni
{\em Proof.} 

The main input of the proof is the above Theorem \ref{th.ALEX} from \cite{alexander2011excursions}. 
Let us set two notations:
\bnum
\item Call the function $g_s(n):= \FK{Q_s}{0}{\tau_0 = n}$.
\item Call the heat-kernel of the simple random walk $p_{\Z}(t,x,y)=p_\Z(t,x-y)$. 
\enum

For any  $k\geq 1$, we have 
\begin{align*}\label{}
\FK{}{(0,0)}{Z^{(s)}_{\tau_{\Z}} = (0,k)} & = \sum_{n = k}^\infty  p_\Z(2n,2k) g_s(2n)\,.
\end{align*}
The factors 2 are here because each odd move goes to non-integer points in $\calD$ while each even move lands at an integer point in $\calD$. 
 
%
%

Let $\eps>0$ be any fixed small parameter. 
Our first observation is that the contributions of integers $n\leq k^{2-\eps}$ is negligible. Indeed, in the sum 
\begin{align*}\label{}
\sum_{n=k}^{k^{2-\eps}} p_\Z(2n,2k) g_s(2n),
\end{align*}
each heat kernel is upper bounded (for small enough $c>0$) by 
\begin{align*}\label{}
\exp(- c (2k)^2/ k^{2-\eps} ) \leq \exp(- c k^\eps)\,.
\end{align*}
As such, the sum $\sum_{n=k}^{k^{2-\eps}}$ is negligible w.r.t. to the power law in $k$ asymptotics that  we are looking after.

%

Now, by the quantitative Local CLT theorem (as can be found for example in \cite{lawler2010random}),  there exists a constant $C>0$ s.t. uniformly in $n$ and in  $k$,  
\begin{align*}\label{}
\left| p_\Z(2n ,2k)  - \frac{1}{\sqrt{2\pi (2n)}} e^{-\frac{(2k)^2}{2*(2n)}}\right| \leq \frac C {n^{3/2}} \,.
\end{align*}

Let us first deal with the upper bound (a matching order lower bound follows from the same analysis in order the conclude the proof of the asymptotics in Proposition \ref{pr.Bessel}). 

For the Bessel term $g_s(2n)$, using Theorem \ref{th.ALEX}, we have that for any small $\delta>0$, then for $n$ large enough, 
\begin{align*}\label{}
g_s(2n) \sim c(s)\, (2n)^{-\frac {3+s} 2}
& \leq (1+\delta)  c(s)  (2n)^{-\frac {3+s} 2} 
\end{align*}
This gives us
\begin{align*}\label{}
& \FK{}{(0,0)}{Z^{(s)}_{\tau_{\Z}} = (0,k)}  \\
& \leq  k^2 e^{-c k^{\eps}} +  \sum_{n\geq k^{2-\eps}}
\frac{1}{\sqrt{4\pi n}} \left(e^{-\frac{k^2}{n}} +  \frac C {n^{3/2}}\right)
(1+\delta)  c(s)  (2n)^{-\frac {3+s} 2} 
\end{align*}

The leading term is given as $k\to \infty$ 
\begin{align*}\label{}
(1+\delta) 2^{-\frac{3+s} 2} \sum_{n\geq k^{2-\eps}} \frac {n^{-(2+\tfrac s 2)}} {\sqrt{4\pi} }  e^{-\frac {k^2} n}
\end{align*}

%
%

We are thus left with controlling the asymptotic of the above series. By applying for example Euler-MacLaurin comparison's formula, we get
\begin{align*}\label{}
&  \sum_{n\geq k^{2-\eps}} \frac {n^{-(2+\tfrac s 2)}} {\sqrt{4\pi} }  e^{-\frac {k^2} n}\\
& \sim_{k\to \infty}
\int_{n=k^{2-\eps}}^\infty
\frac 1 {\sqrt{4\pi}} dn\,   n^{-(2+\frac s 2)} e^{-\frac{ k^2}{n}} \\
& = (n=k^2 u) \,\, k^2 \int_{k^{-\eps}}^\infty k^{- 4 -s}  \frac 1 {\sqrt{4\pi}}  du \frac 1 {u^{2+ \frac s 2}} e^{-\frac 1 u } \\
& \sim \frac 1 {k^{2+s}} \int_0^\infty \frac 1 {\sqrt{4\pi}} du \frac{e^{-\frac 1 u}}{u^{2+\tfrac s 2}}
\end{align*}
With the same analysis, we also obtain for any $\delta>0$ and for $n$ large enough, 
\begin{align*}\label{}
\FK{}{(0,0)}{Z^{(s)}_{\tau_{\Z}} = (0,k)}  \geq (1-\delta) 2^{-\frac{3+s} 2} c(s) \left(  
\int_0^\infty \frac {du} {\sqrt{4\pi}}  \frac{e^{-\frac 1 u}}{u^{2+\tfrac s 2}} \right) 
\frac 1 {k^{2+s}}\,, 
\end{align*} 
which concludes the proof of Proposition \ref{pr.Bessel}. 
\qed

See Appendix \ref{s.BesselZ} for a generalisation of this analysis to the case of Bessel walks on $\Z^2$ as well as $\Z^{d+1}$.

\section{Proof of the main Theorem (Theorem \ref{th.main}) in the case $\alpha\in(2,3)$}\label{s.proof}

\subsection{A $2D$ Gaussian Free Field with inhomogeneous conductances.}

The goal of this subsection is to define a Gaussian field $\varphi_{\alpha,n}$ on the diamond graph $\calD$ with Dirichlet boundary conditions on 
\begin{align*}\label{}
L_n = \left(\{ \ldots, -n-1\} \cup \{n+1,\ldots \}\right)\times \{0\} \subset \calD\,,
\end{align*}
and such that the restriction of this Gaussian field to the $1D$ interval $\Lambda_n \times \{0 \}$ has the desired density. 

 In what follows, we will denote by $D_n$ the two-dimensional slit domain 
\begin{align}\label{e.slit}
D_n:= \calD \setminus L_n  = \calD \setminus \left(\{ \ldots, -n-1\} \cup \{n+1,\ldots \}  \times \{0\}  \right)\,.
\end{align}

For any $s\geq 0$, let us introduce the following field $\mathbf a = (a_{i,j})_{i\sim j}$ of conductances on the graph $\calD$:
\bi
\item All edges which intersect the base line $\Z \times \{0\}$ carry a conductance equal to $\tfrac 1 4$ (this is a way to normalize our conductances which matches with our normalisation in the case $s=0$ from~\eqref{e.gff8}). 
\item Conductances which are at the same height are equal. More precisely, for every edge on $\calD=\frac{e^{i\pi/4}}{\sqrt{2}} \Z^2$, there  exists an integer $r \in \Z$ so that the top vertex of the edge is at height $\frac {r+1} {\sqrt{2}} $ while the bottom vertex is at height $\frac r {\sqrt{2}}$. All the edges which share the same integer $r$ will carry the same conductance which we shall denote with a slight abuse of notation by $a(r,r+1)$. 
\item The field of conductances is symmetric under vertical reflection, i.e. $a(r,r+1) = a(-r-1,-r)$. 
\item We now fix the constraint coming from the vertical discrete Bessel kernel: for any $r\geq 1$,  let 
\begin{align*}\label{}
\frac{a(r,r+1)}{a(r-1,r)+a(r,r+1)} := Q_s(r,r+1)= (\tfrac 1 2  - \frac s 4\, \frac 1 r ) \vee \tfrac 1 4\,. 
\end{align*}
\ei
This readily implies that for any $r\geq r_0(s)$, 
\begin{align*}\label{}
a(r,r+1) = \frac{1  - \frac{s}{2r}}{ 1  + \frac{s}{2r}} a(r-1,r) = (1-\frac{s}{r}+O(r^{-2})) a(r-1,r). 
\end{align*}
This in turn implies (using  $\log(1+x)\leq x$ and the fact $e^x$ is monotone) 
\begin{align*}\label{}
a(r,r+1) &  \leq O(1) \exp( \sum_{i=1}^r \log(1- \frac s i  + O(1/i^{2}))) \\
& \leq O(1) \exp( \sum_{i=1}^r  -  \frac s i + O(1/i^2))) \\
& \leq O(1) r^{-s}\,.
\end{align*}
(N.B. It is not difficult to show using $\log(1+x) \geq x - c x^2$ in the vicinity of the origin, that  $a(r,r+1) \asymp r^{-s}$ but we shall not need this fact).

\medskip

Let us now consider the Gaussian free field $\varphi_{\alpha,n}$ on the diamond graph $\calD$ in the above field of conductances $\mathbf{a}$ and rooted on the set  $L_n$. (See Figure \ref{f.diamond}).  

This means that $\varphi_{\alpha,n}$ is the Gaussian field with  (formal) density 
\begin{align*}\label{}
\exp\left(-\frac \beta 2 \sum_{i\sim j \in \calD} a_{i,j} (\varphi_{\alpha,n}(i) - \varphi_{\alpha,n}(j))^2 \right) 1_{\varphi_{\alpha,n} \equiv 0 \text{  on  } L_n}
\end{align*}
By a straightforward discrete integration by parts, we may then rewrite this density as follows
\begin{align*}\label{}
\exp\left(-\frac \beta 2 \<{\phi, A \phi}\right) 1_{\varphi_{\alpha,n} \equiv 0 \text{  on  } L_n}\,,
\end{align*}
where $A$ is the (non-isotropic) elliptic operator defined for any site $i\in \calD$ by 
\begin{align*}\label{}
A\phi(i):= \sum_{j\sim i\in \calD} a_{i,j}[\phi(j) - \phi(i)]\,.
\end{align*}
It is convenient to decompose this symmetric operator as 
\begin{align*}\label{}
A=D_{\mathbf{a}} (1-P_\mathbf{a})\,,
\end{align*}
where $D_\mathbf{a}$ is the diagonal matrix whose $i^{th}$ component is $\sum_{j \sim i} a_{ij}$ and $1-\P_{\mathbf{a}}$ is the non-symmetric Markov transition matrix of the random walks in conductances $\mathbf{a}$ killed when first hitting the slit $L_n$. We thus obtain that the covariance structure of the field $\varphi_{\alpha,n}$ is given by 
\begin{align*}\label{}
\Eb{\varphi_{\alpha,n}(i) \varphi_{\alpha,n}(j)} = \frac 1 \beta  A^{-1}(i,j)&  = \frac 1 \beta \left[(1-P_{\mathbf{a}})^{-1} D_{\mathbf{a}}^{-1}\right](i,j) \\
& = \frac 1 \beta \frac 1 {\sum_{k\sim j} a_{j,k}} G_{D_n,\mathbf{a}}(i,j)\,,
\end{align*}
where $G_{D_n,\mathbf{a}}$ is the Green's function of the $\mathbf{a}$-random walk, i.e.
\begin{align}\label{e.Greena}
G_{D_n,\mathbf{a}}(i,j):=\EFK{}{i}{\sum_{n=0}^\infty 1_{Z_n^{\mathbf{a}} = j}}\,.
\end{align}
(N.B notice that $G_{D_n,\mathbf{a}}$ is not symmetric but symmetry is of course restored after multiplying by $D_{\mathbf{a}}^{-1}$).
By our normalisation of the conductances intersecting the line $\Z \times \{0\}$, notice that for any point $i$ in this line, we have $\sum_{j\sim i} a_{i,j} = 4*\tfrac 1 4 = 1$.  This implies that for any $x,y\in \Z$, we have 
\begin{align*}\label{}
\Eb{\varphi_{\alpha,n}((x,0)) \varphi_{\alpha,n}((y,0))} = \frac 1 \beta G_{D_n, \mathbf{a}}((x,0), (y,0))
\end{align*}

Recall that we denote with a slight abuse of notation the interval $\Lambda_n := \Lambda_n \times \{0\} = \{-n,\ldots,n\} \times \{0\} \subset \calD$. As in the case $s=0$ (see Lemma \ref{l.gffs0}), we have the following identity in law. 

\begin{lemma}\label{l.gffs}
For any $\alpha\in(2,3)$. the restriction of the $\beta$-GFF field $\varphi_{\alpha,n}$ to the $1D$ interval $\Lambda_n \subset \calD$ is equal in law to the Gaussian vector $\{\bar \varphi_i\}_{i\in \Z}$  with density 
\begin{align*}\label{}
\exp(-\frac {\beta} 2 \sum_{i \neq j} J_{\alpha}(i-j) (\bar \varphi_i - \bar \varphi_j)^2) 1_{\bar \varphi_i = 0 \forall i \notin \Lambda_n}\,,
\end{align*}
and where the coupling constants $J_\alpha(r)$ are defined for any $r\in \Z$ by 
\begin{align}\label{e.Jalpha}
J_\alpha(r) :=  \FK{}{(0,0)}{Z^{(s)}_{\tau_{\Z}} = (0,r)} \sim_{r\to \infty} \frac {b(s)} {r^{\alpha}}  \,\,\, \text{ with $s=\alpha-2$}\,.
\end{align}
\end{lemma}

\ni
{\em Proof.} 
We proceed exactly as in the proof of Lemma \ref{l.gffs0}. Let us  compute the covariance matrix of the field $\{\bar \varphi_i\}_{i\in \Lambda_n}$.  For any $i,j \in \Z$, 
\begin{align*}\label{}
\Eb{\bar \varphi_i \bar \varphi_j} & = \EFK{}{i}{\sharp \text{visits to $j$ for the $J_\alpha$-long-range RW  killed when exiting $\Lambda_n$} } \\
& = \EFK{}{(i,0)}{
\sharp \text{visits to $(j,0)$ for the Bessel RW $Z_n^{(s)}$ in $\calD$ killed when hitting $L_n$}} \\
& = \Eb{\varphi_{\alpha,n}(i) \varphi_{\alpha,n}(j)}\,.
\end{align*}
This ends the proof by definition of the GFF $\varphi_{\alpha,n}$ with $0$-boundary conditions on $L_n:=\{ \ldots, -n-1\} \cup \{n+1,\ldots \}$.
\qed

\subsection{First Green functions estimates.}\label{ss.Green}

We shall only focus here on the  less standard case $\alpha\in(2,3)$ (i.e. $s\in (0,1)$) since the case $\alpha=2$ which corresponds to the GFF in the slit domain $D_n$ is more standard. 
Our first proposition gives an upper bound on the Green function. 

\begin{proposition}\label{pr.Green1}
For any $s\in (0,1)$, there exists a constant $C>0$ such that uniformly in $x,y \in \Lambda_n \times \{0\}$, 
\begin{align*}\label{}
G_{D_n, \mathbf{a}}(x,y) \leq C n^{\alpha-2}\,.
\end{align*}
\end{proposition}

\ni
{\em Proof.}

Let $Q_n$ be the square of radius $10n$ minus the slit, i.e. $B_{\|\cdot\|_\infty}(0,10n)\setminus L_n$. We first claim that it is enough to show that for any $x,y\in \Lambda_n \times \{0\}$, one has 
the bound for the Green function in the bounded domain $Q_n$, i.e. the existence of $\tilde C>0$ such that
\begin{align}\label{e.Qn}
G_{Q_n, \mathbf{a}}(x,y) \leq \tilde C n^{\alpha-2}\,.
\end{align}
Indeed, by Markov property, one can compute the larger Green function $G_{D_n,\mathbf{a}}$ by decomposing the $\mathbf{a}$-random walk into consecutive excursions defined as follows:
\bi
\item Odd excursions start at $\Lambda_n \times \{0\}$ until they first reach the boundary of the square $\p B_{\|\cdot\|_\infty}(0,10n)$.
\item Even excursions start from the boundary  $\p B_{\|\cdot\|_\infty}(0,10n)$ until they first reach $\Lambda_n \times \{0\}$.
\ei 
The classical important point to notice here is that uniformly on the starting point $z\in \p B_{\|\cdot\|_\infty}(0,10n)$,  the $\mathbf{a}$-walk has a positive probability to reach $L_n$ before reaching $\Lambda_n \times \{0\}$. If one has a uniform bound on $G_{Q_n, \mathbf{a}}(x,y)$ in the form $O(1) n^{\alpha-2}$ we thus obtain a uniform bound on $G_{D_n, \mathbf{a}}(x,y)$ by a geometric series involving the above uniformly positive hitting probability as well as the bound $O(1) n^{\alpha-2}$.

\smallskip
We are thus left with proving~\eqref{e.Qn}.
Let us rewrite the Green function as 
\begin{align*}\label{}
G_{Q_n,\mathbf{a}}(x,y)
&= \sum_{t\geq 0}^\infty   \EFK{}{x}{1_{T_n>t}1_{(X_t,Y_t)=y} } \,,
\end{align*}
where $T_n$ is the hitting time of the rooting $\p B_{\|\cdot\|_\infty}(0,10n) \cup L_n$ and where $(X_t,Y_t)$ denotes the horizontal and vertical coordinates of the $\mathbf{a}$-walk $Z_t$.

We will need the following Lemma which handles the probability of return for the vertical direction, namely

\begin{lemma}\label{l.Renew}
For any $s\in [0,1)$ (the parameter of our $2d$ Bessel walks in Section \ref{s.Bessel}), there exists a constant $C>0$ s.t. for any $t\geq 1$, 
\begin{align*}\label{}
\FK{}{0}{Y_t=0} \leq \frac {C}{t^{\frac {1-s} 2}}\,.
\end{align*}
\end{lemma}

\ni
{\em Proof of the Lemma.}

The sequence of times when the vertical component of the Bessel walk returns to the origin corresponds to a \underline{Renewal process} $0 < \tau_1 < \tau_1+\tau_2 < \ldots$  whose law of i.i.d jumps is given by 
\begin{align*}\label{}
\Pb{\tau_i = k } := \FK{}{0}{Y_j \neq 0 \,\, \forall j<k \text{ and } Y_k =0}\,.
\end{align*}
This probability is exactly controlled by Theorem \ref{th.ALEX} from \cite{alexander2011excursions}. This gives us 

\begin{align*}\label{}
\Pb{\tau_i = k } \sim c(s)\, k^{-(\frac {3+s} 2)}\,.
\end{align*}
By applying standard results from renewal processes, see in particular \cite{garsia1962discrete,doney1997one} as well as the discussion in \cite{caravenna2016continuum}, 
this implies that for any time $t\geq 1$, the probability that $t$ belongs to the renewal set is upper bounded by 
\begin{align*}\label{}
C \left( \frac 1 {t} \right)^{\tfrac {1-s} 2} \,\,\, \text{ for some $C>0$. }
\end{align*}

 \qed

\smallskip
We shall now group the times $t$ into two groups. We only sketch the details below and we refer to Sections \ref{ss.Grad1} and \ref{ss.Grad2} where a more elaborate time-decomposition will be used.
\medskip

\medskip

\textbf{First group. $t\leq n^2$.}
In this group, we do not pay attention at the fact the walk may leave the domain $Q_n$ or touch the slit $L_n$. After all we are only looking for an upper bound here. 

By local CLT theorem, the probability that the horizontal walk $X_t$ lands at $y$ at time $t$ is less than $O(t^{-1/2})$. 
Furthermore, by the above Lemma, the probability that the vertical direction $Y_t$ gets back to zero at time $t$ is less than $O(t^{-(1-s)/2})$. 

This implies that the times $t$ in this group contribute at most
\begin{align*}\label{}
\sum_{t=1}^{n^2} C \frac 1 {t^{1 - \tfrac s 2}} \leq O(1) [n^2]^{s/2} = O(1) n^s = O(1) n^{\alpha-2}
\end{align*}

\textbf{Second group. $t>n^2$.} 
We first run the horizontal walk $X_t$ up to time $t/2$. We notice that there exists a constant $a>0$ so that the probability that $X_t$ does not bring $Z_t=(X_t,Y_t)$ outside of $Q_n$ is bounded from above by $a\exp(-a t/n^2)$. This is due to the fact that along each interval of time $n^2$, the walk $X_t$ has positive probability to leave the interval $[-10n, 10n]$. 

In the remaining $t/2$ steps, the local CLT still gives us $O(t^{-1/2})$ probability to reach $y$. The (independent) vertical direction also has probability less than $O(1) t^{-(1 - \tfrac s 2)}$ to be back at $0$, exactly as in the previous group. 

This gives us a contribution less than 
\begin{align*}\label{}
\sum_{t=n^2}^{\infty} C \frac 1 {t^{1 - \tfrac s 2}}  e^{-a t/n^2} \leq O(1) n^{\alpha-2}\,.
\end{align*}
\qed

By using the same technology as in the above proof, we also obtain the following more precise  proposition which gives upper and lower bounds and gives the behaviour of the  Green function close to the boundary:
\begin{proposition}\label{pr.Green2}
For any $s\in (0,1)$, there exist  constants $c_1,c_2>0$ such that uniformly in $x\in \Lambda_n \times \{0\}$, 
\begin{align}\label{e.closeB}
c_1\dist(x,L_n)^{\alpha-2} \leq  G_{D_n, \mathbf{a}}(x,x) \leq c_2\, \dist(x,L_n)^{\alpha-2}\,.
\end{align}
Furthermore for any $\delta>0$, there exists a constant $c_\delta>0$ such that for any $x,y \in \Lambda_n \times \{0\}$ at distance at least $\delta n$ from the boundary, 
\begin{align}\label{e.LBG}
c_\delta n^{\alpha-2} \leq G_{D_n, \mathbf{a}}(x,y) \,\,(\,\leq C\,  n^{\alpha-2})\,.
\end{align}
(The upper bound in parenthesis is the conclusion of Proposition \ref{pr.Green1}). 
\end{proposition}

We only say a few words on its proof: The upper bound in~\eqref{e.closeB} is obtained by the same localisation technique using a box of radius $\frac 1 {10} \dist(x, \p D_n)$ and upper bounding by the induced geometric series. 
The lower bound is obtained by summing all the contributions in the definition of the Green function coming from times $t \leq C \dist(x, \p D_n)^2$. At any such time, the probability to reach $x$ from $x$ in time $t$ is of same order (up to multiplicative constant) as the term $t^{-(1-\frac s 2)}$ used in the proof of Proposition \ref{pr.Green1}. This is because for this range of time, given that $X_t=x$, the conditional probability that the horizontal coordinate went as far as the location of the slit $L_n$ is bounded away from 1, uniformly in $t\leq \dist(x, \p D_n)^2$. Since the first group of times in the proof of Proposition \ref{pr.Green1} already gives the leading order of $n^{\alpha-2}$, this ends the proof. 
The lower bound for $G_{D_n, \mathbf{a}}(x,y)$ follows from the same decomposition into times $t\leq O(\delta^2)$ which induces a positive but potentially low constant $c_\delta$ in~\eqref{e.LBG}.



\subsection{Smoothing and compactifying the domain.}\label{ss.compact}

It turns out that to be rather difficult to handle the Dirichlet energy  (the purpose of Section \ref{s.GRAD}) when we are close to the tip of the slit domain $L_n$.
Indeed, to obtain the desired estimates, one would need to prove a ``sharp'' {\em Beurling estimate} which would quantify the effect of the drift when $s\neq 0$. (See \cite{lawler2010random} for the Beurling estimate in the case $s=0$). This seems rather non-trivial. 
See Remarks \ref{r.tedious1} and  especially \ref{r.tedious2} for an explanation of these technical difficulties.

Instead, we will get rid of the slit type of singularity by first ``regularizing the domain''. 
We shall also need (in Section \ref{s.GRAD}) to work in a bounded region. 

\begin{figure}[!htp]
\begin{center}
\includegraphics[width=\textwidth]{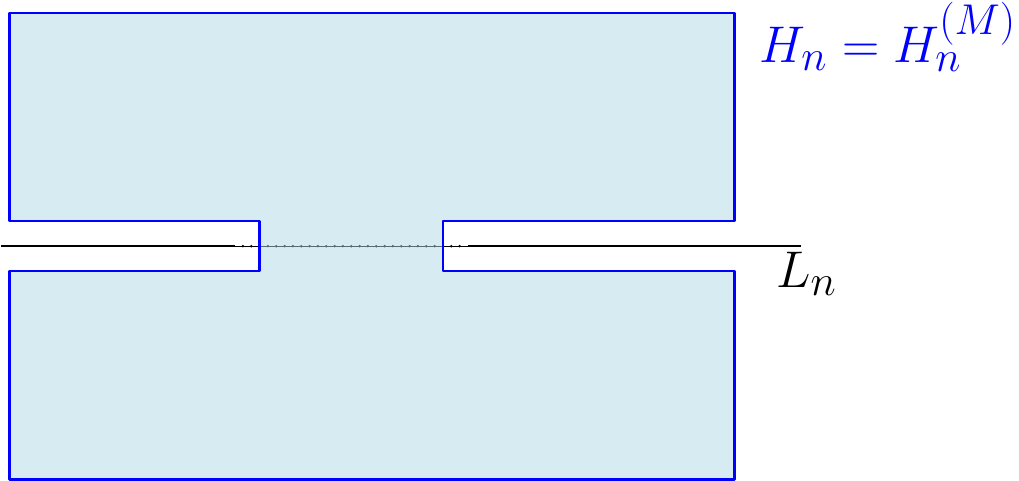}
\end{center}
\caption{The shape of the smoother domain $H_n=H_n^{(M)}$. It is smoother especially near the two tips of the domain $D_n= \calD \setminus L_n$.}\label{f.Hn}
\end{figure}

We will achieve these two regularisations of the domain at once:
for any large parameter $M\geq 1$, we consider the following sub-domain $H_n=H_n^{(M)} \subset D_n$ defined by: 
\begin{align}\label{e.Hn}
H_n = H_n^{(M)}:=  \left\{ z\in D_n \text{ s.t. } \|z\|_\infty \leq M*n \text{ and }  \dist_{\|\cdot\|_\infty}(z,L_n) \geq \frac 1 M n   \right\} \,.
\end{align}
See Figure \ref{f.Hn}. 

%

For any $s \geq 0$ fixed, recall the definition of the Green function $G_{D_n,\mathbf{a}}$ from~\eqref{e.Greena}. We will also denote by $G_{H_n, \mathbf{a},M}$  the Green functions for the random walk in $H_n\subset D_n$ in conductances $\mathbf{a}=\mathbf{a}(s)$ and killed when first hitting the boundary $\p H_n$. 
We shall distinguish two differents parts in this boundary which will play a different role in the proof below:
\bi
\item The points on the square of radius $M\cdot n$. We shall call this part $S_{nM}$. 
\item The points at distance $\frac 1 M n$ from the slit line $L_n$. We shall call this part $\calL_n$. 
\ei
\smallskip

The proposition below shows that both Green functions are very close to each other in the bulk of $\Lambda_n$ when the parameter $M\to \infty$. 
\begin{proposition}\label{pr.DnHn}
For any $s\geq 0$ and any $\delta>0$, there exists a positive function $\xi_\delta(M)$ which goes to zero when $M\to \infty$ and which is such that for any $n\geq 1$ and any  $x,y \in \Lambda_n \times \{0\}$ at distance $\geq \delta n$ from $L_n$, we have 
\begin{align*}\label{}
|G_{D_n, \mathbf{a}}(x,y) - G_{H_n, \mathbf{a},M}(x,y)|  \leq \xi_\delta(M) G_{D_n,\mathbf{a}}(x,y)\,.
\end{align*}
\end{proposition}

\ni
{\em Proof.}
By the monotony properties of Green functions, we clearly have $G_{H_n,\mathbf{a},M}(x,y) \leq G_{D_n,\mathbf{a}}(x,y)$. 
Let $\tau$ be the stopping time  when the $s$-random walk $Z_t$ first hits the boundary  $\p H_n = \calL_n \cup S_{nM}$. Using Markov property, we may write
\begin{align*}\label{}
& G_{D_n,\mathbf{a}}(x,y) -  G_{H_n,\mathbf{a},M}(x,y)\\
& =   \EFK{}{x}{ G_{D_n,\mathbf{a}}(Z_\tau,y) } \\
& = \FK{}{x}{Z_{\tau} \in S_{nM}}  
\EFK{}{x}{ G_{D_n, \mathbf{a}}(Z_\tau,y)  \md Z_{\tau} \in S_{nM}} \\
& + \FK{}{x}{Z_{\tau} \in \calL_{nM}}  \EFK{}{x}{G_{D_n, \mathbf{a}}(Z_\tau,y)  \md Z_{\tau} \in \calL_{nM}} 
\end{align*}

Let us first deal with the first term in this sum: by classical random walk estimates (N.B. the drift induced by $s\geq 0$ is only helping us here),
we obtain a function $\zeta(M)$ going to zero with $M$, such that uniformly in $n$ and the starting point $x$, we have 
\begin{align*}\label{}
\FK{}{x}{Z_{\tau} \in S_{nM}}   \leq \zeta(M)
\end{align*}

We conclude this case by using the following two facts:
\bnum
\item If $s=0$, (i.e. if we consider the SRW on the grid $\calD$), then we have for any point $z$ in the boundary of the square of radius $nM$, i.e.  $S_{nM}$ and for any point $y\in \Lambda_n \times \{0\}$, 
\begin{align*}\label{}
G_{D_n,\mathbf{a}}(z,y) \leq O(1)  \leq \Omega(G_{D_n,\mathbf{a}}(x,y))\,. 
\end{align*}
Here we use the assumption that $x$ and $y$ are at macroscopic distance from the boundary. (Otherwise one would need a slightly more detailed analysis). 

\item If $s\in (2,3)$ as we have seen in Proposition \ref{pr.Green1}, there exists a constant $C=C(s)>0$ such that uniformly in $y\in \Lambda_n\times \{0\}$ and  $z$ in the boundary of the square $S_{nM}$, $G_{D_n,\mathbf{a}}(z,y) \leq C n^{\alpha-2}$. Furthemore, as $x$ and $y$ are at distance $\delta$ from $L_n$, we also have $G_{D_n,\mathbf{a}}(x,y) = \Omega(n^{\alpha-2})$.  
\enum

The second case is handled similarly with one essential difference:
as opposed to the previous case, it is not true in this case that the probability $\FK{}{x}{Z_{\tau} \in \calL_{nM}}$ degenerates as $M\to \infty$. What is helping us in this second case is that if $z$ lies anywhere on the part of the boundary $\calL_n$, then by definition, $z$ is $\frac 1 M n$ close to the slit $L_n$. 
In particular, it is not difficult to see that we can find a function $\hat \zeta_\delta(M)$ which goes to zero when $M\to \infty$ and which is such that uniformly in $n$ as well as in $z\in \calL_n$, one has 
\begin{align*}\label{}
\FK{D_n, \mathbf{a}}{z}{\tau_y <\infty} < \hat \zeta_\delta(M)\,,
\end{align*}
where the above probability is the probability in $D_n$ starting from $z$ to reach $y$ before intersecting $\p D_n$. Note that the above functions will tend to zero only when the parameter $M$ becomes $\gg \delta^{-1}$, otherwise it is not true that this probability is small. 

Now when $s=0$, by using another stopping time on a ball of radius $\delta^2 n$ around $y$, we obtain the desired result.

When $s\in (0,1)$, we use the combined facts below in order to conclude:  
\begin{align*}\label{}
G_{D_n,\mathbf{a}}(z,y) = \FK{D_n, \mathbf{a}}{z}{\tau_y <\infty} \times G_{D_n,\mathbf{a}}(y,y)\,, 
\end{align*}
and 
\begin{align*}\label{}
G_{D_n,\mathbf{a}}(y,y) \leq O(1) G_{D_n,\mathbf{a}}(x,y)\,.
\end{align*}
(the second point above corresponds to the estimate~\eqref{e.LBG}. Notice that it would be wrong when $s=0$, whence the additional stopping time mentioned above when $s=0$). 
\qed


\subsection{Fröhlich-Spencer Coulomb gas expansion in an inhomogeneous field of conductances and proof of Theorem \ref{th.main}.}\label{ss.TheProof}
$ $

We start by stating a more precise version (the first half) of our main Theorem \ref{th.main}. Our invariance principle is in the same sense as the invariance principle for the discrete fractional Gaussian field from \cite{lodhia2016fractional} and which we stated in Proposition \ref{pr.sheffield} (i.e. in the sense of the field tested against continuous test functions).  

As in this latter proposition, we need to set some notations and some appropriate rescaling for our integer-valued field $\Psi_n$.
As in~\eqref{e.rescaled}, 
for any $\beta>0$, $n\geq 0$, let us consider the \uline{rescaled field} on $\frac 1 n \Z$:
\begin{align}\label{e.rescaledPsi}
\bar \Psi_n(t):=  \frac 1 {n^{H}}  \Psi_n(n t)\,\,, \,\,\, \forall t\in \frac 1 n \Z\,.
\end{align}

We may now state our more precise version of our main theorem:

\begin{theorem}[Proposition 12.2 in \cite{lodhia2016fractional}]\label{th.main2}
For any $\alpha\in[2,3)$, let $J_\alpha(r)\sim c r^{-\alpha}$ be the coupling constants built either in Section \ref{s.Bessel} or in Appendix \ref{s.BesselZ}. 
Then there exists an inverse temperature $\beta^*$ and a constant $K=K(J_\alpha)>0$ such that the following holds.
For any $\beta\leq \beta^*$, the distribution $\mathbf{H}_n:=\sum_{t\in \frac 1 n \Z} \bar \Psi_n(t)  \frac 1 n \delta_t$ converges in law to  $\frac{K}{\sqrt{\beta}} B_t^H$ (where  $B_t^H$ is the Dirichlet $H$-fBm on $D=(-1,1)$ defined in Definition \ref{d.fBm3}). The convergence is in the sense that for any test functions $f_1,\ldots,f_k \in C_c(D)$, 
\begin{align*}\label{}
(\<{\mathbf{H}_n,f_1},\ldots,\<{\mathbf{H}_n,f_k}) \,\,\, \overset{(d)}\longrightarrow_{n\to \infty} \,\,\,\, \frac{K}{\sqrt{\beta}} (\<{B^H,f_1},\ldots,\<{B^H,f_k})\,.
\end{align*}
\end{theorem}

We now have all the necessary ingredients (plus with the key gradient estimates which shall be obtained in the next Section) to complete the proof. 
\medskip
\medskip

\ni
{\em Proof.}

Let us fix $\alpha\in [2,3)$. The proof is divided into the following different steps.

\bnum
\item \textbf{Caffarelli-Silvestre extension of the discrete fractional Laplacien.}
We use here the Bessel random walks from Section \ref{s.Bessel}. This allows us to view the field $\Psi_n$ as the Gaussian free field  on $\calD$ (with inhomogeneous conductances $\mathbf{a}$) conditioned to take integer values on the $1d$ line $\Lambda_n \times \{0\}$. 

\item \textbf{Identification of the suitable Laplace transforms.}
If we are given test functions $f_1,\ldots,f_k$ in $C_c(D)$ where $D=(-1,1)$. To characterise our limiting distribution, it is sufficient to  
show that for any $\lambda_1,\ldots, \lambda_k \in \R$, one has the following convergence (for small enough $\beta$), 
\begin{align*}\label{}
&\lim_{n\to \infty}\EFK{\Lambda_n,J_\alpha,\beta}{dGC,0}{\exp\left( \lambda_1 \<{\mathbf{H}_n,f_1} + \ldots + \lambda_k \<{\mathbf{H}_n,f_k}\right)} \\
& = \EFK{D,\beta}{H,0}{\exp K^2 \left( \lambda_1 \<{B^H,f_1} + \ldots + \lambda_k \<{B^H,f_k}\right)} \\
& = \lim_{n\to \infty}\EFK{\Lambda_n,J_\alpha,\beta/K^2}{GFF,0}{\exp\left( \lambda_1 \<{h_n,f_1} + \ldots + \lambda_k \<{h_n,f_k}\right)}\,,
\end{align*}
where the second equality is the result of Proposition \ref{pr.sheffield} from \cite[proposition 12.2]{lodhia2016fractional}. 
It will thus be sufficient in order to prove our main result to have an asymptotic comparison between the first and third line in the above displayed equality. 

It will be useful below to ``unscale'' the test functions $(f_i)$ and to work on domains $D_n= \calD \setminus L_n$ with one-unit spacing.
For this notice that for any given $f_1,\ldots, f_k \in C_c(D)$ and any scalars $\lambda_1,\ldots, \lambda_k$, the Laplace transform we are after is 
\begin{align*}\label{}
\lambda_1 \<{\mathbf{H}_n,f_1} + \ldots + \lambda_k \<{\mathbf{H}_n,f_k} 
&\overset{\text{By Def.}}= \sum_{t \in \frac 1 n \Z} \frac 1 {n^{1+H}} \Psi_n(nt) \sum_{i=1}^k \lambda_i f_i(t) \\
& = \sum_{ x\in \Lambda_n  } \frac 1 {n^{1+H}} \Psi_n(x) \sum_{i=1}^k \lambda_i f_i(\frac x n) \\
& =  \<{\Psi, g}\,, 
\end{align*}
for the test function $g: \Lambda_n \to \R$ defined by 
\begin{align}\label{e.defg}
g(x):=  \frac 1 {n^{1+H}} \sum_{i=1}^k \lambda_i f_i(\frac x n)\,.
\end{align}

Notice that since $f_i\in \calC_c(D), \forall 1\leq i \leq k$, one can find some $\delta>0$ such that 
\begin{align*}\label{}
\mathrm{Supp}(g) \subset \{-(1-\delta) n , \ldots, (1-\delta) n \}\,.
\end{align*}

\item \textbf{Ginibre.}
Ginibre inequality (Proposition \ref{pr.ginibre}) immediately gives us the easy bound, namely that the Laplace transform of the integer-valued Gaussian is smaller than the Gaussian one. 
We thus only need to find a matching lower bound.  

\item \textbf{Inequality for Gaussians on lattices \cite{regev2017inequality}}. By using this inequality (which we stated in Theorem \ref{th.regev}), we obtain that for any scale parameter $M$, the Laplace  transform of a test function $g\in \calC_c(D)$ satisfies
\begin{align*}\label{}
\EFK{H_n^{(M)},\mathbf{a},\beta}{dGC,0}{e^{\<{\Psi, g}}}  \leq \EFK{D_n,\mathbf{a},\beta}{dGC,0}{e^{\<{\Psi, g}}}\,.
\end{align*}
(Recall the definition of the sub-domain $H_n^{(M)}$ from~\eqref{e.Hn} and Figure \ref{f.Hn}). 

To see why this holds from Theorem \ref{th.regev}, this corresponds to setting infinite coupling constants on each edge of the graph $D_n \setminus H_n^{(M)}$. 
 
 \item \textbf{Small loss in the Laplace transforms as the smoothing effect $M\to \infty$.}
 
 In the steps to follow, we will show if $M$ is chosen large enough (in particular one needs that $\frac 1 M \leq \frac 1 4 \dist(\mathrm{supp(f)}, \p D)$) then,  as $n\to \infty$, one has 
 \begin{align}\label{e.convenient1}
\lim_{n\to \infty}
\frac
{
\EFK{H_n^{(M)},\mathbf{a},\beta}{dGC,0}{e^{\<{\Psi, g}}}  
}
{
\EFK{H_n^{(M)},\mathbf{a},\beta}{GFF,0}{e^{\<{\varphi, g}}} 
}
=1\,.
\end{align}
Assuming the above holds, in order to conclude the proof of the main Theorem, we need to prove purely on the ``Gaussian free field'' side that
 \begin{align}\label{e.convenient2}
\liminf_{M\to \infty}  
\liminf_{n\to \infty}
\frac
{
\EFK{H_n^{(M)},\mathbf{a},\beta}{dGC,0}{e^{\<{\varphi, g}}}  
}
{
\EFK{D_n,\mathbf{a},\beta}{GFF,0}{e^{\<{\varphi, g}}} 
}
=1\,.
\end{align}
It is clear that such a limit is $\leq 1$ (by direct comparison of the quadratic forms involved).

To prove the other direction, notice that since both the numerator and denominator are the Laplace transforms of Gaussian vectors, there are given in terms of their variance and thus, we  only need to show that 
\begin{align*}\label{}
\liminf_{M\to \infty}  
\liminf_{n\to \infty}
\left(\EFK{H_n^{(M)},\mathbf{a},\beta}{dGC,0}{\<{\varphi, g}^2}  
-
\EFK{D_n,\mathbf{a},\beta}{GFF,0}{\<{\varphi, g}^2} \right)
=0\,.
\end{align*}

We proceed in two steps:

\ni
$(i)$ First, we show that both  variance are bounded as $n\to \infty$. By the same monotony as above, it is sufficient to show that 
\begin{align*}\label{}
\limsup_{n\to \infty}\EFK{D_n,\mathbf{a},\beta}{GFF,0}{\<{\varphi, g}^2}  < \infty\,.
\end{align*}

Since our test function $g$ is supported on the line $\Lambda_n \times \{0\} \subset D_n$, we find by the definition of the test function $g$ from \eqref{e.defg}
\begin{align*}\label{}
\EFK{D_n,\mathbf{a},\beta}{GFF,0}{\<{\varphi, g}^2} 
& = \sum_{x,y \in \Lambda_n} \frac{f(\frac x n)}{n^{1+H}} \frac  {f(\frac y n)}{n^{1+H}}  G_{D_n,\mathbf{a}}(x,y) \\
& \leq O(1) \|f\|_\infty^2\,,
\end{align*}
using the upper bound from Proposition \ref{pr.Green1}. Note that we do not use the fact that $f$ is compactly supported here. 

\medskip

\ni
$(ii)$ We now use the fact $(f_i)_{1\leq i \leq k}$ are compactly supported in $D=(-1,1)$. This implies in particular that there exists $\delta>0$ such that our  test function $g=g_n$ is supported for all $n\geq 1$ in $\{-(1-\delta)n, \ldots, (1-\delta)n \} \times \{0\} \subset D_n$. 
We fix such a $\delta$ once and for all in what follows. (This will affect constants in the propositions which will be needed next, namely Proposition \ref{pr.DnHn}, Proposition \ref{pr.Green1} and Proposition \ref{pr.Green2}). 

Of immediate interest, we apply Proposition \ref{pr.DnHn}, this gives us for $M$ large enough:
\begin{align*}\label{}
\EFK{H_n^{(M)},\mathbf{a},\beta}{GFF,0}{\<{\varphi, g}^2} 
& = \sum_{x,y \in \Lambda_n} \frac{f(\frac x n)}{n^{1+H}} \frac  {f(\frac y n)}{n^{1+H}}  G_{H_n^{(M)},\mathbf{a}}(x,y) \\
& \geq (1-\xi_\delta(M)) \sum_{x,y \in \Lambda_n} \frac{f(\frac x n)}{n^{1+H}} \frac  {f(\frac y n)}{n^{1+H}}  G_{D_n,\mathbf{a}}(x,y) \\
& \geq (1-\xi_\delta(M)) \EFK{D_n,\mathbf{a},\beta}{GFF,0}{\<{\varphi, g}^2} \,.
\end{align*}
Since at fixed $\delta>0$, we have from Proposition \ref{pr.DnHn} that $\lim_{M\to \infty} \xi_\delta(M) =0$ and since we checked that both variances remain bounded as $n\to \infty$, this ends the proof that the limit~\eqref{e.convenient2} holds.

\medskip
For the remaining steps, we may then fix some arbitrary large parameter $M$. We will only work in the sequence of ``smooth'' domains $H_n^{(M)}$ (defined in~\eqref{e.Hn}, see also Figure \ref{f.Hn}). This is important as our gradient estimates from the next Section \ref{s.GRAD} only apply in such domains. 

Our goal from then on is to establish that the limit~\eqref{e.convenient1} holds, which would then finish the proof. 

\medskip

As we did in $2d$ domains in Section \ref{s.proof2}, we will run a Coulomb-gas expansion in the spirit of Fröhlich-Spencer in the ``smooth'' domain $H_n^{(M)}$ with Coulomb charges restricted to the $1d$ line $\Z\times \{0\} \cap H_n^{(M)}$. 
The proof, modulo the Gradient estimates from the next Section, will be close to what we have done in Section \ref{s.proof2}. We will shortly explain how to adapt Fröhlich-Spencer to the geometry of $H_n^{(M)}$. This will be mostly harmless except one key point in Wirth's argument to control the Dirichlet rooting. 
\medskip

\item \textbf{Setting up the Coulomb-gas expansion.}
As in Section \ref{s.proof2}, we assign Coulomb charges to each point on the line $\Z \times \{0\} \cap H_n$. We work with Dirichlet boundary conditions on $\p H_n$. 

Recall that before conditioning the field to be integer-valued, we work with the GFF in inhomogeneous conductances:
\begin{align*}\label{}
\exp\left(-\frac \beta 2  \sum_{i\sim j \in \calD} a_{i,j} (\varphi_{\alpha,n}(i) - \varphi_{\alpha,n}(j))^2 \right) 1_{\varphi_{\alpha,n} \equiv 0 \text{  on  } L_n}
\end{align*}
Since we localised the analysis, we shall then work with 
\begin{align*}\label{}
\exp\left(-\frac \beta 2 \sum_{i\sim j \in \calD} a_{i,j} (\varphi_{\alpha,n}(i) - \varphi_{\alpha,n}(j))^2 \right) 1_{\varphi_{\alpha,n} \equiv 0 \text{  on  } \p H_n}
\end{align*}

\item \textbf{Adapting Fröhlich/Spencer -- Wirth expansion to the geometry of the slit domain and the inhomogeneous conductances.}

Before looking for a quantitative version which is precise enough to lead to the invariance principle and to identifies $\beta_{eff}(\beta)$, we first need to run the  Fröhlich/Spencer -- Wirth Coulomb gas expansion. Namely the two main steps we described previously, i.e. the charge expansion and then the spin-wave analysis to reduce the ``activities'' $z(\rho,\beta,\calN)$. 
 
We claim, that the same proof goes through, the only parts which require some attention are the following two points:

(a). First,  In the case of a slit domain, we need to adjust the notion of distance used when we decide to combine two charges $\rho_1$ and $\rho_2$. 
The metric we should use here is the induced graph distance inside $D_n$ ($H_n$ in fact), rather than the euclidean distance in the plane. Indeed otherwise a charge above the ``slit'' would be considered close to a charge just below the ``slit'' and this would make the algorithm behind the gas expansion fail.  
\smallskip

(b). We also need to work with spin-waves in the presence of the inhomogeneous field of conductances $\mathbf{a}$. Here we claim that this is harmless for the estimations of the activities $z(\rho,\beta,\calN)$. The reason why one can accommodate with such conductances is the property that they vanish as the vertical directions goes to infinity (as $d^{-s}$ if $d$ is the distance to the line $\Z \times \{0\}$). The difficulty with the analysis of spin waves is the fact they may carry a lot of Dirichlet energy. Those conductances only make these Dirichlet energies smaller and this is sufficient to run the Fröhlich-Spencer expansion.  

\item \textbf{Obtaining the up to constants estimates from Theorem \ref{th.main}.}

At this stage, we may already obtain the up to constants estimates~\eqref{e.UpToC} from Theorem \ref{th.main} by using the power of this Coulomb gas expansion together with Proposition \ref{pr.Green2}. 
Given Theorem \ref{pr.Green2}, this is very similar as in \cite{FS,wirth2019maximum}. 

\item\textbf{Extending Wirth's argument all the way to the boundary. }

Since we will need to control the Dirichlet energy of Green functions all the way the left/right boundaries of $H_n$, we need to extend Wirth's analysis all the way to the boundary. This can be done exactly as we did in Section \ref{ss.NN}.

\item \textbf{Controling the lack of Gaussian behaviour via the Dirichlet energy of the line.}

Now, given our test function $g$ (defined in~\eqref{e.defg}), we proceed as in the Proof of Theorem \ref{th.line2}. Namely we define the function 
\begin{align}\label{e.sigma4}
	\sigma:= \frac 1 \beta [- \Delta_{H_n,\mathbf{a}}]^{-1} g 
\end{align} 
 As in the computation~\eqref{e.FRACT}, the deviation (in the sense of Laplace transforms) between the integer-valued case and the Gaussian case is controlled by the Dirichlet energy of $\sigma$ on the $1D$ line, namely 
\begin{align*}\label{}
 \sum_{j\sim l \in \Z \times \{0\}} (\sigma_j - \sigma_l)^2 \,. \nn 
\end{align*}
As in~\eqref{e.FRACT}, one can express this by using in the present setting the Green function $G_{H_n, \mathbf{a}, M}$. 

\item \textbf{Conclusion.}
The next Section \ref{s.GRAD} is precisely providing us with the needed analytic control on the gradients $\nabla_i G_{H_n, \mathbf{a}, M}(x,i)$ in order to quantitatively upper bound the quantity $\sum_{j\sim l \in \Z \times \{0\}} (\sigma_j - \sigma_l)^2$. As such we obtain our invariance principle together with the {\em invisibility of integers} phenomenon.
\enum

\qed

\section{Controlling the Dirichlet Energy supported on a line when $\alpha<3$}\label{s.GRAD}



This section provides the key technical estimates in order to prove the invariance principle towards an $H=H(\alpha)$-fractional Brownian motion when $\alpha\in[2,3)$. It shows that most of the Dirichlet Energy of the $2d$ Bessel random walk is spread over $\calD\setminus (\Z \times \{0\})$.

We consider the Bessel random walk on the diamond graph $\calD$ as defined in the Section \ref{s.Bessel}. The proof also works for Bessel walks on $\Z^2$ discussed in Appendix \ref{s.BesselZ}.
Recall  that these walks are parametrised by the parameter $s>-1$ and that the relationship with $\alpha$ reads as follows:
\begin{align*}\label{}
s=\alpha-2\,.
\end{align*}

From then on, we will fix $s\geq 0$, some small $\delta>0$ as well as some large  parameter $M$ as in Subsection \ref{ss.compact}. 
Recall from this subsection that $H_n=H_n^{(M)}$ is the bounded ``smoothed'' version of the slit domain defined in~\eqref{e.Hn}. Recall also  that
$G_{H_n,\mathbf{a},M}$  is the Green function for the random walk in conductances $\mathbf{a}=\mathbf{a}(s)$ and killed when first hitting either the $\frac 1 M n$ boundary $\calL_n$ of the ``slit'' $L_n$ or the square of radius $nM$, $S_{nM}$.  In what follows, for simplicity,  we will only denote it by $G_{H_n,\mathbf{a}}$.

 Our first main estimate is to control the gradient of the Green function $y \mapsto G_{H_n,\mathbf{a}}(x,y)$  uniformly in $x$ and $y$ on the line $\Lambda_n=\Lambda_n\times \{0\}$ and at macroscopic distance from $L_n$. 
 \begin{proposition}\label{pr.Grad}
 For any $s\in [0,1)$ any $\delta>0$ and any $\eps>0$, there exists a constant $c>0$ such that for all $n\geq 1$, $x\in \{-(1-\delta)n , \ldots, (1-\delta)n \}$ and $y\in \{-(1-\tfrac \delta 2)n , \ldots, (1- \tfrac \delta 2)n \}$  one has 
%

\begin{align*}\label{}
 |G_{H_n,\mathbf{a}}(x,y+1)-G_{H_n,\mathbf{a}}(x,y)|  \leq   c\,  \left( \frac 1 {|x-y|+1} \right)^{1-s - \eps}\,.
\end{align*}
 \end{proposition}

We will then prove the Proposition below which handles the case where points $y$ are close to the killing boundary $L_n$. 
It will rely on a version of {\em Beurling's estimate} adapted to the present setting of discrete Bessel walks when $s\geq 0$. (See Remark \ref{r.tedious2}). 

\begin{proposition}\label{pr.Grad2}
 For any $s\in [0,1)$ any $\delta>0$
 and any $\eps>0$, there exists a constant $c>0$ such that for all $n\geq 1$, $x \in \{-(1-\delta)n , \ldots, (1-\delta)n \}$  and $y\in\{-n,\ldots, -(1- \tfrac \delta 2) n\} \cup \{(1-\tfrac \delta 2) n, \ldots,  n\}$ one has 
\begin{align*}\label{}
 |G_{H_n,\mathbf{a}}(x,y+1)-G_{H_n,\mathbf{a}}(x,y)|  \leq  
   c\, \left( \frac 1 {n} \right)^{1-s - \eps} \,.  
\end{align*}
 \end{proposition}

Using Proposition \ref{pr.Grad} and \ref{pr.Grad2}, one readily obtains the estimate below which is the main technical input for the proof of "invisibility of integers" when $\alpha\in[2,3)$ in Section \ref{ss.TheProof}.
\begin{corollary}
If $0\leq s<1$, then for any base point $x$ at macroscopic distance from the boundary, the Dirichlet energy of $y\in H_n \mapsto G_{H_n,\mathbf{a}}(x,y)$ restricted to the line $\{-n,\ldots,n\}\times \{0\}$ is controlled by 
\begin{align*}\label{}
& \sum_{i\in \Z\times \{0\} \cap \Lambda_n}  | \nabla_i G_{H_n,\mathbf{a}}(x,i)|^2 \leq O(1) \sum_{k=1}^{n} \frac 1 {k^{2-2s - 2\eps}}\,,
\end{align*}
\end{corollary}
%
%

\smallskip

\subsection{Proof of Proposition \ref{pr.Grad}.}\label{ss.Grad1}
By using the symmetry property of the Green function 
(N.B. which holds on the middle line, otherwise this is only correct modulo the action of $D_\mathbf{a}^{-1}$), we may write the above gradient of two Green functions as follows: 
\begin{align*}\label{}
G_{H_n,\mathbf{a}}(y+1,x)-G_{H_n,\mathbf{a}}(y,x)
&= \sum_{t\geq 0}^\infty   \EFK{}{}{1_{T_n^{1}>t}1_{(X_t^{1},Y_t^{1})=(x,0)} - 1_{T_n^0>t}1_{(X_t^{0},Y_t^{0})=(x,0)} } \,,
\end{align*}
for any coupling of a Bessel walk $(X_t^{1},Y_t^{1})$ starting at $(y+1,0)$ with a Bessel walk $(X_t^{0},Y_t^{0})$ starting at $(y,0)$ and where $T_n^{1}$ (resp. $T_n^0$) are the stopping times of these walks when first hitting the boundary of the bounded slit domain $\p H_n$.

The first (classical) coupling which comes to mind is as follows: the vertical coordinate (in the diamond graph $\calD$) are identical for both walks while the horizontal coordinates $X_t^{1}$ and $X_t^0$ are coupled via independent walks until 
\bi
\item Either one of the two walks exits $\p H_n$ in which case it stops and the second walk keeps moving independently. In this case, we let $\tau:= \infty$  
\item Or they merge before any of the two walks visits $\p H_n$, in which case they continue the same journey from then on. In this case we let $\tau\in \N$ to be the stopping time when they first merge. 
\ei


This classical coupling has the following slightly unpleasant property (which would not be visible in the continuum limit): the marginal $\{(X_s^0)_{s\geq 0}\}$ is \underline{not} independent of the value of the stopping time $\tau$. (This is because the process $s\mapsto X_s^{1}-X_s^0$ is not independent of $s\mapsto X_s^0$).

To overcome this small issue, we proceed slightly differently below for the coupling of the horizontal coordinates (the vertical coordinates follow the exact same trajectory). Let us also fix some small $\hat \eps=\hat\eps(\eps) < \eps$ whose value will be fixed further below. 

\smallskip
\ni
We shall break the series 
\begin{align*}\label{}
\sum_{t\geq 0}^\infty   \EFK{}{}{1_{T_n^{1}>t}1_{(X_t^{1},Y_t^{1})=(x,0)} - 1_{T_n^0>t}1_{(X_t^{0},Y_t^{0})=(x,0)} } 
\end{align*}
into four groups as outlined below. 

\textbf{First group. $t\leq |x-y|^{2-\hat \eps}$.}
By choosing $c$ sufficiently large, it is immediate to see that the times $t$ in this first group contribute at most $e^{-c^{-1}|x-y|^{\hat \eps}}$ to the total sum. (This is very similar to the analysis carried in Section \ref{s.Bessel}).   

\textbf{Second group. $|x-y|^{2- \hat \eps} \leq t \leq n^{2- \hat \eps}$.}

For each such time $t$, we first apply the standard coupling argument until time $\tfrac t 2 $. If we disregard the effect of the boundary $\p H_n$, the probability that the coupling has not succeeded (i.e. that $\tau> T_0$) is easily seen to be upper bounded by $O(1) t^{-1/2}$. 
Now by taking into account the presence of the boundary $\p H_n$, the event $\{\tau > \tfrac t 2 \}$ can be written as the union of 2 events: $A)$ exactly one of the 2 walks (starting from $y+1$ or $y$) visits $\p H_n$ before time $\tfrac t 2 $  (recall we assumed $y$  to be at distance at least $\tfrac \delta 2 n$ from $\p H_n$). And $B)$ two independent walks on $\tfrac \Z 2$ (without boundaries) and starting at $y+1$ and $y$ do not merge by time $\tfrac t 2 $ (the horizontal walks are along $\tfrac 1 2 \Z$ because of the projection of the diamond graph $\calD$). The probability of the event $B)$ is less than $O(1) t^{-1/2}$ as already discussed while the probability of the event $A)$ is easily seen to be less than $\exp(-c^{-1}n^{\hat \eps})$ by the same observation as above. All together, with probability at most $O(1) t^{-1/2}$, at least one of the particules $X_{\tfrac t 2 }^1$ and $X_{\tfrac t 2 }^0$ is still "active" (did not reach the boundary) and they are not coupled yet. 
Conditioned on this event, we find $X_{\tfrac t 2 }^1$ and $X_{\tfrac t 2 }^0$ at some random conditional (different) positions. (We may need to follow only one active particle below but this leads to the same upper bound). From then on, we let the trajectories behave independently of each other (we do not even dare to check whether they merge or not). The only observation we need is that uniformly in the positions of $X_{\tfrac t 2 }^1$ and $X_{\tfrac t 2}^0$, the probability that either of these walks terminate precisely at $y$ at time $t$ is upper bounded (again by Local CLT) by $O(1) t^{-1/2}$. Here we use the fact that the walks have a remaining time $\tfrac t 2$ to run between their conditioned position $X_{\tfrac t 2 }^1, X_{\tfrac t 2 }^0$ and their "target" position  $y$ at time $t$.  
Note that for this part of the analysis, the boundary $\p H_n$ is only helping us, while for the coupling part ($s< \tfrac t 2$), it was creating interferences against us. 

We thus find that each time $t$ in this regime contributes at most 
\begin{align*}\label{}
O(1) t^{-1/2} \times  t^{-1/2} \times \Pb{Y_t^1=Y_t^0=0}\,.
\end{align*}
The last term is controlled by Lemma \ref{l.Renew}. 
This gives us a contribution for the time $t$ bounded by 
\begin{align*}\label{}
O(1) t^{-1}  \times t^{-(1-s)/2}\,.
\end{align*}
The sum over all times in this second group is therefore upper bounded by 
\begin{align*}\label{}
& \sum_{|x-y|^{2-\hat \eps} \leq t \leq n^{2-\hat \eps}} O(1) t^{-1} \times t^{-(1-s)/2}  \leq O(1) (|x-y|+1)^{s-1 + \tfrac {\hat \eps} 2 (1-s)}\,.
\end{align*}

We now introduce our third group of times contributing to the gradient of Green functions:

\textbf{Third group. $n^{2- \hat \eps} \leq t \leq n^{2 + \hat \eps}$.}

The reason why we need this third group is because it is more difficult to control the probability that the coupling fails before time $\tfrac t 2$. We managed to have a good control for times in the second group because particles did not have the time (except with a stretch exponentially small cost) to reach $\p H_n$. This is no longer the case here. 

To deal with this technical issue, we only try to couple the two walks until $T_0:= \tfrac 1 2 n^{2-\hat \eps}$. Since $T_0 \ll n^2$, it still requires a stretch exponential cost to reach $\p H_n$ before time $T_0$. By the same analysis as in the second group but applied with a fixed $T_0$ rather than with $\tfrac t 2$, we get that with probability at most $O(1) T_0^{-1/2}$, we still have one or two active particles. Now, since $T_0$ was chosen as half of $n^{2-\hat \eps}$, for any time $t$ in this third group, we still have at least $\tfrac t 2$ steps before reaching the target $y$. The same local CLT bound gives us that uniformly in the position of the active(s) particles at time $T_0$, this gives an additional costs of $t^{-1/2}$. Finally, the vertical direction provides in the exact same way an additional $t^{-(1-s)/2}$ term. 

As such, we obtain
\begin{align*}\label{}
& \sum_{n^{2-\hat \eps} \leq t \leq n^{2+\hat \eps}}   \EFK{}{}{1_{T_n^{1}>t}1_{(X_t^{1},Y_t^{1})=(x,0)} - 1_{T_n^0>t}1_{(X_t^{0},Y_t^{0})=(x,0)} }  \\ 
& \leq O(1) \sum_{n^{2-\hat \eps} \leq t \leq n^{2+\hat \eps}} n^{-\tfrac 1 2 (2-\hat \eps)} t^{-1/2} \times t^{-(1-s)/2} \\
& \leq O(1) n^{-1 + \tfrac 1 2 \hat \eps} n^{2+\hat \eps}* n^{-(2-\hat \eps)*(1-\tfrac s 2)} \\
& \leq O(1) n^{s-1 + (\tfrac{5-s} 2) \hat \eps} 
\end{align*}

\textbf{Fourth group. $t > n^{2 + \hat \eps}$.}
In this group, we proceed as follows:  we do not even need to dare coupling the particles in an efficient way. We may just let both particles undergo an independent simple random walk along the horizontal axes. 
The only easy observation we make (similarly as in the proof of Proposition \ref{pr.Green1}) is that there exists $a=a(M)>0$ so that uniformly in $x$, for any $t\geq n^2$, 
\begin{align*}\label{}
\FK{}{x}{X^1_u \text{ stays inside }[-M\cdot n, M \cdot n]\,\,,\,\, \forall 0\leq u \leq t} \leq a \exp(-a \frac t {n^2})\,.
\end{align*}
(Recall that $M$ is the scale parameter in the definition of the bounded slit domain $H_n=H_n^{(M)}$ in~\eqref{e.Hn}). 
This is again because, on each consecutive time interval of length $n^2$, the horizontal walk has a positive probability (depending on $M$) to exit the ball $B_{\|\cdot\|_\infty}(0,M \cdot n)$. 

As such this fourth group contributes at most 
\begin{align*}\label{}
 \sum_{ t > n^{2+\hat \eps}}  a \exp(- a t/n^2) \leq  O(1) n^2 e^{- a n^{\hat \eps}}\,,
\end{align*}
which is negligible.

By summing the contributions from each of the four groups above and taking $\hat \eps = \hat \eps(\eps)$ sufficiently small, this ends our proof of Proposition \ref{pr.Grad} in the case where the "target" point $y$ is at macroscopic distance from $\p \Lambda_n$.  \qed

\begin{remark}\label{r.tedious1}
The previous fourth group of (large) times $t$ is the reason why we need to compactify our space in Subsections \ref{ss.compact} and \ref{ss.TheProof} before applying the Coulomb gas expansion. Indeed without this technical step, we would need to analyse the efficiency of the coupling of two walks for large times but in a complicated geometry given by a slit domain $D_n$. In particular if walks survive for a long time $\gg n^2$ and did not couple yet, this is not straightforward to analyze, as their conditional position inside the slit domain is complicated. 
\end{remark}

\subsection{Proof of Proposition \ref{pr.Grad2}.}\label{ss.Grad2}
$ $

Let us now deal with the second situation where the target point $y$ is now at small distance (less than $\frac \delta 2 n$) from the slit $L_n$ and is at far distance from $x$ (at least $\frac \delta 2 n$ given our assumptions).


Let us call $d:= \dist(y, \p \Lambda_n)$ which may take any value in  $\{1,\ldots, \frac \delta 2 \, n\}$. 

As before we need to upper bound using suitable (possibly $t$-dependent) couplings the series 
\begin{align*}\label{}
\sum_{t\geq 0}^\infty   \EFK{}{}{1_{T_n^{1}>t}1_{(X_t^{1},Y_t^{1})=(x,0)} - 1_{T_n^0>t}1_{(X_t^{0},Y_t^{0})=(x,0)} }\,. 
\end{align*}

One advantage of this present case is that we assumed $|x-y| \geq \frac \delta  2 n$. One additional difficulty on the other hand is that the walks may start very close to $\p \Lambda_n$ which may affect the efficiency of the coupling. 
We will then decompose the sum in the following groups:
\bnum
\item \textbf{Small times.} $t\leq n^{2-\hat \eps}$. For these times it requires a stretch exponential cost ($e^{-c n^{\hat \eps}}$) for particles to reach $x$ without any considerations of couplings or boundary effects. 
\item \textbf{Main contributing times.} $n^{2-\hat \eps} \leq t \leq n^{2+\hat \eps}$.

For these times, we proceed differently as in Section \ref{ss.Grad1}. We only take benefit of the coupling between times $0$ and $d^{2-\hat \eps}$. 
The reason for such a time scale is that for times much smaller than $d^{2 -\hat \eps}$, the walks do not feel the presence of the boundary, not the presence of the target point $x$ which is at distance at least $\frac \delta 2 n$. 
With probability $1 - O(1)d^{-1+\tfrac 1 2 \hat \eps} - e^{-c d^{\hat \eps}}$, the walks have coupled and the contribution to the gradient of the Green functions is zero (we did not have time to visit $x$ yet). 
We are thus left with an event of probability $O(1) d^{-1+\tfrac 1 2 \hat \eps}$ where at least one particle is still active and needs to reach $x$ at time $t$. 

Now, we need to take into account the effect of the nearby boundary. 
We argue as follows: up to another stretch exponentially small term $O(e^{-c d^{\hat \eps}})$, the remaining active particles are at distance at most $d$ from $y$ and as such at distance at most $2d$ from $\p H_n^{(M)}$. (Recall the definition of this domain from \ref{e.Hn} and see Figure \ref{f.Hn}).  
In order to reach $x$, the particles will need to avoid this boundary on a long-time period. 
We only quantify this effect on the time interval $[d^{2-\bar \eps},  \tfrac t 2 \wedge n^{2-\tfrac {\hat \eps} 2}]$. On this time interval, we use the fact the boundary of $H_n^{(M)}$ looks exactly like a full vertical line at distance $O(d)$ from the position at time $d^{2- \bar 
\eps}$. (see Figure \ref{f.Hn}). This is because the height is this wall is $\frac 1 M n$ and because there is only a stretch exponentially small probability by time $n^{2- \frac {\hat \eps} 2}$ for the walks to reach the corners of $H_n$. As such, by standard Gambler ruin estimates, the probability that active particles do not intersect $\p H_n$ is upper bounded by 
\begin{align*}\label{}
O(1) \frac d {n^{1- \frac{\hat \eps} 4}} \wedge 1\,.  
\end{align*}

\begin{remark}\label{r.tedious2}
This estimate in fact plays a crucial role in the analysis. 
If we had not restricted the chain to leave in a smoother domain $H_n^{(M)}$ in Section \ref{ss.compact}, we would end up here with a much less good control, as the walks would not see a ``wall'' of heigth $\frac 1 M n$ but would see instead the slit half-line $L_n$. Using a Beurling type estimate adapted to $s$-Bessel walks, this would give us an upper bound of the order 
\begin{align*}\label{}
 \sqrt{\frac d {n^{1- \frac{\hat \eps} 4}} \wedge 1} \text{    instead of   }  
\frac d {n^{1- \frac{\hat \eps} 4}} \wedge 1 \,.
\end{align*}
Such a bound would simply not allow us to conclude when  $s\geq 1/2$ (i.e. $\alpha\geq \frac 5 2$) !
\end{remark}
\medskip

At this stage, we only need to check the horizontal positions of the walks here so that one does not impose any conditioning on the behaviour of the vertical position. (Up to the stretch exponential probability mentioned above, the horizontal coordinates in this time window are able to tell us whether the $2d$ random walks exists the domain or not). 
 
\medskip
 
Now, as we did in Section \ref{ss.Grad1}, we let evolve (freely) the particles up to time $t/2$ and we condition on the horizontal positions of the active particles at time $\tfrac t 2 $ (there could be one or two such active particles). In the remaining $\tfrac t 2 $ steps, by local CLT theorem and uniformly on the position of $x$ (close to $\p H_n$ or not), we get another term $O(1) t^{-1/2}$. Finally the vertical process does not care about the horizontal vicinity of $\p H_n$ and gives us exactly the same contribution as in Section \ref{ss.Grad1}. 

As such this interval of times gives a contributions smaller than 
\begin{align*}\label{}
&\sum_{n^{2-\hat \eps} \leq t \leq n^{2+\hat \eps}} O(1) d^{-1+\tfrac 1 2 \hat \eps} *\left( \frac d {n^{1-\tfrac {\hat \eps} 4}} \wedge 1\right) *  t^{-1/2} *  t^{-(1-s)/2} \\
&\leq O(1) n^{2+\tfrac {\hat \eps} 2 }  d^{\frac {\hat \eps} 2} \frac 1 {n^{1 - \tfrac {\hat \eps} 4} * n^{1-\tfrac {\hat \eps} 2} *  n^{1-  s  - \hat \eps(\tfrac {1-s} 2)} } \\
& \leq O(1) n^{s- 1 + \frac {11} 4 \hat \eps - \frac s 2 \hat \eps} \\
& \leq O(1) n^{s-1 + \frac {11} 4 \hat \eps}\,.
\end{align*}

\item \textbf{Large times $t > n^{2 + \hat \eps}$.}

This case is handled exactly as in the fourth group of times in Section \ref{ss.Grad1}, i.e. by exploiting the fact our smoother domain $H_n=H_n^{(M)}$ is of diameter less than $M \cdot n$ (see Figure \ref{f.Hn}).

%
%
%
\enum

The above analysis thus finishes the proof of Proposition \ref{pr.Grad2} when $y$ is close to boundary and $x$ is at distance $\delta n$ from $y$. Notice that this second case requires us to lower a bit the value of $\hat \eps$.   \qed

%

\section{Proof of Propositions \ref{pr.comp1} and \ref{pr.comp2}}
\label{s.visible}

The cases $\alpha>3$ in $d=1$ and $\alpha>4$ in $d=2$ share the following common feature: in each case the random walks on $\Z$ (resp. $\Z^2$) with long-range weights $\frac 1 {\mathrm{dist}^\alpha}$ have a \uline{finite} second moment. 

This immediately calls for the  following well-known Lemma from \cite{varopoulos1986theorie}, see Lemma 2.1. from \cite{pittet2000stability}.
We state it in the particular case of random walks on the lattice $\Z^d$, but this holds in much more general settings.
\begin{lemma}[Lemma 2.1. in \cite{pittet2000stability}]\label{l.varo}
If one considers a random walk on $\Z^d$ with long-range i.i.d increments which have a finite second moment, then there exists a constant $C>0$ such that 
\begin{align*}\label{}
\frac 1 C (-\Delta_{\Z^d}) \leq (-\Delta^{LR}) \leq C (-\Delta_{\Z^d})\,,
\end{align*}
where $\Delta^{LR}$ is the associated long-range Laplacian and where the inequalities are understood in the sense of quadratic forms. 
\end{lemma}

\ni
{\em Proof of Proposition  \ref{pr.comp1} and \ref{pr.comp2}.} 
Let us give the (short) proof in the case of the long-range integer-valued GFF with Dirichlet boundary conditions outside $\Lambda_n^2 \subset \Z^2$.
When $\alpha>4$, it is easy to check that the induced $2d$ random walk has finite $L^2$ moment. 

On can then apply the above Lemma. The inequality of interest is the less obvious one, namely 
\begin{align*}\label{}
(-\Delta^{J_\alpha}) \leq C_\alpha (-\Delta_{\Z^d})\,
\end{align*}
where the coupling constants of the induced random walk are, say,  $J_\alpha(r) = r^{-\alpha}$ or any other such decaying kernel. 

Now, we may use the (non-trivial) Theorem \ref{th.regev} to deduce that the Laplace transform of the (in principle complicated) long-range integer-valued field $\Psi_n^{d=2,J_\alpha}$ are \uline{bounded from below} by the Laplace transform of the more classical nearest-neighbor integer-valued GFF. 

We thus extract, thanks to this comparison principle,  $\log n$ variance bounds for these long-range models out of the $\log n$ variance bounds for the nearest neighbour case (\cite{FS,wirth2019maximum,RonFS, lammers2022dichotomy, park2022central}). 

Still, note that we do not obtain an invariance principle using such comparison techniques.

When $d=2$ and $\alpha=\alpha_c=4$, the $\log \log(n)$ upper bound on the variance in Proposition \ref{pr.comp2} follows by combining Ginibre inequality, Theorem \ref{th.regev} as well as the $\log \log(n)$ upper  bounds on the effective resistance of $\alpha=\alpha_c$ walks on $\Z^2$ from the work \cite{caputo2009recurrence} (N.B. matching lower bounds are also obtained in \cite{caputo2009recurrence} and the walks start being transient as soon as $\alpha<\alpha_c=4$, see also the work \cite{baumler2023recurrence}).

See Open Problem \ref{OP1} below, where we discuss further the case $\alpha_c=4$ in $d=2$ which seems particularly interesting to us. 
\medskip

Finally, the same idea works also in the $1d$ case, i.e. for the discrete Gaussian Chain when $\alpha>3$. This allows us to bound from below the fluctuations of the $\alpha>3$-discrete Gaussian chain  by the fluctuations of a $\Z$-valued random walk. Since the latter one is always $\sqrt{n}$ fluctuating whatever $\beta$ is, this shows in a soft way the delocalisation at all $\beta$ of the 
$\alpha>3$-discrete Gaussian chain.

To end the proof of Proposition \ref{pr.comp1}, we now discuss the case $d=1$, $\alpha=3$. The upper bound is obtained exactly as the previous $\log \log(n)$ upper bound by relying on another estimate from \cite{caputo2009recurrence}, namely in appendix B.2 of this paper. The lower bound (at high enough temperature) is more involved but can be obtained following the same strategy as in our main theorem (namely via a Caffarelli-Silverstre extension coupled with a Fröhlich-Spencer Coulomb gas analysis) except we do not optimise the proof here in order to control the invariance principle nor the effective temperature.  
Note that we do not obtain matching upper and lower bounds here. This is due to the fact that, only in this specific case, matching and lower bounds are not obtained for the associated effective resistances in \cite{caputo2009recurrence}. 
See Open Problem \ref{OP6}. 
\qed

\begin{remark}\label{r.NonM}
Interestingly,  using the (easier) first  inequality of Lemma \ref{l.varo} together with Theorem \ref{th.regev}, it is not difficult to show that at least when $\beta$ is large, then the diffusivity constant needs to be smaller than 
\begin{align*}\label{}
\exp( - \frac {c(\alpha)} \beta) \ll \frac 1 \beta\,.
\end{align*}
\textbf{This shows that in the regime of low temperatures, it is impossible to expect $\beta_{eff}(\beta) = \beta$}.

Since we expect $\beta_{eff}(\beta)=\beta$ at all $\beta>0$ in the regime $\alpha\in(2,3)$, we find this ``non-monotony'' behaviour of $\beta_{eff}$ rather intriguing.  
\end{remark}

\section{The hierarchical integer-valued GFF always delocalises}\label{s.9}

The purpose of this short section is to prove that the hierarchical integer-valued GFF, as studied for example in the recent work \cite{biskup2023limit} by Biskup and Huang, always delocalises. This is in contrast with the case of the non-hierarchical long-range integer-valued GFF (at $\alpha=\alpha_c=2$) as shown by Fröhlich and Zegarlinsky in \cite{frohlich1991phase}.

\begin{proposition}\label{pr.HGFF}
The hierarchical integer-valued Gaussian Free Field delocalises at all inverse temperatures $\beta>0$. 
\end{proposition}
We refer to \cite{biskup2023limit} for a precise definition of this model. In words, the hierarchical Gaussian Free Field is the Gaussian free field on the infinite binary tree $T_2$, rooted at the origin. The integer-valued hierarchical Gaussian free field (in \cite{biskup2023limit}) is defined for any $n\geq 1$ as the  hierarchical Gaussian Free Field conditioned to take integer values on the $2^n$ vertices at distance $n$ from the root. Let us call this (conditioned) field $\Psi_n$. 

\ni
{\em Proof of Proposition \ref{pr.HGFF}.}

Let $\hat \Psi_n$ be the hierarchical GFF (on the binary tree $T_2$) conditioned to take integer values on all the leaves at distance less or equal than $n$ from the root $0$. As such, $\hat \Psi_n$ is conditioned to take integer values on more vertices than the field $\Psi_n$. 

We may then argue as in the proof of the Ginibre-type inequality, Proposition \ref{pr.ginibre} (and by following the interpolation proof scheme from \cite{RonFS}) to conclude that for any point $x$ at distance $n$ from the root, one has the inequality 
\begin{align}\label{e.fin}
\Var{\Psi_n(x)} \geq \Var{\hat \Psi_n(x)}\,.
\end{align}

Now, while the law of $\Psi_n$ is not easy to handle (see \cite{biskup2023limit}), we notice that the law of $\hat \Psi_n$ is straightforward: indeed, due to the absence of cycles, it is easy to see that when restricted to vertices at distance less than $n$ from the root, the field $\hat \Psi_n$ is precisely a hierarchical ``random-walk'' free field where the increments up to level-$n$-leaves are Gaussian random variables $\calN(0,\tfrac 1 \beta)$ conditioned to take integer values (which we denoted by $\calN_\Z(0, \tfrac 1 \beta)$). 

In particular, we obtain a constant $c>0$ such that for large $\beta$ (the low temperature regime) and for any point $x$ at distance $n$ from $0$, we have 
\begin{align*}\label{}
\Var{\hat \Psi_n(x)} \geq \exp(-c \beta)*n\,.
\end{align*}
Using the comparison~\eqref{e.fin} this ends the proof of delocalisation.

\qed

\section{Concluding remarks and Open problems}\label{s.OP}

\begin{remark}\label{r.BC}
Other boundary conditions may also be considered. The most natural one we did not discuss would be the {\em periodic boundary conditions}. This amounts to considering the discrete Gaussian chain $\Psi$ on a $1D$ torus $\T_n=\Z/n\Z$ by rooting the field $\Psi$ say at 0. 
The techniques of this paper also work in this case, but one would still need an invariance principle in the Gaussian case to compare with, i.e. a periodic version of Proposition \ref{pr.sheffield}.  


Similarly, our analysis also works in the case of free boundary conditions for the discrete Gaussian chain on $\Z$ (rooted for example at the origin). The infinite volume limit of this chain is briefly discussed in Section \ref{ss.IVL}. Yet, to obtain an invariance principle, we would need the analog of Proposition \ref{pr.sheffield} for free boundary conditions on $\R$. 
%
Let us stress in this case that the very precise asymptotics  on the  {\em harmonic potential} of the associated $\alpha$-long-range random walks on $\Z$ obtained in  \cite{chiarini2023fractional} may be of great help. 
\end{remark}

\begin{OP}\label{OP1}
What happens for the long-range integer-valued GFF on $\Z^2$ when $\alpha_c=4$ ? 
Gaussian domination indicates that fluctuations are at most in $\sqrt{\log \log n}$. It would be tempting to analyse such a field by using a Caffarelli-Silvestre extension in $\Z^3$ with suitable conductances. Most of the proof goes through except one essential feature: the  ``dipoles'' formed in  $\Z^2\times \{0\} \subset \Z^3$  no longer accumulate enough energy cost to overcome the entropy terms precisely at $\alpha_c=4$. 

We find this situation quite intriguing and do not know what to conjecture here, namely delocalisation or not at high temperature. (localisation at low temperature is straightforward). 
\end{OP}

\begin{OP}\label{OP2}
Show an invariance principle (at high temperature) for the long-range IV-GFF in $2d$ when $\alpha>4$ towards a continuum GFF. 

The analogous case for the nearest-neighbour case has been established recently in the breakthrough works  \cite{bauerschmidt2022discrete,bauerschmidt2022discreteN2}.
\end{OP}

\begin{OP}\label{OP3}
We conjecture that our invariance principle holds at all $\beta>0$, (including in the low temperature regime) when 
\begin{align*}\label{}
\alpha \in (2,3)\,. 
\end{align*}
And that furthermore the effective inverse temperature $\beta \mapsto \beta{eff}(\beta)$ is linear on the entire $\R_+$. 
\end{OP}

\begin{remark}\label{}
It is known that it is not the case at $\alpha_c=2$ thanks to \cite{frohlich1991phase}. 
\end{remark}

\begin{OP}\label{op.disc}
Show that at $\alpha=\alpha_c=2$, there is a discontinuity in the delocalisation of the discrete Gaussian Chain. Slurink and Hilhorst predicted in \cite{slurink1983roughening} that the integers should be invisible for all $T>T_c$ while the critical case $T=T_c$ should sit in between. 

We refer to \cite{lammers2022dichotomy} for a a proof of such a result in the $2D$ case using RSW techniques.
\end{OP}


An easier Open Problem than the above one would be  
\begin{OP}\label{OP5}
Prove the delocalisation of the discrete Gaussian chain at all temperatures when $\alpha\in(2,3)$, possibly with quantitative estimates on the variance.  A heuristics is discussed in \cite{frohlich1991phase}, but as far as we know, a proof is still missing. 
\end{OP}

({\em Update:} Delocalisation at arbitrary low temperatures has now been proved in this regime $\alpha\in(2,3)$ in the recent preprint \cite{loren}. Their proof relies on nice relative entropy techniques and builds in particular on \cite{frohlich1981absence,coquille2018absence}. It leads to non-quantitative bounds. As such, obtaining quantitative bounds and possibly an invariance principle towards a fractional Brownian motion is still an interesting open problem). 
%

\begin{OP}\label{OP6}
What happens (in terms of invariance principle / effective temperature / low temperature phase) in the boundary case $d=1$ and $\alpha=3$ ? Proposition \ref{pr.comp1} provides some partial information and shows that some logarithmic corrections do appear. 
\end{OP}

\begin{OP}\label{OP7}
When $\alpha>4$, show delocalisation at high temperature for the long-range SOS model on $\Z^2$, i.e. the interface model with Hamiltonian 
\begin{align*}\label{}
H(\psi) = \sum_{x\neq y} \frac {|\psi_x - \psi_y|}{\|x-y\|_2^\alpha},\,
\end{align*}
 
\end{OP}

\bibliographystyle{alpha}
\bibliography{biblioIVGFF}

\begin{thebibliography}{CvELNR18}

\bibitem[AABK16]{amir2016one}
Gideon Amir, Omer Angel, Itai Benjamini, and Gady Kozma.
\newblock One-dimensional long-range diffusion-limited aggregation i.
\newblock 2016.

\bibitem[ABJL23]{aidekon2023multiplicative}
{\'E}lie A{\"\i}d{\'e}kon, Nathana{\"e}l Berestycki, Antoine Jego, and Titus
  Lupu.
\newblock Multiplicative chaos of the brownian loop soup.
\newblock {\em Proceedings of the London Mathematical Society},
  126(4):1254--1393, 2023.

\bibitem[ACCN88]{aizenman1988discontinuity}
Michael Aizenman, Jennifer~T Chayes, Lincoln Chayes, and Charles~M Newman.
\newblock Discontinuity of the magnetization in one-dimensional $1/|x-y|^2$
  ising and potts models.
\newblock {\em Journal of Statistical Physics}, 50:1--40, 1988.

\bibitem[AGS22]{aru2022percolation}
Juhan Aru, Christophe Garban, and Avelio Sep{\'u}lveda.
\newblock Percolation for 2d classical heisenberg model and exit sets of vector
  valued gff.
\newblock {\em arXiv preprint arXiv:2212.06767}, 2022.

\bibitem[AHPS21]{aizenman2021depinning}
Michael Aizenman, Matan Harel, Ron Peled, and Jacob Shapiro.
\newblock Depinning in integer-restricted gaussian fields and bkt phases of
  two-component spin models.
\newblock {\em arXiv preprint arXiv:2110.09498}, 2021.

\bibitem[Ale11]{alexander2011excursions}
Kenneth Alexander.
\newblock Excursions and local limit theorems for bessel-like random walks.
\newblock 2011.

\bibitem[AYH70]{anderson1970exact}
Philip~W Anderson, G~Yuval, and DR~Hamann.
\newblock Exact results in the kondo problem. ii. scaling theory, qualitatively
  correct solution, and some new results on one-dimensional classical
  statistical models.
\newblock {\em Physical Review B}, 1(11):4464, 1970.

\bibitem[B{\"a}u23]{baumler2023recurrence}
Johannes B{\"a}umler.
\newblock Recurrence and transience of symmetric random walks with long-range
  jumps.
\newblock {\em Electronic Journal of Probability}, 28:1--24, 2023.

\bibitem[BGR61]{blumenthal1961distribution}
Robert~M Blumenthal, Ronald~K Getoor, and DB~Ray.
\newblock On the distribution of first hits for the symmetric stable processes.
\newblock {\em Transactions of the American Mathematical Society},
  99(3):540--554, 1961.

\bibitem[BH23]{biskup2023limit}
Marek Biskup and Haiyu Huang.
\newblock A limit law for the maximum of subcritical {DG}-model on a
  hierarchical lattice.
\newblock {\em arXiv preprint arXiv:2309.09389}, 2023.

\bibitem[Bis20]{biskup2020extrema}
Marek Biskup.
\newblock Extrema of the two-dimensional discrete gaussian free field.
\newblock In {\em Random Graphs, Phase Transitions, and the Gaussian Free
  Field: PIMS-CRM Summer School in Probability, Vancouver, Canada, June 5--30,
  2017}, pages 163--407. Springer, 2020.

\bibitem[Bol99]{bolthausen1999large}
Erwin Bolthausen.
\newblock Large deviations and interacting random walks and random surfaces.
\newblock {\em Lecture notes St-Flour}, 1999.

\bibitem[BPR22a]{bauerschmidt2022discrete}
Roland Bauerschmidt, Jiwoon Park, and Pierre-Fran{\c{c}}ois Rodriguez.
\newblock The discrete gaussian model, {I.} renormalisation group flow at high
  temperature.
\newblock {\em arXiv preprint arXiv:2202.02286}, 2022.

\bibitem[BPR22b]{bauerschmidt2022discreteN2}
Roland Bauerschmidt, Jiwoon Park, and Pierre-Fran{\c{c}}ois Rodriguez.
\newblock The discrete gaussian model, {II.} infinite-volume scaling limit at
  high temperature.
\newblock {\em arXiv preprint arXiv:2202.02287}, 2022.

\bibitem[Car81]{cardy1981one}
John~L Cardy.
\newblock One-dimensional models with 1/r2 interactions.
\newblock {\em Journal of Physics A: Mathematical and General}, 14(6):1407,
  1981.

\bibitem[CCH16]{chiarini2016extremes}
Alberto Chiarini, Alessandra Cipriani, and Rajat~Subhra Hazra.
\newblock Extremes of some gaussian random interfaces.
\newblock {\em Journal of Statistical Physics}, 165:521--544, 2016.

\bibitem[CD00]{caputo2000large}
P~Caputo and J-D Deuschel.
\newblock Large deviations and {$\P$} variational principle for harmonic
  crystals.
\newblock {\em Communications in Mathematical Physics}, 209:595--632, 2000.

\bibitem[CD01]{caputo2001critical}
P~Caputo and J-D Deuschel.
\newblock Critical large deviations in harmonic crystals with long-range
  interactions.
\newblock {\em The Annals of Probability}, 29(1):242--287, 2001.

\bibitem[CFG09]{caputo2009recurrence}
Pietro Caputo, Alessandra Faggionato, and Alexandre Gaudilliere.
\newblock Recurrence and transience for long-range reversible random walks on a
  random point process.
\newblock 2009.

\bibitem[CJR23]{chiarini2023fractional}
Leandro Chiarini, Milton Jara, and Wioletta~M Ruszel.
\newblock Fractional edgeworth expansions for one-dimensional heavy-tailed
  random variables and applications.
\newblock {\em Electronic Journal of Probability}, 28:1--42, 2023.

\bibitem[CL83]{caldeira1983quantum}
Amir~O Caldeira and Anthony~J Leggett.
\newblock Quantum tunnelling in a dissipative system.
\newblock {\em Annals of physics}, 149(2):374--456, 1983.

\bibitem[CLRV15]{ciaurri2015connection}
{\'O}scar Ciaurri, Carlos Lizama, Luz Roncal, and Juan~Luis Varona.
\newblock On a connection between the discrete fractional laplacian and
  superdiffusion.
\newblock {\em Applied Mathematics Letters}, 49:119--125, 2015.

\bibitem[CS07]{caffarelli2007extension}
Luis Caffarelli and Luis Silvestre.
\newblock An extension problem related to the fractional laplacian.
\newblock {\em Communications in partial differential equations},
  32(8):1245--1260, 2007.

\bibitem[CSS03]{caux2003two}
J-S Caux, H~Saleur, and F~Siano.
\newblock The two-boundary sine-gordon model.
\newblock {\em Nuclear Physics B}, 672(3):411--461, 2003.

\bibitem[CSZ16]{caravenna2016continuum}
Francesco Caravenna, Rongfeng Sun, and Nikos Zygouras.
\newblock The continuum disordered pinning model.
\newblock {\em Probability theory and related fields}, 164:17--59, 2016.

\bibitem[CvELNR18]{coquille2018absence}
Loren Coquille, Aernout~CD van Enter, Arnaud Le~Ny, and Wioletta~M Ruszel.
\newblock Absence of dobrushin states for 2 d long-range ising models.
\newblock {\em Journal of Statistical Physics}, 172:1210--1222, 2018.

\bibitem[CvENR24]{loren}
Loren Coquille, Aernout~CD van Enter, Arnaud~Le Ny, and Wioletta~M Ruszel.
\newblock Absence of shift-invariant {G}ibbs states (delocalisation) for
  one-dimensional $\mathbb{Z}$-valued fields with {L}ong-{R}ange interactions.
\newblock {\em arXiv preprint arXiv:2401.17722}, 2024.

\bibitem[DCGT20]{duminil2020long}
Hugo Duminil-Copin, Christophe Garban, and Vincent Tassion.
\newblock Long-range models in 1d revisited.
\newblock {\em arXiv preprint arXiv:2011.04642}, 2020.

\bibitem[DNPV12]{di2012hitchhiker}
Eleonora Di~Nezza, Giampiero Palatucci, and Enrico Valdinoci.
\newblock Hitchhiker?s guide to the fractional {S}obolev spaces.
\newblock {\em Bulletin des sciences math{\'e}matiques}, 136(5):521--573, 2012.

\bibitem[Don97]{doney1997one}
Ronald~A Doney.
\newblock One-sided local large deviation and renewal theorems in the case of
  infinite mean.
\newblock {\em Probability theory and related fields}, 107:451--465, 1997.

\bibitem[Dys69]{dyson1969non}
Freeman~J Dyson.
\newblock Non-existence of spontaneous magnetization in a one-dimensional ising
  ferromagnet.
\newblock 1969.

\bibitem[Fey55]{feynman1955slow}
Richard~Phillips Feynman.
\newblock Slow electrons in a polar crystal.
\newblock {\em Physical Review}, 97(3):660, 1955.

\bibitem[FP81]{frohlich1981absence}
J{\"u}rg Fr{\"o}hlich and Charles Pfister.
\newblock On the absence of spontaneous symmetry breaking and of crystalline
  ordering in two-dimensional systems.
\newblock {\em Communications in Mathematical Physics}, 81:277--298, 1981.

\bibitem[FS81]{FS}
J{\"u}rg Fr{\"o}hlich and Thomas Spencer.
\newblock The {Kosterlitz-Thouless transition in two-dimensional abelian spin
  systems and the Coulomb gas}.
\newblock {\em Communications in Mathematical Physics}, 81(4):527--602, 1981.

\bibitem[FS82]{FS82}
J{\"u}rg Fr{\"o}hlich and Thomas Spencer.
\newblock The phase transition in the one-dimensional ising model with $1/r^2$
  interaction energy.
\newblock 1982.

\bibitem[FZ85]{fisher1985quantum}
Matthew~PA Fisher and Wilhelm Zwerger.
\newblock Quantum brownian motion in a periodic potential.
\newblock {\em Physical Review B}, 32(10):6190, 1985.

\bibitem[FZ91]{frohlich1991phase}
J~Fr{\"o}hlich and B~Zegarlinski.
\newblock The phase transition in the discrete gaussian chain with $1/r^2$
  interaction energy.
\newblock {\em Journal of statistical physics}, 63:455--485, 1991.

\bibitem[Gin70]{Ginibre}
Jean Ginibre.
\newblock General formulation of griffiths' inequalities.
\newblock {\em Communications in mathematical physics}, 16(4):310--328, 1970.

\bibitem[GL62]{garsia1962discrete}
Adriano Garsia and John Lamperti.
\newblock A discrete renewal theorem with infinite mean.
\newblock {\em Commentarii Mathematici Helvetici}, 37(1):221--234, 1962.

\bibitem[GS23a]{GS2}
Christophe Garban and Avelio Sep{\'u}lveda.
\newblock Quantitative bounds on vortex fluctuations in $2d$ coulomb gas and
  maximum of the integer-valued {G}aussian free field.
\newblock {\em Proceedings of the London Mathematical Society},
  127(3):653--708, 2023.

\bibitem[GS23b]{GS1}
Christophe Garban and Avelio Sep{\'u}lveda.
\newblock Statistical reconstruction of the {GFF} and {KT} transition.
\newblock {\em Journal of the European Mathematical Society}, 26(2):639--694,
  2023.

\bibitem[HS13]{hammond2013power}
Alan Hammond and Scott Sheffield.
\newblock Power law pólya?s urn and fractional {B}rownian motion.
\newblock {\em Probability Theory and Related Fields}, 157(3-4):691--719, 2013.

\bibitem[KH82]{kjaer1982discrete}
KH~Kjaer and HJ~Hilhorst.
\newblock The discrete gaussian chain with $1/r^n$ interactions: Exact results.
\newblock {\em Journal of Statistical Physics}, 28:621--632, 1982.

\bibitem[Kol40]{kolmogorov1940wienersche}
Andrei~N Kolmogorov.
\newblock Wienersche spiralen und einige andere interessante kurven in
  hilbertscen raum, cr (doklady).
\newblock {\em Acad. Sci. URSS (NS)}, 26:115--118, 1940.

\bibitem[KP17]{RonFS}
Vital Kharash and Ron Peled.
\newblock The {F}röhlich-{S}pencer proof of the
  {B}erezinskii-{K}osterlitz-{T}houless transition.
\newblock {\em arXiv preprint arXiv:1711.04720}, 2017.

\bibitem[Lam22a]{lammers2022dichotomy}
Piet Lammers.
\newblock A dichotomy theory for height functions.
\newblock {\em arXiv preprint arXiv:2211.14365}, 2022.

\bibitem[Lam22b]{lammers2022height}
Piet Lammers.
\newblock Height function delocalisation on cubic planar graphs.
\newblock {\em Probability Theory and Related Fields}, 182(1-2):531--550, 2022.

\bibitem[Law08]{lawler2008conformally}
Gregory~F Lawler.
\newblock {\em Conformally invariant processes in the plane}.
\newblock Number 114. American Mathematical Soc., 2008.

\bibitem[LL10]{lawler2010random}
Gregory~F Lawler and Vlada Limic.
\newblock {\em Random walk: a modern introduction}, volume 123.
\newblock Cambridge University Press, 2010.

\bibitem[LO23]{lammers2023delocalisation}
Piet Lammers and S{\'e}bastien Ott.
\newblock Delocalisation and absolute-value-fkg in the solid-on-solid model.
\newblock {\em Probability Theory and Related Fields}, pages 1--25, 2023.

\bibitem[LSSW16]{lodhia2016fractional}
Asad Lodhia, Scott Sheffield, Xin Sun, and Samuel~S Watson.
\newblock Fractional gaussian fields: A survey.
\newblock {\em Probability Surveys}, 13:1--56, 2016.

\bibitem[LT20]{lammers2020macroscopic}
Piet Lammers and Martin Tassy.
\newblock Macroscopic behavior of lipschitz random surfaces.
\newblock {\em arXiv preprint arXiv:2004.15025}, 2020.

\bibitem[MVN68]{mandelbrot1968fractional}
Benoit~B Mandelbrot and John~W Van~Ness.
\newblock Fractional brownian motions, fractional noises and applications.
\newblock {\em SIAM review}, 10(4):422--437, 1968.

\bibitem[Par22]{park2022central}
Jiwoon Park.
\newblock Central limit theorem for multi-point functions of the 2d discrete
  gaussian model at high temperature.
\newblock {\em arXiv preprint arXiv:2211.14367}, 2022.

\bibitem[PSC00]{pittet2000stability}
Christophe Pittet and Laurent Saloff-Coste.
\newblock On the stability of the behavior of random walks on groups.
\newblock {\em Journal of Geometric Analysis}, 10(4), 2000.

\bibitem[RSD17]{regev2017inequality}
Oded Regev and Noah Stephens-Davidowitz.
\newblock An inequality for gaussians on lattices.
\newblock {\em SIAM Journal on Discrete Mathematics}, 31(2):749--757, 2017.

\bibitem[RY13]{revuz2013continuous}
Daniel Revuz and Marc Yor.
\newblock {\em Continuous martingales and Brownian motion}, volume 293.
\newblock Springer Science \& Business Media, 2013.

\bibitem[SD85]{spohn11985quantum}
Herbert Spohn and Rolf D{\"u}mcke.
\newblock Quantum tunneling with dissipation and the {I}sing model over
  $\mathbb{R}$.
\newblock {\em Journal of statistical physics}, 41:389--423, 1985.

\bibitem[SH83]{slurink1983roughening}
J~Slurink and HJ~Hilhorst.
\newblock Roughening behavior of a one-dimensional crystal surface with inverse
  square potential.
\newblock {\em Physica A: Statistical Mechanics and its Applications},
  120(3):627--634, 1983.

\bibitem[She05]{sheffield2005random}
Scott Sheffield.
\newblock Random surfaces.
\newblock {\em Ast{\'e}risque}, 304, 2005.

\bibitem[Sla18]{slade2018critical}
Gordon Slade.
\newblock Critical exponents for long-range o (n) o (n) models below the upper
  critical dimension.
\newblock {\em Communications in Mathematical Physics}, 358:343--436, 2018.

\bibitem[Spi58]{spitzer1958some}
Frank Spitzer.
\newblock Some theorems concerning 2-dimensional brownian motion.
\newblock {\em Transactions of the American Mathematical Society},
  87(1):187--197, 1958.

\bibitem[Spi13]{spitzer2013principles}
Frank Spitzer.
\newblock {\em Principles of random walk}, volume~34.
\newblock Springer Science \& Business Media, 2013.

\bibitem[Spo86a]{spohn1986roughening}
Herbert Spohn.
\newblock Roughening and pinning transitions for the polaron.
\newblock {\em Journal of Physics A: Mathematical and General}, 19(4):533,
  1986.

\bibitem[Spo86b]{spohn1986statistical}
Herbert Spohn.
\newblock Statistical mechanics of the effective mass of the polaron.
\newblock {\em Physical Review B}, 33(12):8906, 1986.

\bibitem[Spo87]{spohn1987effective}
Herbert Spohn.
\newblock Effective mass of the polaron: A functional integral approach.
\newblock {\em Annals of Physics}, 175(2):278--318, 1987.

\bibitem[Spo05]{spohn2005models}
Herbert Spohn.
\newblock Models of statistical mechanics in one dimension originating from
  quantum ground states.
\newblock In {\em Statistical Mechanics and Field Theory: Mathematical Aspects:
  Proceedings of the International Conference on the Mathematical Aspects of
  Statistical Mechanics and Field Theory Held in Groningen, The Netherlands,
  August 26--30, 1985}, pages 209--233. Springer, 2005.

\bibitem[SS95]{leclair1995boundary}
H~Saleur and S~Skorik.
\newblock Boundary energy and boundary states in integrable quantum field
  theories.
\newblock {\em Nuclear Physics B}, 453(3):581--618, 1995.

\bibitem[SSW95]{saleur1995boundary}
H~Saleur, S~Skorik, and NP~Warner.
\newblock The boundary sine-gordon theory: classical and semi-classical
  analysis.
\newblock {\em Nuclear Physics B}, 441(3):421--436, 1995.

\bibitem[Tho69]{thouless1969long}
DJ~Thouless.
\newblock Long-range order in one-dimensional ising systems.
\newblock {\em Physical Review}, 187(2):732, 1969.

\bibitem[Var86]{varopoulos1986theorie}
Nicholas~Th Varopoulos.
\newblock Th{\'e}orie du potentiel sur des groupes et des vari{\'e}t{\'e}s.
\newblock {\em CR Acad. Sci. Paris S{\'e}r. I Math}, 302(6):203--205, 1986.

\bibitem[Vel06]{velenik2006localization}
Yvan Velenik.
\newblock Localization and delocalization of random interfaces.
\newblock 2006.

\bibitem[Whi02]{whitt2002stochastic}
Ward Whitt.
\newblock Stochastic-process limits: an introduction to stochastic-process
  limits and their application to queues.
\newblock {\em Space}, 500:391--426, 2002.

\bibitem[Wir19]{wirth2019maximum}
Mateo Wirth.
\newblock Maximum of the integer-valued gaussian free field.
\newblock {\em arXiv preprint arXiv:1907.08868}, 2019.

\end{thebibliography}

\appendix

\section{Caffarelli-Silvestre extension of the discrete Laplacian in $\Z^d$ and Bessel random walks on $\Z^{d+1}$}\label{s.BesselZ}

The extension we used in Section \ref{pr.Bessel} was based on a random walk on the diamond graph $\calD$ (see Figure \ref{f.diamond}). Though it was sufficient for our application, the diamond graph $\calD$ is very much specific to the two-dimensional case (where by two-dimensional we mean here $2=1+1$) and it does not easily extend to higher dimensions. The purpose of this appendix is to introduce a ``Caffarelli-Silvestre'' extension of the discrete fractional Laplacians from $\Z^d$ to  $\Z^{d+1}$. 
\medskip

As we shall see below, the estimates are very similar as in Section \ref{s.Bessel} except some further computations are needed due to the combinatorial terms arising from choosing which of the vertical/horizontal directions are moving.  
\medskip

We start with the  $d=1$ case in Subsection \ref{aa.2d} below, where we introduce suitable Bessel-type walks on $\Z^{1+1}$.
From this analysis, the extension to any $\Z^{d+1}$ is not difficult as we shortly explain in Subsection \ref{aa.dim}. 




\subsection{Bessel random walks on $\Z^2$.}\label{aa.2d}

Similarly as in Section \ref{s.Bessel}, 
for each $s\geq 0$, we build a Markov processes $(Z_n^{(s)})_{n\geq 0}$ on the grid $\Z^2$ as follows.

 When $s=0$, we will consider the simple random walk on $\Z^2$ (or equivalently the reflected simple random walk on $\Z\times \N$). When $s>0$ we introduce the same type of confining force as on $\calD$. Namely, at each time $n\geq 0$, the Markov process  $(Z_n^{(s)})_{n\geq 0}$ will:
\bi
\item either move vertically with probability $1/2$ (up or down) according to a {\em discrete Bessel random walk} as studied for example in \cite{alexander2011excursions}. 
\item or move horizontally with probability $1/2$ according to a symmetric simple random walk. 
\ei

The vertical moves follow the exact same discrete Bessel process as in Section \ref{s.Bessel}, i.e. the Markov transition kernel $Q_s: \N \times \N \to [0,1]$. 

\bi
\item $Q_s(0,1)=1$.
\item $Q_s(r,r-1)+Q(r,r+1)=1$ whenever $r\geq 1$.
\item $Q_s(r,r+1)= (\tfrac 1 2  - \frac s 4\, \frac 1 r ) \vee \tfrac 1 4$ whenever $r\geq 1$.
\ei

\ni
The main result of this appendix is the analogue of Proposition \ref{pr.Bessel} on the grid  $\Z^2$:
\begin{proposition}\label{pr.BesselZ}
For any $s>-1$, there exists a constant $c= c(s)>0$ such that as $k\to \infty$ 
\begin{align*}\label{}
\FK{}{(0,0)}{Z^{(s)}_{\tau_{\Z}} = (0,k)} \sim \frac {c(s)} {k^{2+s}} \,.
\end{align*}
\end{proposition}

\ni
{\em Proof.} 

The main input will again be Theorem \ref{th.ALEX} from \cite{alexander2011excursions}, and we will use the same notations: 
\bnum
\item Call the function $g_s(n):= \FK{Q_s}{0}{\tau_0 = n}$.
\item Call the heat-kernel of the simple random walk $p_{\Z}(t,x,y)=p_\Z(t,x-y)$. 
\enum

Let us analyse the case where $k=2x > 0 $ is even. (The odd case is obtained in the same fashion with same asymptotics). 
For any such $x$, we have 
\begin{align*}\label{}
\FK{}{(0,0)}{Z^{(s)}_{\tau_{\Z}} = (0,2x)} & = \sum_{n = x}^\infty  \sum_{m=x}^n \binom{2n}{2m}  2^{-2n} p_\Z(2m,2x) g_s(2(n-m)) \\
& = \sum_{n = x}^\infty  \sum_{m=0}^n \binom{2n}{2m}  2^{-2n} p_\Z(2m,2x) g_s(2(n-m))
\end{align*}
The first equality is because at least $2x$ steps are required to reach $2x$. 
For the second, we just use the fact the heat kernel is zero if there is not enough time to reach $2x$ in time $2m$. 

Let us first change of variable to center the sum $\sum_{m=0}^n$ around  $n/2$. We will use the following slight abuse of notations by writing 
\begin{align*}\label{}
\sum_{m=0}^n \binom{2n}{2m}  2^{-2n} \text{  as  }
\sum_{l=-n/2}^{n/2} \binom{2n}{n+2 l}  2^{-2n}
\end{align*}
whether $n$ is odd or not. (For example, if $n=3$, then we are summing over the set $\{-3/2,-1/2,1/2,3/2\}$). 

We thus have 
\begin{align*}\label{}
\FK{}{(0,0)}{Z^{(s)}_{\tau_{\Z}} = (0,2x)} 
& = \sum_{n = x}^\infty  \sum_{l=-n/2}^{n/2} \binom{2n}{n+2 l}  2^{-2n} p_\Z(n+2l ,2x) g_s(n- 2l)
\end{align*}

Let $\eps>0$ be any fixed small parameter. 
Our first observation is that the contributions of integers $n\leq x^{2-\eps}$ is negligible. Indeed, 
\begin{align*}\label{}
\sum_{n=x}^{n^{2-\eps}} \sum_{l=-n/2}^{n/2}  \binom{2n}{n+2 l}  2^{-2n} p_\Z(n+2l ,2x) g_s(n- 2l)
& \leq \sum_{n=x}^{n^{2-\eps}}
 \sum_{l=-n/2}^{n/2}  \binom{2n}{n+2 l}  2^{-2n} p_\Z(n+2l ,2x)\,,
\end{align*}
and each heat kernel in this sum is upper bounded (for small enough $c>0$) by 
\begin{align*}\label{}
\exp(- c (2x)^2/ (2x)^{2-\eps} ) \leq \exp(- c k^\eps)\,.
\end{align*}
This readily implies that the $\sum_{n=x}^{x^{2-\eps}}$ is negligible w.r.t. to the power law in $k$ asymptotics that  we are looking after.

Furthermore, exactly as we argued above for the terms $n\leq x^{2-\eps}$, one can easily show that the $l$ satisfying $|l|\geq n^{1/2+\eps}$ also lead to a stretch exponential contribution. We thus get 
\begin{align}\label{e.GS}
\FK{}{(0,0)}{Z^{(s)}_{\tau_{\Z}} = (0,2x)} 
& =  \sum_{n=x^{2-\eps}}^\infty \sum_{l= -n^{1/2+\eps}}^{n^{1/2+\eps}}  \binom{2n}{n+2 l}  2^{-2n}
 \\
& \hskip 3 cm \, \times p_\Z(n+2l ,2x) g_s(n-2l) + \text{ Strech. Exp.}\,,  \nn
\end{align}
where the sum $-n^{1/2+\eps} \leq l \leq n^{1/2+\eps}$ is understood as a sum over integers if $n$ is even and instead as a sum over $\Z+\tfrac 1 2 $.

Now, again by the quantitative Local CLT theorem (\cite{lawler2010random}),  there exists a constant $C>0$ s.t. uniformly in $n$ and in  $l\in \{-n/2,\ldots, n/2\}$,  
\begin{align*}\label{}
\left| p_\Z(n+2l ,2x)  - \frac{1}{\sqrt{2\pi (n+2l)}} e^{-\frac{(2x)^2}{2*(n+2l)}}\right| \leq \frac C {n^{3/2}} \,.
\end{align*}
Using that $\frac l n = O(n^{-1/2+\eps})$, we obtain 
\begin{align*}\label{}
p_\Z(n+2l ,2x) & =  \frac 1 {\sqrt{2\pi n}} (1+ O(\frac l n)) e^{-\frac{2 x^2}{n} (1+ O(\frac l n))} + O(n^{-3/2})\,.
\end{align*}
Let us deal with the upper bound (a matching order lower bound follows from the same analysis in order the conclude the proof of the asymptotics in Proposition \ref{pr.BesselZ}). We get 
\begin{align*}\label{}
p_\Z(n+2l ,2x) & \leq   \frac 1 {\sqrt{2\pi n}} e^{-\frac{2 x^2}{n} (1 -  c n^{-1/2+\eps})} + O(n^{-1+\eps})\,,
\end{align*}
(where we assume we have chosen $\eps\leq 1/2$ here). 

For the Bessel term $g_s(n-2l)$, using Theorem \ref{th.ALEX}, we have that for any small $\delta>0$, then for $n$ large enough and uniformly in $-n^{1/2+\eps} \leq l \leq  n^{1/2+\eps}$, 
\begin{align*}\label{}
g_s(n-2l) \sim c(s)\, \frac 1{(n-2l)^{\frac {3+s} 2}}
& \leq (1+\delta)  c(s)  n^{-\frac {3+s} 2} 
\end{align*}
Going back to~\eqref{e.GS}, we find 
\begin{align*}\label{}
& \FK{}{(0,0)}{Z^{(s)}_{\tau_{\Z}} = (0,2x)}  \\
& \leq   \sum_{n=x^{2-\eps}}^\infty  \sum_{l=-n^{1/2+\eps}}^{n^{1/2+\eps}}\binom{2n}{n+2l}2^{-2n}   \,\,    \left( \frac 1 {\sqrt{2\pi n}} e^{-\frac{2 x^2}{n} (1 -  c n^{-1/2+\eps})}  +  O(n^{-1+\eps})\right)  \\
& \hskip 6cm   \times (1+\delta) c(s)  n^{-\frac {3+s} 2}  \\
& \leq  \sum_{n=x^{2-\eps}}^\infty 
 \left( \frac 1 {\sqrt{2\pi n}} e^{-\frac{2 x^2}{n} (1 -  c n^{-1/2+\eps})}  + O(n^{-1+\eps})\right)   \times  (1+\delta) c(s)  n^{-\frac {3+s} 2}   \\
& 
=  (1+\delta) c(s) \sum_{n=x^{2-\eps}}^\infty
\left(  \frac 1 {\sqrt{2\pi}}  n^{-(2+\frac s 2)} e^{-\frac{2 x^2}{n}  +  c \frac{2 x^2}
{n^{3/2-\eps}}} + O(n^{-(\tfrac 5 2 + \tfrac s 2 -\eps)}) \right) \\
&\leq 
(1+\delta) c(s) \sum_{n=x^{2-\eps}}^\infty
\frac 1 {\sqrt{2\pi}}  n^{-(2+\frac s 2)} e^{-\frac{2 x^2}{n}  +  c \frac{2 x^2}
{n^{3/2-\eps}}}     + O(x^{-(2-\eps)(\tfrac 3 2 + \frac s 2 - \eps)})\,.
\end{align*}

Note that in the error term, the exponent $(2-\eps)(\tfrac 3 2 + \frac s 2 - \eps)$ is equal (for small $\eps$) to $3+s - O(\eps)$ which is stricly larger than the desired exponent $2+s$ in Proposition \ref{pr.BesselZ} (when $\eps$ is chosen to be small enough). 

We are thus left with controlling the asymptotic of the above series. By applying Euler-Maclaurin comparison's formula, we get
\begin{align*}\label{}
&  \sum_{n=x^{2-\eps}}^\infty
\frac 1 {\sqrt{2\pi}}  n^{-(2+\frac s 2)} e^{-\frac{2 x^2}{n}  +  c \frac{2 x^2}
{n^{3/2-\eps}}} \\
& \sim_{x\to \infty}
\int_{n=x^{2-\eps}}^\infty
\frac 1 {\sqrt{2\pi}} dn\,   n^{-(2+\frac s 2)} e^{-\frac{2 x^2}{n}  +  c \frac{2 x^2}
{n^{3/2-\eps}}} \\
& = (n=x^2 u) x^2 \int_{x^{-\eps}}^\infty x^{- 4 -s}   du \frac 1 {u^{2+ \frac s 2}} e^{-\frac 2 u + c \frac 1 {x^{1-\eps} u^{3/2-\eps}}} \\
& \sim \frac 1 {x^{2+s}} \int_0^\infty du \frac{e^{-\frac 2 u}}{u^{2+\tfrac s 2}}
\end{align*}
With the same analysis, we also obtain for any $\delta>0$ and for $n$ large enough, 
\begin{align*}\label{}
\FK{}{(0,0)}{Z^{(s)}_{\tau_{\Z}} = (0,2x)}  \geq (1-\delta) \left( c(s) 
\int_0^\infty du \frac{e^{-\frac 2 u}}{u^{2+\tfrac s 2}} \right) 
\frac 1 {x^{2+s}}\,, 
\end{align*} 
which concludes the proof of Proposition \ref{pr.BesselZ}. 
\qed

\subsection{Bessel walks in $\Z^{d+1}$ for any $d\geq 1$.}\label{aa.dim}
The above proof easily extends to any dimension $d\geq 1$. The only main difference is that one needs to rely instead on the 
quantitative local CLT theorem on $\Z^d$ (instead of $\Z$) which is proved for example in \cite{lawler2010random}. This allows us to obtain the following estimate. 

\begin{proposition}\label{l.keyD}

For any dimension $d\geq 1$ and any $s>-1$,  if $Z^{(s)}$ is the Markov chain on $\Z^{d+1}$ with SRW steps along the first $d$ coordinates and the $s$-discrete Bessel process (defined as in Section \ref{s.Bessel}) along the $d+1^{th}$ coordinate, then there exists a constants $C(s)>0$ such that for any $a,b \in \Z^d \times \{0\}$, as $\|a-b\|_2 \to \infty$, 
\begin{align}\label{}
 \FK{}{a}{Z^{(s)}_{\tau_{\Z}} = b} \sim  \frac {C(s)} {\|a-b\|^{d+1+s}}\,.
\end{align}
\end{proposition}

\end{document}